\documentclass[11pt, a4paper, reqno]{amsart}

\usepackage[T1]{fontenc}
\usepackage[utf8]{inputenc}
\usepackage[english]{babel}
\usepackage{mathpazo}
\usepackage{microtype}
\usepackage{comment}

\usepackage{amsmath, amssymb, amsthm, amsfonts, mathrsfs}
\usepackage{mathtools}

\usepackage{pgfplots}
\pgfplotsset{compat=1.18}
\usepackage{tikz}
\usetikzlibrary{calc}

\usepackage[margin=1.25in]{geometry}
\usepackage{setspace}
\usepackage{xcolor}
\usepackage[colorlinks=true, linkcolor=blue!70!black, citecolor=green!50!black, urlcolor=blue!70!black]{hyperref}

\theoremstyle{plain}
\newtheorem{theorem}{Theorem}[section]
\newtheorem{lemma}[theorem]{Lemma}
\newtheorem{proposition}[theorem]{Proposition}
\newtheorem{corollary}[theorem]{Corollary}
\newtheorem{conjecture}[theorem]{Conjecture}

\theoremstyle{definition}
\newtheorem{definition}[theorem]{Definition}
\newtheorem{remark}[theorem]{Remark}

\numberwithin{equation}{section}

\newcommand{\Sbb}{\mathbb{S}}
\newcommand{\B}{\mathbb{B}}
\newcommand{\R}{\mathbb{R}}
\newcommand{\vol}{\operatorname{vol}}

\title[Maximal Steklov Metrics on Higher-Dimensional Cylinders]{Maximal Steklov Metrics on Higher-Dimensional Cylinders and the Critical \texorpdfstring{$n$}{n}-Catenoid}

\author[A. D. Noutchegueme]{Alain Didier Noutchegueme}
\address[Alain Didier Noutchegueme]{D\'epartement de math\'ematiques et de statistique, Universit\'e de Montr\'eal, CP 6128 succ Centre-Ville, Montr\'eal, QC H3C 3J7, Canada}
\email{\href{mailto:alain.didier.noutchegueme@umontreal.ca}{alain.didier.noutchegueme@umontreal.ca}}

\begin{document}

\begin{abstract}
In this article, we study the maximization of the normalized first nonzero Steklov eigenvalue over rotationally symmetric metrics and boundary weights on the $n$-dimensional cylinder, a problem that reduces to a one-dimensional analysis. We identify the stable regime for the conformal length of the cylinder, providing the criterion that guarantees the existence of maximal metrics within a given conformal class. Moreover, by optimizing globally over the moduli space of cylinders, we prove that there exists a unique global maximizer (up to scaling and diffeomorphism), which is realized by the metric induced on the critical $n$-catenoid, the higher-dimensional analogue of the critical catenoid of Fraser and Schoen. 
\end{abstract}

\maketitle

\section{Introduction}
\label{sec:intro}

The Steklov eigenvalue problem, first introduced by V. Steklov in 1902, describes the vibrational modes of a manifold where the mass is concentrated entirely on the boundary \cite{GP17}. Let $(\Omega, g)$ be a compact Riemannian manifold of dimension $n$ with smooth boundary $\partial \Omega$, and let $\rho > 0$ be a smooth weight function on $\partial \Omega$. The Steklov eigenvalue problem is given by
\begin{equation} \label{eq:steklov_problem}
    \begin{cases}
    \Delta_g u = 0 & \text{in } \Omega, \\
    \partial_\nu u = \sigma \rho u & \text{on } \partial \Omega,
    \end{cases}
\end{equation}
where $\Delta_g$ is the Laplace-Beltrami operator and $\nu$ is the outward unit normal to $\partial \Omega$. The spectrum is discrete, provided the trace operator $H^1(\Omega) \to L^2(\partial \Omega)$ is compact (see, e.g., \cite{AM}), forming a sequence $0 = \sigma_0 < \sigma_1 \le \sigma_2 \le \dots \nearrow \infty$.

A central theme in geometric spectral theory is the maximization of these eigenvalues under various geometric constraints. The classical theory for surfaces, anchored by the foundational works of Fraser and Schoen \cite{FS11, FS12, FS16} and extended in \cite{Kokarev, Pet19}, reveals that the topology controls the spectrum and that extremal metrics are closely tied to free boundary minimal surfaces. When a smooth global maximizer of $\sigma_1(\Omega, g) L_g(\partial \Omega)$ exists over all Riemannian metrics on a given smooth surface (such as for genus $0$ with an arbitrary number of boundary components \cite{KKMS}, and more generally via \cite{Pet24}), it is realized, up to a conformal factor equal to $1$ on $\partial \Omega$, as the induced metric of a free boundary minimal immersion in a unit ball $\B^m$. A crucial ingredient of this two-dimensional theory is the conformal invariance of the Dirichlet energy.

For dimensions $n \ge 3$, however, one encounters existence issues. The work of Colbois, El Soufi, and Girouard \cite{CEG}, alongside Girouard and Polymerakis \cite{GP2}, demonstrated that on any fixed topology, geometric degeneration can make eigenvalues arbitrarily large while keeping the boundary volume fixed, showing that spectral control is impossible without further geometric constraints. To address this instability, Karpukhin and M\'etras \cite{KM} introduced the scale-invariant normalization:
\begin{equation} \label{eq:KM_normalization}
    \bar{\sigma}_k(\Omega, g, \rho) = \sigma_k(\Omega, g, \rho) \left( \int_{\partial \Omega} \rho \, d\mathcal{H}_g \right) \vol(\Omega, g)^{\frac{2-n}{n}}.
\end{equation}
This establishes the natural normalization in higher dimensions. It is the only one involving both the interior volume and the boundary area, for which one can formally hope for a smooth critical metric. For this functional, critical metrics in a conformal class correspond to free boundary $n$-harmonic maps $\Phi\colon(\Omega,\partial \Omega)\to\mathbb{B}^{m+1}$, satisfying the Euler--Lagrange equation:
\[
\left\{
\begin{array}{rcll}
\Delta_{n, g} \Phi & = & 0 & \text{ in } \Omega, \\
\partial_{\nu} \Phi & = & \sigma_k \rho_\Phi \Phi & \text{ on } \partial \Omega,
\end{array}
\right.
\]
where $\Delta_{n, g} \Phi = -\mathrm{div}_g(|d\Phi|_g^{n-2} d\Phi)$ is the geometric $n$-Laplacian associated with the metric $g$, and $\rho_\Phi = \sigma_k^{-1}(\partial_{\nu} \Phi \cdot \Phi) > 0$. Regularity theory for the $p$-Laplacian implies that $\Phi$ is of class $C^{1,\alpha}$ in $\Omega$ \cite{Lin}.

Yet, even with this normalization, the unconstrained optimization problem remains unstable. As shown in Proposition \ref{prop:steklov_blowup_variable}, the spectrum continues to degenerate in general. 

To provide a clear perspective of $\bar{\sigma}_1$-optimization, Table \ref{tab:landscape} summarizes various regimes.

\begin{table}[htpb]
\centering
\renewcommand{\arraystretch}{1.5}
\begin{tabular}{p{0.22\linewidth} p{0.34\linewidth} p{0.38\linewidth}}
\hline
\textbf{Regime} & \textbf{Analytical properties} & \textbf{Steklov optimization outcome} \\
\hline
\textbf{Dimension $n=2$,} Unconstrained \cite{FS16, Kokarev, Pet19} & Spectral bounds determined by topology. & Smooth global maximizers, when they exist, are realized by free boundary minimal immersions in $\B^m$. \\
\textbf{Dimension $n=2$,} Equivariant \cite{FS16, FS11, Vin1} & $O(2)$-symmetry on the annulus. & The unique global optimum is achieved explicitly by the 2D critical catenoid. \\
\textbf{Dimension $n\ge 3$,} Unconstrained \cite{Bri1, CEG, GKL, GP2, KM}, Proposition \ref{prop:steklov_blowup_variable} & Unbounded spectrum; scale-invariant normalization. & Smooth global maximizers are generally non-existent. \\
\textbf{Dimension $n \ge 3$,} Equivariant (This Paper) & $O(n)$-symmetry on the cylinder; reduction to 1D. & Unique global maximizer realized by the critical $n$-catenoid. \\
\hline
\end{tabular}
\vspace{0.2cm}
\caption{Panorama of geometric regimes in Steklov eigenvalue optimization.}
\label{tab:landscape}
\end{table}

\subsection{Formulation of the Problem and Main Results}

We resolve these analytical challenges through symmetry-enforced dimensional reduction. While this strategy has proven successful in the context of surfaces (see, e.g., \cite{KKMS, KN, Mor, Vin1}), we extend the approach to higher dimensions. We study the maximization of $\bar{\sigma}_1$ on the fixed $n$-dimensional cylinder $\Omega = [-1, 1] \times \Sbb^{n-1}$. We impose an $O(n)$-invariance constraint on the metrics and boundary weights, which reduces the $n$-dimensional PDE to a system of one-dimensional ODEs on the conformal models $\Omega_L = [-L, L] \times \Sbb^{n-1}$ parameterized by the conformal modulus $L \in (0, \infty)$. The multiplicity of $\sigma_1$ in this symmetric setting is uniformly bounded by $n+1$, preventing the unbounded multiplicity issues appearing in general high-dimensional optimization.

\begin{figure}[htpb]
    \centering
    \begin{tikzpicture}[scale=0.9, >=stealth]

    \begin{scope}[shift={(-4.5, 0)}]
        \draw[thick] (-1.2, -2.4) -- (-1.2, 2.4);
        \draw[thick] (1.2, -2.4) -- (1.2, 2.4);

        \draw[thick, fill=blue!5] (0,2.4) ellipse (1.2 and 0.3);
        \node[above] at (0, 2.9) {$\Gamma_L \ (t = L)$};

        \draw[thick, fill=blue!5] (1.2,-2.4) arc (0:-180:1.2 and 0.3);
        \draw[dashed, color=gray!80, thick] (1.2,-2.4) arc (0:180:1.2 and 0.3);
        \node[below] at (0, -2.9) {$\Gamma_{-L} \ (t = -L)$};

        \draw[dashed, opacity=0.4] (1.2,0) arc (0:-180:1.2 and 0.3);
        \draw[dashed, opacity=0.4] (1.2,0) arc (0:180:1.2 and 0.3);

        \node[fill=white, fill opacity=0.8, text opacity=1, inner sep=2pt] at (0,0) {$\Omega_L = [-L, L] \times \mathbb{S}^{n-1}$};
    \end{scope}

    \draw[->, thick, shorten >= 10pt, shorten <= 10pt] (-2.8, 0.5) to[bend left=20] node[midway, above] {$\Phi$} (-0.2, 0.5);

    \begin{scope}[shift={(3.5, 0)}]
        \def\R{2.8844} 

        \draw[thick, densely dashed, fill=gray!5, fill opacity=0.5] (0,0) circle (\R);
        \draw[densely dashed, opacity=0.5] (\R,0) arc (0:-180:\R\space and 0.6);
        \draw[densely dashed, opacity=0.5] (\R,0) arc (0:180:\R\space and 0.6);
        \node[above left, xshift=-2pt, yshift=2pt] at (135:\R) {$\partial \mathbb{B}^{n+1}$};

        \fill[blue!10, opacity=0.6] 
            (-1.6, 2.4) arc (-180:0:1.6 and 0.3) 
            .. controls (0.8, 1.2) and (0.8, 0.4) .. (0.8, 0)
            .. controls (0.8, -0.4) and (0.8, -1.2) .. (1.6, -2.4) 
            arc (0:-180:1.6 and 0.3) 
            .. controls (-0.8, -1.2) and (-0.8, -0.4) .. (-0.8, 0)
            .. controls (-0.8, 0.4) and (-0.8, 1.2) .. (-1.6, 2.4); 

        \draw[thick, color=blue!70!black] (1.6, 2.4) .. controls (0.8, 1.2) and (0.8, 0.4) .. (0.8, 0) .. controls (0.8, -0.4) and (0.8, -1.2) .. (1.6, -2.4);
        \draw[thick, color=blue!70!black] (-1.6, 2.4) .. controls (-0.8, 1.2) and (-0.8, 0.4) .. (-0.8, 0) .. controls (-0.8, -0.4) and (-0.8, -1.2) .. (-1.6, -2.4);

        \draw[dashed, color=blue!70!black, thick] (1.6, 2.4) arc (0:180:1.6 and 0.3);
        \draw[thick, color=blue!70!black] (1.6, 2.4) arc (0:-180:1.6 and 0.3);

        \draw[dashed, color=blue!70!black, thick] (1.6, -2.4) arc (0:180:1.6 and 0.3);
        \draw[thick, color=blue!70!black] (1.6, -2.4) arc (0:-180:1.6 and 0.3);

        \node[fill=blue!10, fill opacity=0.7, text opacity=1, inner sep=1pt, color=blue!80!black] at (0, 1.4) {$\mathcal{C}_n$};

        \filldraw (0,0) circle (1.5pt) node[right, scale=0.8, xshift=2pt] {$0$};

        \coordinate (O) at (0,0);
        \coordinate (P) at (1.6, 2.4);

        \draw[thick, red!80!black, ->] (O) -- (P);
        \draw[thick, red!80!black, ->] (P) -- ($(O)!1.2!(P)$) node[above right, scale=0.95] {$\nu_{\mathbb{B}} = X$};

        \coordinate (T1) at ($(P) + 0.3*(-2.4, 1.6)$);
        \coordinate (T2) at ($(P) - 0.3*(-2.4, 1.6)$);
        \draw[thick, green!50!black] (T1) node[above left, scale=0.8, color=black] {$T_X \partial\mathbb{B}^{n+1}$} -- (T2); 

        \coordinate (Ndir) at ($(P)!0.2cm!(O)$); 
        \coordinate (Tdir) at ($(P)!0.2cm!90:(O)$); 
        \coordinate (Corner) at ($(Ndir)+(Tdir)-(P)$); 
        \draw[thick, red!80!black] (Ndir) -- (Corner) -- (Tdir);

        \node[anchor=north west, scale=0.85, text width=3.8cm, align=left] at (3.0, 2.4) {
            Profile tangent aligns \\
            with position vector $X$ \\
            $\implies t_* L'(t_*) = L(t_*)$
        };

    \end{scope}
    \end{tikzpicture}
    \caption{Conformal embedding of the cylinder $\Omega_L$ into the Euclidean unit ball $\mathbb{B}^{n+1}$ as the critical $n$-catenoid $\mathcal{C}_n$, which meets $\partial \mathbb{B}^{n+1}$ orthogonally at $t = \pm t_*$.}
    \label{fig:catenoid_mapping}
\end{figure}
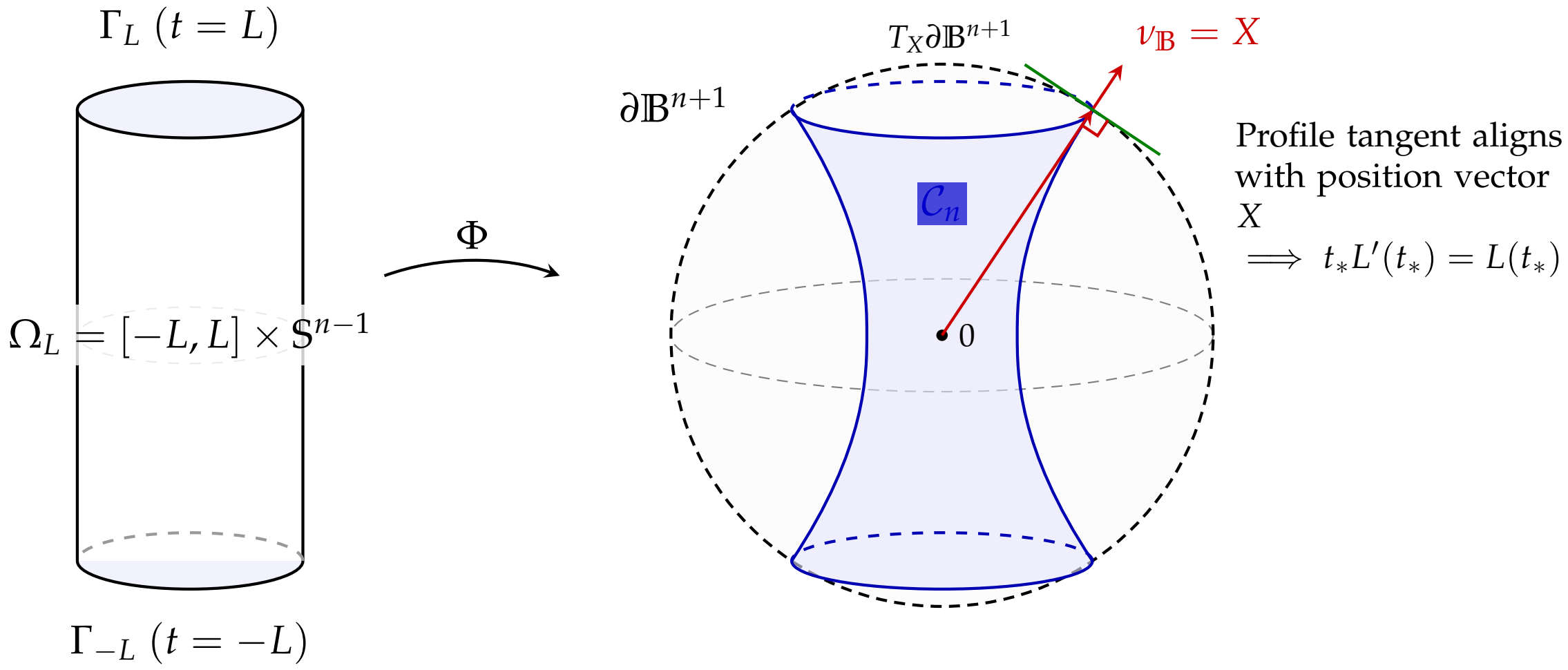

The maximal metric emerging from this dimensional reduction corresponds to the $n$-dimensional critical catenoid \cite{SSTZ} (see Figure \ref{fig:catenoid_mapping}). As the conformal modulus $L$ stretches toward infinity, the Steklov spectrum is bounded above by the normalized eigenvalue of the unit ball $\B^n$, which is strictly dominated by that of the critical $n$-catenoid. This proves that the global supremum must be achieved at a finite conformal modulus $L^* \in (0, \infty)$. Global criticality then shows that the critical $n$-catenoid uniquely realizes this supremum.

Our main result establishes both the existence and the precise geometric characterization of the maximizer for this problem.

\begin{theorem}[Main Theorem] \label{thm:main_theorem}
Let $\Omega = [-1, 1] \times \Sbb^{n-1}$ be the $n$-dimensional cylinder, where $n \ge 3$. The supremum of the normalized first Steklov eigenvalue $\bar{\sigma}_1(\Omega, g, \rho)$ over the space of all $O(n)$-invariant metrics and weights is finite and is achieved by a smooth metric $g_{\max}$ and weight $\rho_{\max}$. 

Furthermore, the maximal metric $g_{\max}$ is uniquely characterized, up to homothety and diffeomorphism, as the metric induced by the critical $n$-catenoid: an $O(n)$-equivariant free boundary minimal hypersurface in the Euclidean ball $\B^{n+1}$.
\end{theorem}

\subsection{Plan of the Paper}
The paper is organized as follows. Section \ref{sec:analytical_framework} establishes the analytical framework: we define $O(n)$-invariant metrics on $\Omega = [-1, 1] \times \Sbb^{n-1}$, show that each such metric is isometric to a conformally flat metric on $\Omega_L = [-L, L] \times \Sbb^{n-1}$ for a uniquely determined conformal modulus $L \in (0, \infty)$, and decompose the Steklov problem into a family of one-dimensional Sturm-Liouville ODEs. Section \ref{sec:bounds} proves uniform upper bounds for $\bar{\sigma}_1$ independent of the conformal modulus $L$, the conformal factor, and the boundary weights.

The proof of existence and geometric characterization of the global maximizer then unfolds in four steps across Sections \ref{sec:existence}--\ref{sec:criticality}. First, for a fixed conformal modulus $L$, a maximizing sequence can a priori degenerate by concentrating all its boundary mass onto a single boundary component, mimicking a punctured ball with asymptotic value $\bar{\sigma}_1(\B^n)$; Section \ref{sec:existence} therefore establishes a conditional existence result, proving via Ekeland's Variational Principle and Clarke's subdifferential calculus that whenever the strict threshold $\Sigma_1(L) > \bar{\sigma}_1(\B^n)$ holds, boundary collapse is ruled out and $\Sigma_1(L)$ is achieved by a smooth metric and weight associated with a free boundary $n$-harmonic embedding. Second, to show that this threshold is actually met at a finite optimal modulus, Section \ref{sec:moduli_space} studies the continuity of the supremum map $L \mapsto \Sigma_1(L)$ over the moduli space $(0, \infty)$ and proves that it vanishes as $L \to 0$ and is bounded above by $\bar{\sigma}_1(\B^n)$ as $L \to \infty$. Third, Section \ref{sec:spectral_domination} explicitly constructs the critical free boundary $n$-catenoid $\mathcal{C}_n$ as a competitor and proves directly that $\bar{\sigma}_1(\mathcal{C}_n) > \bar{\sigma}_1(\B^n)$; combined with the asymptotic bounds of Section \ref{sec:moduli_space}, this traps the global supremum at a finite modulus $L^* \in (0, \infty)$ satisfying $\Sigma_1(L^*) \ge \bar{\sigma}_1(\mathcal{C}_n) > \bar{\sigma}_1(\B^n)$, which in turn triggers the existence theorem of Section \ref{sec:existence} at $L^*$. Finally, Section \ref{sec:criticality} establishes the global criticality conditions in the space of all $O(n)$-invariant metrics, proving that since the maximizer at $L^*$ is critical under all $O(n)$-invariant variations (not merely conformal ones), its associated $n$-harmonic embedding turns into a free boundary minimal embedding into $\B^{n+1}$ and is therefore congruent to the critical $n$-catenoid. Appendix \ref{sec:numerics_3D} provides numerical data and asymptotic behavior as $n \to \infty$, and Appendix \ref{sec:maximization_of_Bn} describes the maximization in the conformal class of $\B^n$.

\subsection{Acknowledgements}
The author would like to thank Mikhail Karpukhin for suggesting this problem. The author would also like to thank Iosif Polterovich, Bruno Colbois, Luigi Provenzano, and Denis Vinokurov for fruitful discussions and valuable comments. This project is part of the author's PhD thesis at Universit\'e de Montr\'eal under the supervision of Iosif Polterovich and Mikhail Karpukhin. The author is supported by the ISM and the Banque Nationale ESP scholarship.

\subsection{LLM Usage Disclosure}
The development of this work used a workflow utilizing Gemini Pro 3.1, including for discussing proof strategies, editing, coding images and assisting with calculations. All mathematical statements, proofs, and verifications were subsequently completed and rigorously checked by the author, who takes full responsibility for the results.

\section{Geometric Preliminaries \& Separation of Variables}
\label{sec:analytical_framework}

To solve the Steklov maximization problem on the $n$-dimensional cylinder, we first establish the geometric and analytical framework that governs our optimization. We start from the fixed smooth manifold
\[
    \Omega := [-1, 1] \times \Sbb^{n-1},
\]
and show that any $O(n)$-invariant metric on $\Omega$ is isometric to a conformally flat metric on a standard cylinder $\Omega_L := [-L, L] \times \Sbb^{n-1}$, where $L \in (0, \infty)$ is the conformal modulus (conformal half-length) determined by the metric. Subsequently, we exploit this cylindrical symmetry to decompose the multidimensional Steklov boundary value problem into a family of one-dimensional Sturm-Liouville equations.

\subsection{Geometric Equivalence of Maximization Problems}

\begin{definition}[$O(n)$-invariant metrics and weights] \label{def:On_invariant}
Consider the natural action of the orthogonal group $O(n)$ on the cylinder $\Omega = [-1, 1] \times \Sbb^{n-1}$ given by
\[
    \kappa \cdot (r, \xi) := (r, \kappa \xi), \quad \forall \kappa \in O(n), \; (r, \xi) \in [-1, 1] \times \Sbb^{n-1}.
\]
A smooth Riemannian metric $g$ on $\Omega$ is said to be \emph{$O(n)$-invariant} if the diffeomorphism $R_\kappa(r, \xi) := \kappa \cdot (r, \xi)$ is an isometry of $(\Omega, g)$ for every $\kappa \in O(n)$, that is, $R_\kappa^* g = g$. Similarly, a boundary weight $\rho \in C^\infty_{>0}(\partial \Omega)$ is \emph{$O(n)$-invariant} if $\rho \circ R_\kappa = \rho$ for all $\kappa \in O(n)$, which means that $\rho$ takes constant positive values $\rho_-$ and $\rho_+$ on the two boundary components $\{-1\} \times \Sbb^{n-1}$ and $\{1\} \times \Sbb^{n-1}$.
\end{definition}

\begin{proposition} \label{prop:moduli_equivalence}
Let $n \ge 3$. Any $O(n)$-invariant smooth Riemannian metric $g$ on $\Omega = [-1, 1] \times \Sbb^{n-1}$ is isometric to a conformally flat metric of the form
\[
    g_L = f(t)^2 \left( dt^2 + g_{\Sbb^{n-1}} \right) \quad \text{on } \Omega_L := [-L, L] \times \Sbb^{n-1},
\]
for a uniquely determined conformal modulus $L \in (0, \infty)$ and a smooth positive function $f \in C^\infty([-L, L])$. Consequently, the supremum of $\bar{\sigma}_1(\Omega, g, \rho)$ over all $O(n)$-invariant metrics and boundary weights on $\Omega$ is equal to the supremum over $L \in (0, \infty)$ of $\bar{\sigma}_1(\Omega_L, g_L, \rho)$ among conformal metrics $g_L = f(t)^2(dt^2 + g_{\Sbb^{n-1}})$ and $O(n)$-invariant boundary weights on $\Omega_L$.
\end{proposition}

\begin{proof}
Let $g$ be an $O(n)$-invariant metric on $\Omega = [-1, 1] \times \Sbb^{n-1}$. The principal orbits of the $O(n)$ action are the cross-sections $\{r\} \times \Sbb^{n-1}$ for $r \in [-1, 1]$. Relative to this product structure, $g$ takes the general form
\[
    g = A(r)dr^2 + 2 dr \odot \omega_r + h_r,
\]
where $A(r) > 0$ is a smooth function on $[-1, 1]$ (since $O(n)$ acts transitively on $\Sbb^{n-1}$), $\omega_r$ is a smooth family of 1-forms on $\Sbb^{n-1}$, and $h_r$ is a smooth family of Riemannian metrics on $\Sbb^{n-1}$.

Since $g$ is $O(n)$-invariant, both $\omega_r$ and $h_r$ must be invariant under the action of $O(n)$ on $\Sbb^{n-1}$. For $n \ge 3$, the isotropy group fixing a point $\xi \in \Sbb^{n-1}$ is isomorphic to $O(n-1)$, and the standard representation of $O(n-1)$ on $T_\xi \Sbb^{n-1} \cong \R^{n-1}$ admits no non-zero fixed vectors. Consequently, there exist no non-zero $O(n)$-invariant 1-forms on $\Sbb^{n-1}$, which forces $\omega_r \equiv 0$ (in contrast with the $n=2$ case, where the isotropy group on $SO(2)$ is trivial and allows the invariant 1-form $d\theta$).

Furthermore, since the isotropy representation of $O(n-1)$ on $T_\xi \Sbb^{n-1}$ is irreducible, Schur's lemma implies that the only $O(n)$-invariant symmetric $(0,2)$-tensor on $\Sbb^{n-1}$ is, up to a positive scalar factor, the standard round metric $g_{\Sbb^{n-1}}$. Thus, $h_r = B(r)g_{\Sbb^{n-1}}$ for a smooth positive function $B(r) > 0$ on $[-1, 1]$, and the metric diagonalizes as
\[
    g = A(r)dr^2 + B(r)g_{\Sbb^{n-1}}.
\]

Define the conformal modulus (conformal half-length) of $(\Omega, g)$ by
\[
    L := \frac{1}{2} \int_{-1}^1 \sqrt{\frac{A(s)}{B(s)}} \, ds \in (0, \infty),
\]
and introduce the longitudinal coordinate change $\tau \colon [-1, 1] \to [-L, L]$ given by
\[
    t = \tau(r) := \int_{-1}^r \sqrt{\frac{A(s)}{B(s)}} \, ds - L.
\]
Then $\tau$ is a smooth diffeomorphism from $[-1, 1]$ onto $[-L, L]$ satisfying $dt^2 = \frac{A(r)}{B(r)}dr^2$. Setting $f(t) := \sqrt{B(\tau^{-1}(t))}$, the diffeomorphism $\Psi(r, \xi) := (\tau(r), \xi)$ from $\Omega$ onto $\Omega_L = [-L, L] \times \Sbb^{n-1}$ pulls back the conformally flat metric $g_L = f(t)^2(dt^2 + g_{\Sbb^{n-1}})$ to $g$:
\[
    \Psi^* g_L = B(r) \left( \frac{A(r)}{B(r)} dr^2 + g_{\Sbb^{n-1}} \right) = g.
\]
Since the normalized Steklov functional $\bar{\sigma}_1$ is invariant under isometries, maximizing $\bar{\sigma}_1$ over all $O(n)$-invariant metrics and weights on $\Omega$ is equivalent to maximizing $\bar{\sigma}_1$ over all $L \in (0, \infty)$, all conformal factors $f \in C^\infty_{>0}([-L, L])$, and all invariant weights $\rho_\pm > 0$ on $\partial \Omega_L$.
\end{proof}

To compute and bound eigenvalues, we explicitly evaluate the Steklov spectrum of these conformal cylinders using separation of variables.

\subsection{Spectral Analysis of $O(n)$-Invariant Metrics}
\label{sec:computation}

We now provide an explicit characterization of the Steklov spectrum for the conformal cylinders $\Omega_L = [-L, L] \times \Sbb^{n-1}$. By exploiting the symmetries of the manifold, we reduce the Steklov boundary value problem from a multidimensional partial differential equation to a family of one-dimensional Sturm-Liouville problems along the longitudinal profile.

On $\Omega_L = [-L, L] \times \Sbb^{n-1}$, an $O(n)$-invariant conformal metric $g$ takes the warped product form:
\[
    g = f(t)^2 \left( dt^2 + g_{\Sbb^{n-1}} \right), \quad t \in [-L, L],
\]
where $f > 0$ is a smooth warping function. The boundary $\partial \Omega_L = \Gamma_{-L} \cup \Gamma_L$ consists of the two parallel spheres $\{-L\} \times \Sbb^{n-1}$ and $\{L\} \times \Sbb^{n-1}$. We consider an $O(n)$-invariant boundary weight $\rho$ taking constant positive values $\rho_-$ and $\rho_+$ on $\Gamma_{-L}$ and $\Gamma_L$, respectively.

Let $0 = \lambda_0 < \lambda_1 \le \lambda_2 \le \dots$ be the distinct eigenvalues of the Laplace-Beltrami operator $\Delta_{\Sbb^{n-1}}$ on the standard unit sphere. It is a classical result that $\lambda_k = k(k+n-2)$ (see, e.g., \cite{LMP}), with eigenspaces spanned by the spherical harmonics $\beta^{k, i}(\xi)$, $i=1,\dots,m_{\lambda_k}$, where $m_{\lambda_k}$ is the multiplicity of $\lambda_k$.

\begin{theorem} \label{thm:spectrum_computation}
The Steklov spectrum of $(\Omega_L, g, \rho)$ decomposes into discrete branches $\sigma^{(k)}$ corresponding to the spherical harmonics of $\Sbb^{n-1}$. There is a Steklov eigenbasis of the separated form
\[
    u_{k, i}(t, \xi) = Y_k(t) \beta^{k, i}(\xi), \quad t \in [-L, L], \;\; \xi \in \Sbb^{n-1},
\]
where $Y_k(t)$ is a smooth solution to the second-order ordinary differential equation
\begin{equation} \label{eq:profile_ode}
    Y''(t) + (n-2)\frac{f'(t)}{f(t)}Y'(t) - \lambda_k Y(t) = 0,
\end{equation}
subject to the Steklov boundary conditions:
\begin{equation} \label{eq:profile_bc}
    \begin{cases}
    Y'(L) = \sigma^{(k)} \rho_+ f(L) Y(L), \\
    -Y'(-L) = \sigma^{(k)} \rho_- f(-L) Y(-L).
    \end{cases}
\end{equation}
The first non-zero Steklov eigenvalue is $\sigma_1(g, \rho) = \min\{\sigma^{(0)}_1, \sigma^{(1)}_1\}$, where the first non-zero longitudinal eigenvalue is
\begin{equation} \label{eq:sigma_0_exact}
    \sigma^{(0)}_1 = \frac{\rho_+^{-1}f(L)^{1-n} + \rho_-^{-1}f(-L)^{1-n}}{\int_{-L}^L f(t)^{2-n} \, dt},
\end{equation}
and the first transverse eigenvalue is characterized variationally by
\begin{equation}\label{eq:rayleigh_inf}
    \sigma^{(1)}_1 = \inf_{\psi \not\equiv 0} \frac{\int_{-L}^L \left(\psi'(t)^2 + (n-1)\psi(t)^2\right) f(t)^{n-2} \, dt}{\rho_+ f(L)^{n-1} \psi(L)^2 + \rho_-f(-L)^{n-1} \psi(-L)^2}.
\end{equation}
\end{theorem}

\begin{proof}
We search for harmonic functions $u \in C^\infty(\Omega_L)$ satisfying $\Delta_g u = 0$. Using the formula for the Laplacian of a warped product, substituting the separated ansatz $u(t, \xi) = Y(t)\beta(\xi)$ yields:
\[
    \frac{\left( f^{n-2}(t) Y'(t) \right)'}{f^{n-2}(t) Y(t)} = -\frac{\Delta_{\Sbb^{n-1}} \beta(\xi)}{\beta(\xi)} = \lambda.
\]
For $u$ to be harmonic, $\lambda$ must be an eigenvalue of $\Sbb^{n-1}$ ($\lambda = \lambda_k$), forcing $\beta$ to be a spherical harmonic of degree $k$. Expanding the left-hand side differential and dividing by $f^{n-2}$ yields the governing equation \eqref{eq:profile_ode}. Evaluating the outward normal derivative $\partial_\nu = \pm f(\pm L)^{-1}\partial_t$ on $\Gamma_{\pm L}$ yields the Steklov boundary conditions \eqref{eq:profile_bc}. Because $g$ is a warped product, $\Delta_g$ commutes with the spherical Laplacian $\Delta_{\Sbb^{n-1}}$, yielding the orthogonal decomposition of $H^1(\Omega_L)$ into subspaces indexed by the spectrum of $\Sbb^{n-1}$:
\[
    H^1(\Omega_L) = \bigoplus_{k=0}^\infty \mathcal{H}_k, \quad \text{where } \mathcal{H}_k = \text{span}\{ \psi(t)\beta^{k,i}(\xi) \mid \psi \in H^1([-L,L]),\, i=1,\dots,m_{\lambda_k} \}.
\]
This decomposition captures the entire Steklov spectrum of $\Omega_L$, which splits into orthogonal branches. For any $k \ge 1$, the spherical harmonics are $L^2$-orthogonal to constants on $\Sbb^{n-1}$. Thus, any separated function $u(t,\xi) = \psi(t)\beta^{k,i}(\xi)$ with $k \ge 1$ is orthogonal on $\partial \Omega_L$ (with respect to $\rho \, d\mathcal{H}_g$) to the constant Steklov eigenfunctions corresponding to $\sigma_0 = 0$.

To characterize the first eigenvalue of the $k$-th branch, $\sigma_1^{(k)}$ ($k \ge 1$), we evaluate the Steklov Rayleigh quotient $R_{g,\rho}(u) = \|\nabla_g u\|_{L^2(\Omega_L)}^2 / \|u\|_{L^2(\partial \Omega_L, \rho)}^2$ on $\mathcal{H}_k$. Using the volume form $dv_g = f(t)^n \, dt \, dv_{\Sbb^{n-1}}$, the Dirichlet energy evaluates to:
\[
    \int_{\Omega_L} |\nabla_g u|_g^2 \, dv_g = \int_{-L}^L f(t)^{n-2} \left[ \psi'(t)^2 \|\beta^{k,i}\|_{L^2(\Sbb^{n-1})}^2 + \psi(t)^2 \|\nabla \beta^{k,i}\|_{L^2(\Sbb^{n-1})}^2 \right] \, dt.
\]
Since $\|\nabla \beta^{k,i}\|_{L^2}^2 = \lambda_k \|\beta^{k,i}\|_{L^2}^2$ and the boundary measure is $d\mathcal{H}_g = f(\pm L)^{n-1} \, dv_{\Sbb^{n-1}}$, the spherical integrals cancel and the Rayleigh quotient reduces to the unconstrained one-dimensional quotient:
\[
    R_k(\psi) = \frac{\int_{-L}^L \left(\psi'(t)^2 + \lambda_k \psi(t)^2\right) f(t)^{n-2} \, dt}{\rho_+ f(L)^{n-1} \psi(L)^2 + \rho_-f(-L)^{n-1} \psi(-L)^2}.
\]

Since $\lambda_k = k(k+n-2)$ is strictly increasing with $k$, we have $\lambda_k > \lambda_1$ for all $k \ge 2$. Let $\psi_k \not\equiv 0$ be an eigenprofile realizing $\sigma_1^{(k)}$ for some $k \ge 2$. Substituting $\psi_k$ into $R_1$ lowers the numerator and yields the strict inequality:
\[
    \sigma_1^{(k)} = R_k(\psi_k) > R_1(\psi_k) \ge \inf_{\psi \not\equiv 0} R_1(\psi) = \sigma_1^{(1)}.
\]
Thus, $\sigma_1^{(k)} > \sigma_1^{(1)}$ for all $k \ge 2$, and the first non-zero Steklov eigenvalue is given by $\sigma_1(g, \rho) = \min\{\sigma^{(0)}_1, \sigma^{(1)}_1\}$.
\vspace{0.25cm}

\paragraph{\textbf{The Longitudinal Branch ($\lambda_0 = 0$).}}
Here, the spherical harmonic is constant, and the ODE reduces to $(f^{n-2} Y')' = 0$. Integrating directly yields $Y'(t) = C f(t)^{2-n}$. Integrating again from $-L$ to $L$ gives:
\[
    Y(L) - Y(-L) = C \int_{-L}^L f(t)^{2-n} \, dt.
\]
Applying the boundary conditions \eqref{eq:profile_bc} gives $Y(L) = C (\sigma_1^{(0)})^{-1} \rho_+^{-1} f(L)^{1-n}$ and $Y(-L) = -C (\sigma_1^{(0)})^{-1} \rho_-^{-1} f(-L)^{1-n}$. Substituting these into the relation above and cancelling $C \neq 0$ yields \eqref{eq:sigma_0_exact}.
\vspace{0.25cm}

\paragraph{\textbf{The Transverse Branch ($\lambda_1 = n-1$).}}
For $k = 1$, we have $\lambda_1 = n-1$. Substituting $\lambda_1$ into the one-dimensional quotient $R_1(\psi)$ yields the variational characterization \eqref{eq:rayleigh_inf}.
\end{proof}

\begin{corollary}\label{cor:multiplicity}
For any $O(n)$-invariant metric $g$ and invariant weight $\rho$ on $\Omega_L$, the multiplicity of the first Steklov eigenvalue $\sigma_1(\Omega_L, g, \rho)$ is bounded above by $n+1$.
\end{corollary}

\begin{proof}
By Theorem \ref{thm:spectrum_computation}, the first non-zero Steklov eigenvalue is achieved by the longitudinal ($k=0$) or transverse ($k=1$) branches. The longitudinal branch ($k=0$) corresponds to constant functions on $\Sbb^{n-1}$ ($m_{\lambda_0} = 1$). Its radial ODE $(f^{n-2}Y')'=0$ has a two-dimensional solution space spanned by constants (giving $\sigma_0 = 0$) and one non-constant profile, yielding a one-dimensional eigenspace for $\sigma_1^{(0)}$. The transverse branch ($k=1$) corresponds to the first spherical harmonics on $\Sbb^{n-1}$ ($m_{\lambda_1} = n$). By standard ODE theory (see Coddington \cite[Chapter 3]{Cod}), the second-order equation \eqref{eq:profile_ode} has a two-dimensional solution space on $[-L, L]$, and $\sigma_1^{(1)}$ is one of the corresponding eigenvalues. For a fixed $\sigma$, the space of solutions satisfying the boundary conditions \eqref{eq:profile_bc} is one-dimensional: indeed, if there were two linearly independent solutions, a non-trivial linear combination $Y$ would satisfy $Y(-L) = 0$, which by \eqref{eq:profile_bc} implies $Y'(-L) = 0$, forcing $Y \equiv 0$ by Cauchy-Lipschitz uniqueness. Thus, the transverse branch spans an $n$-dimensional eigenspace. When the two branches cross ($\sigma^{(0)}_1 = \sigma^{(1)}_1$), their eigenspaces sum directly, bounding the multiplicity by $1 + n = n+1$.
\end{proof}

\begin{remark}[Low-Regularity Spectrum via Sturm-Liouville Theory]
The spectral decomposition of Theorem \ref{thm:spectrum_computation} and the multiplicity bound of Corollary \ref{cor:multiplicity} remain valid for conformal factors $f \in L^n([-L,L])$ satisfying a uniform lower bound $f \ge \epsilon > 0$ almost everywhere. Setting $\alpha = f^{n-2}$, the ODE \eqref{eq:profile_ode} writes in divergence form as $(\alpha Y')' = \lambda_k \alpha Y$. Under the lower bound $f \ge \epsilon$, we have $\alpha \in L^{\frac{n}{n-2}}([-L, L]) \subset L^1([-L, L])$ and $\alpha^{-1} \in L^\infty([-L, L]) \subset L^1([-L, L])$. By Sturm-Liouville theory (see Zettl \cite[Chapter 2]{Zet}), the integrability of both $\alpha$ and $\alpha^{-1}$ guarantees a regular Sturm-Liouville problem, preserving the discrete spectrum and the multiplicity bound along the low-regularity maximizing sequence constructed in Section \ref{sec:existence}.
\end{remark}

Notice that in the Rayleigh quotient \eqref{eq:rayleigh_inf}, the conformal factor $f$ appears as $f^{n-2}$ in the bulk integral and as $f(\pm L)^{n-1}$ in the boundary terms. To formulate a well-posed variational problem in a Banach space without requiring point evaluations of $f$ at $\pm L$, we decouple the bulk density from the boundary measure.

\begin{definition}[Space of admissible densities and measures $\mathcal{M}_L$] \label{def:space_ML}
Let $L > 0$ and set $q := \frac{n}{n-2}$. We define the space of admissible $O(n)$-invariant bulk densities and boundary measures on $[-L, L]$ by
\[
    \mathcal{M}_L := \left\{ (\alpha, \mu) \in L^q([-L, L]) \times \left(\R_{\ge 0}^2 \setminus \{(0,0)\}\right) \;\middle|\; \operatorname*{ess\,inf}_{t \in [-L, L]} \alpha(t) > 0 \right\},
\]
where $\mu = (\mu_{-L}, \mu_L)$ represents the effective boundary weights on $\Gamma_{-L}$ and $\Gamma_L$ (corresponding to $\mu_{\pm L} = \rho_\pm f(\pm L)^{n-1}$ and $\alpha(t) = f(t)^{n-2}$ in the smooth setting). For any $(\alpha, \mu) \in \mathcal{M}_L$, the first non-zero Steklov eigenvalue is given by the expression: $\sigma_1(\alpha, \mu) := \min\{\sigma_1^{(0)}(\alpha, \mu), \sigma_1^{(1)}(\alpha, \mu)\}$, where
\[
    \sigma_1^{(0)}(\alpha, \mu) := \inf_{\substack{\psi \in W^{1,n}([-L, L]) \setminus \{0\} \\ \mu_{-L}\psi(-L) + \mu_L\psi(L) = 0}} \frac{\int_{-L}^L \psi'(t)^2 \alpha(t) \, dt}{\mu_L \psi(L)^2 + \mu_{-L}\psi(-L)^2} = \frac{\mu_L^{-1} + \mu_{-L}^{-1}}{\int_{-L}^L \alpha(t)^{-1} \, dt}
\]
(with the convention that $\sigma_1^{(0)}(\alpha, \mu) = +\infty$ if $\mu_L = 0$ or $\mu_{-L} = 0$), and
\[
    \sigma_1^{(1)}(\alpha, \mu) := \inf_{\psi \in W^{1,n}([-L, L]) \setminus \{0\}} \frac{\int_{-L}^L \left(\psi'(t)^2 + (n-1)\psi(t)^2\right) \alpha(t) \, dt}{\mu_L \psi(L)^2 + \mu_{-L}\psi(-L)^2}.
\]
The scale-invariant normalized Steklov functional on $\mathcal{M}_L$ then reads
\[
    \bar{\sigma}_1(\alpha, \mu) := \sigma_1(\alpha, \mu) \left( \mu_{-L} + \mu_L \right) \|\alpha\|_{L^q([-L, L])}^{-1} \vol(\Sbb^{n-1})^{\frac{2}{n}},
\]
and the supremum for a fixed conformal modulus $L \in (0, \infty)$ is defined by
\begin{equation} \label{eq:sigma_L_supremum}
    \Sigma_1(L) := \sup \left\{ \bar{\sigma}_1(\alpha, \mu) \;\middle|\; (\alpha, \mu) \in \mathcal{M}_L \right\}.
\end{equation}
\end{definition}

\section{A Priori Spectral Bounds}
\label{sec:bounds}

The Steklov spectrum is sensitive to both dimension and metric symmetry. While normalized eigenvalues are uniformly bounded for surfaces ($n=2$), Colbois, El Soufi, and Girouard \cite{CEG} demonstrated that for $n \ge 3$, geometric degeneration can drive eigenvalues arbitrarily large while keeping the boundary volume fixed. We show here that $O(n)$-invariance ensures the spectrum remains uniformly bounded.

\subsection{Conformal upper bounds via the conformal volume}
\label{sec:conformal_upper_bounds}

The foundational work of Li and Yau \cite{LY1, LY2}, generalized by El Soufi and Ilias \cite{ESvol} and Korevaar \cite{Kor}, establishes that for any conformal class $[g]$, Laplace eigenvalues satisfy uniform bounds related to the conformal volume (see also \cite{ColMae, GNY}). This approach, adapted by Hassannezhad \cite{Hass} using $n$-$n$ capacitors and the conformal volume $V([g])$, provides a general framework for measure-based eigenvalue problems. In the Steklov setting, where the measure is supported on $\partial \Omega$, this machinery guarantees that the supremum of the normalized eigenvalues is finite within a fixed conformal class (see also \cite{CG}).

\begin{proposition}[\cite{Hass}]\label{prop:steklov_isop_bounds}
Let $\Omega^n$ be a compact manifold with smooth boundary $\partial \Omega$ ($n \geq 2$). There exist universal constants $A_n, B_n > 0$ such that for every conformal class $[g]$ on $\Omega$, every positive boundary density $\rho \in C^\infty_{>0}(\partial \Omega)$, and every integer $k \geq 1$:
\[
\bar{\sigma}_k(\Omega, g, \rho) \leq A_n (V[g])^{\frac{2}{n}} + B_n k^{\frac{2}{n}}.
\]
\end{proposition}

This bound ensures that maximizing sequences within a fixed conformal class remain bounded from above.

\subsubsection*{Arbitrarily large normalized eigenvalues}

In dimensions $n \ge 3$, it is well established that for several normalizations, one can construct metric sequences where normalized eigenvalues blow up (see, e.g., \cite{Ble, CD, Bri1, Bri2, GL, GKL, Mat1, Mat2}). Below we start by noting that even for the scale-invariant normalization \eqref{eq:KM_normalization}, the spectrum remains unbounded even when the number of boundary components is fixed.

\begin{proposition}
\label{prop:steklov_blowup_variable}
Let $M$ be a fixed closed Riemannian manifold of dimension $n \geq 3$. For any integer $b \geq 1$, there exists a sequence of metrics $(g_\varepsilon)_\varepsilon$ on $M$ and a sequence of subdomains $(\Omega_\varepsilon)_\varepsilon \subset M$ with exactly $b$ boundary components such that for all $k \ge 1$:
\[
\bar{\sigma}_k(\Omega_\varepsilon, g_\varepsilon, 1) \longrightarrow +\infty \quad \text{as } \varepsilon \to 0.
\]
\end{proposition}

\begin{proof}
By \cite{CD}, there exist metrics $(g_\delta)_\delta$ on $M$ such that $\bar{\lambda}_1(M, g_\delta) \to +\infty$. We approximate the Laplace spectrum of $M$ via the Steklov spectrum of perforated domains $(\Omega_{\varepsilon, \delta})_\varepsilon \subset M$. By \cite[Proposition 1.12]{GKL}, for each fixed $\delta$, we have $\bar{\sigma}_1(\Omega_{\varepsilon, \delta}, g_\delta, 1) \to \bar{\lambda}_1(M, g_\delta)$ as $\varepsilon \to 0$.

To restrict the geometry to exactly $b$ boundary components, for each fixed $(\varepsilon, \delta)$ we apply the surgery construction of Fraser and Schoen \cite[Theorem 1.2]{FS19}: by removing sufficiently thin tubular neighborhoods of radius $\eta(\varepsilon, \delta) \ll \varepsilon$ around paths connecting the excess boundary components, we extract a smooth subdomain $\Omega'_{\varepsilon, \delta} \subset \Omega_{\varepsilon, \delta}$ with $b$ boundary components such that $\sigma_1(\Omega'_{\varepsilon, \delta}, g_\delta, 1)$, $\vol(\Omega'_{\varepsilon, \delta}, g_\delta)$, and $\mathcal{H}_{g_\delta}(\partial \Omega'_{\varepsilon, \delta})$ are arbitrarily close to those of $\Omega_{\varepsilon, \delta}$, ensuring that
\[
    \bar{\sigma}_1(\Omega'_{\varepsilon, \delta}, g_\delta, 1) \ge \bar{\sigma}_1(\Omega_{\varepsilon, \delta}, g_\delta, 1) - \varepsilon.
\]
Selecting a diagonal subsequence $\delta(\varepsilon) \to 0$ yields for all $k \ge 1$:
\[
    \bar{\sigma}_k(\Omega'_{\varepsilon, \delta(\varepsilon)}, g_{\delta(\varepsilon)}, 1) \ge \bar{\sigma}_1(\Omega'_{\varepsilon, \delta(\varepsilon)}, g_{\delta(\varepsilon)}, 1) \longrightarrow +\infty \quad \text{as } \varepsilon \to 0.
\]
\end{proof}

While Proposition \ref{prop:steklov_blowup_variable} controls the number of boundary components, the interior topology of $\Omega_\varepsilon$ varies with $\varepsilon$. Nevertheless, one expects that unboundedness persists on any fixed smooth domain:

\begin{conjecture}
\label{conj:steklov_blowup_fixed}
Let $\Omega$ be a fixed compact manifold of dimension $n \ge 3$ with $b \ge 1$ boundary components. There exists a sequence of metrics $h_\varepsilon$ on $\Omega$ such that for all $k \ge 1$:
\[
\bar{\sigma}_k(\Omega, h_\varepsilon, 1) \longrightarrow +\infty \quad \text{as } \varepsilon \to 0.
\]
\end{conjecture}

\begin{remark}
If one smoothly embeds $\Omega$ into a closed $n$-dimensional manifold $M$ (for instance, by capping off the boundary components), applying Proposition \ref{prop:steklov_blowup_variable} to $M$ yields subdomains $\Omega_\varepsilon \subset M$ with $b$ boundary components such that $\bar{\sigma}_k(\Omega_\varepsilon, g_\varepsilon, 1) \to +\infty$. The difficulty is that $\Omega_\varepsilon$ is not a priori diffeomorphic to $\Omega$; proving Conjecture \ref{conj:steklov_blowup_fixed} requires performing the metric surgery while preserving the diffeomorphism type of $\Omega$.
\end{remark}

\subsection{Uniform Bound for $O(n)$-Invariant Metrics}

Our objective is now to establish a uniform upper bound for $\bar{\sigma}_1$ over all $O(n)$-invariant metrics and weights on the cylinder $\Omega = [-1, 1] \times \Sbb^{n-1}$, or equivalently, over all conformal moduli $L \in (0, \infty)$ and all pairs $(\alpha, \mu) \in \mathcal{M}_L$.

\begin{theorem} \label{thm:uniform_bound}
For $n \ge 3$, there exists a constant $C_n > 0$ depending only on $n$ such that for every conformal modulus $L > 0$ and every $(\alpha, \mu) \in \mathcal{M}_L$,
\[
    \bar{\sigma}_1(\alpha, \mu) \le C_n.
\]
In particular, $\bar{\sigma}_1(\Omega, g, \rho) \le C_n$ for all $O(n)$-invariant metrics $g$ and invariant weights $\rho$ on $\Omega = [-1, 1] \times \Sbb^{n-1}$.
\end{theorem}

\begin{proof}
By Definition \ref{def:space_ML}, for any $L > 0$, any $(\alpha, \mu) \in \mathcal{M}_L$, and any non-zero test profile $\psi \in W^{1,n}([-L, L])$, we have $\sigma_1(\alpha, \mu) \le \sigma_1^{(1)}(\alpha, \mu)$, and therefore
\[
    \bar{\sigma}_1(\alpha, \mu) \le \frac{\int_{-L}^L \left(\psi'(t)^2 + \lambda_1 \psi(t)^2\right) \alpha(t) \, dt}{\mu_L \psi(L)^2 + \mu_{-L} \psi(-L)^2} \cdot (\mu_L + \mu_{-L}) \, \|\alpha\|_{L^q([-L, L])}^{-1} \, \vol(\Sbb^{n-1})^{\frac{2}{n}},
\]
where $\lambda_1 = n-1$ and $q = \frac{n}{n-2}$ (note that in the smooth setting $\alpha = f^{n-2}$ and $\mu_{\pm L} = \rho_\pm f(\pm L)^{n-1}$, the factor $\mu_L + \mu_{-L} = \rho_+ f(L)^{n-1} + \rho_- f(-L)^{n-1}$ comes from the boundary mass $\int_{\partial \Omega_L} \rho \, d\mathcal{H}_g$).

\textbf{Case 1: $L < 1$.} Choosing the constant test profile $\psi \equiv 1$, the boundary factor $\mu_L + \mu_{-L}$ cancels between the numerator and the denominator. Applying H\"older's inequality with conjugate exponents $q = \frac{n}{n-2}$ and $p = \frac{n}{2}$ gives:
\[
    \int_{-L}^L \alpha(t) \, dt \le \|\alpha\|_{L^q([-L, L])} (2L)^{\frac{2}{n}}.
\]
Cancelling $\|\alpha\|_{L^q([-L, L])}$ yields $\bar{\sigma}_1(\alpha, \mu) \le \lambda_1 \vol(\Sbb^{n-1})^{\frac{2}{n}} (2L)^{\frac{2}{n}} \le \lambda_1 \vol(\Sbb^{n-1})^{\frac{2}{n}} 2^{\frac{2}{n}}$ for $L < 1$.

\textbf{Case 2: $L \ge 1$.} We localize $\psi$ near $\partial \Omega_L$ by choosing the piecewise linear cutoff satisfying $\psi(\pm L) = 1$, $\psi(t) = 0$ on $[-L+1, L-1]$, and $|\psi'(t)| = 1$ on $I := (-L, -L+1) \cup (L-1, L)$. Since $\psi(\pm L) = 1$, the boundary sum $\mu_L + \mu_{-L}$ again cancels out. Applying H\"older's inequality over $I$ (which has Lebesgue measure $|I| = 2$) gives:
\[
    \int_I \alpha(t) \, dt \le \|\alpha\|_{L^q(I)} 2^{\frac{2}{n}} \le \|\alpha\|_{L^q([-L, L])} 2^{\frac{2}{n}}.
\]
Since $\psi'^2 + \lambda_1 \psi^2 \le 1 + \lambda_1$ on $I$ and vanishes outside $I$, this yields
\[
    \bar{\sigma}_1(\alpha, \mu) \le (1 + \lambda_1) \vol(\Sbb^{n-1})^{\frac{2}{n}} 2^{\frac{2}{n}},
\]
which completes the proof with $C_n := n 2^{\frac{2}{n}} \vol(\Sbb^{n-1})^{\frac{2}{n}}$.
\end{proof}

\section{Existence of a Maximizer for a Fixed Conformal Modulus}
\label{sec:existence}

We now work on the standard cylinder $\Omega_L = [-L, L] \times \Sbb^{n-1}$ for a fixed conformal modulus $L \in (0, \infty)$, with background metric $g_0 = dt^2 + g_{\Sbb^{n-1}}$. In this section, we prove the existence of a smooth maximizing $O(n)$-invariant configuration $(\alpha_0, \mu_0) \in \mathcal{M}_L$ (and its associated free boundary $n$-harmonic map) whenever $L$ lies in the spectrally stable regime.

Because Steklov maximization on a cylinder is susceptible to boundary mass concentration onto a single boundary component (mimicking a topological ball $\B^n$), a smooth cylindrical maximizer need not exist for all lengths $L$. As shown in Section \ref{sec:moduli_space}, as $L \to \infty$, $\Sigma_1(L)$ is bounded above by $\bar{\sigma}_1(\B^n)$. We therefore establish our existence theorem for conformal moduli $L$ satisfying the strict inequality $\Sigma_1(L) > \bar{\sigma}_1(\B^n)$. Since the global supremum strictly exceeds this threshold (as realized by the critical $n$-catenoid in Section \ref{sec:spectral_domination}), this condition holds at the optimal modulus $L^*$.

Maximizing eigenvalues within a fixed conformal class generally lacks compactness (see, e.g., \cite{KNPP, Pet14, Vin1, Vin2}). To extract a regular limit without assuming a priori smoothness, we combine Ekeland's Variational Principle with Clarke's subdifferential calculus and the one-dimensional regularity of $O(n)$-equivariant $n$-harmonic maps.

\subsection{Preliminaries on Non-Smooth Variational Tools}
\label{sec:nonsmooth_preliminaries}

Maximizing $\bar{\sigma}_1$ presents two analytical difficulties: the lack of compactness in the space of conformal factors, and the non-differentiability of $\sigma_1$ at points of multiplicity. We address these using Ekeland's Variational Principle and Clarke's generalized differential.

\subsubsection*{Ekeland's Variational Principle}

In the absence of compactness, arbitrary maximizing sequences carry no differential information. Ekeland's Variational Principle \cite{Ekeland} extracts a maximizing sequence that is also \emph{almost-critical}.

\begin{theorem}[Ekeland's Variational Principle]
Let $(X, d)$ be a complete metric space, and $F \colon X \to \mathbb{R} \cup \{+\infty\}$ a lower semi-continuous functional bounded from below with $F\not\equiv+\infty$. For every $\gamma > 0$, $\delta > 0$, and $u \in X$ satisfying $F(u) \le \inf_X F + \delta$, there exists $v \in X$ such that:
\begin{itemize}
    \item[(i)] $F(v) \le F(u)$;
    \item[(ii)] $d(u, v) \le \gamma$;
    \item[(iii)] For all $w \neq v$, $F(v) < F(w) + \frac{\delta}{\gamma} \, d(v, w)$.
\end{itemize}
\end{theorem}

We apply this principle in Section \ref{sec:ekeland_application}. Property (iii) characterizes $v$ as a strict global minimizer of the perturbed functional $w \mapsto F(w) + \frac{\delta}{\gamma} d(v, w)$, yielding approximate Euler--Lagrange equations that relate the metric density to the energy density of the Steklov eigenfunctions.

\subsubsection*{Clarke Subdifferential Calculus}

To differentiate the locally Lipschitz functional $\bar{\sigma}_1$ and the norm penalty, we use Clarke's generalized differential \cite{Clarke}. 

\begin{definition}    
Let $E$ be a Banach space and $f \colon E \to \mathbb{R}$ be locally Lipschitz on a neighborhood of $x \in E$. The \emph{Clarke directional derivative} of $f$ at $x$ in the direction $v \in E$ is defined as:
\[
    f^\circ(x; v) := \limsup_{\tilde{x} \to x, \, t \downarrow 0} \frac{f(\tilde{x} + tv) - f(\tilde{x})}{t}.
\]
The \emph{Clarke subdifferential} $\partial_C f(x)$ is the non-empty, convex, weak-$*$ compact set of continuous linear functionals given by:
\[
    \partial_C f(x) := \left\{ \xi \in E^* \mid \langle \xi, v \rangle \le f^\circ(x; v) \quad \forall v \in E \right\}.
\]
\end{definition}

To differentiate the Ekeland distance penalty, we use the subdifferential of the norm.
\begin{proposition}[{\cite[Proposition 4.6.2]{Sch}}] \label{prop:subdiff_norm}
Let $E$ be a Banach space, and $\omega_z(x) := \|x - z\|_E$. Then
\[
    \partial_C \omega_z(x) = \left\{ x^* \in E^* \mid \|x^*\|_{E^*} = 1, \; \langle x - z, x^* \rangle = \|x - z\|_E \right\} \quad \text{if } x \neq z,
\]
and $\partial_C \omega_z(z) = B_{E^*}$, where $B_{E^*}$ is the closed unit ball in the dual space $E^*$.
\end{proposition}

To differentiate composite functionals involving the Steklov eigenvalues, we use the chain rule:
\begin{proposition}[Chain Rule, {\cite[Proposition 2.2]{Vin1}}] \label{prop:subdiff_chain}
Let $h \colon F \to E$ be a continuously Fr\'echet differentiable map of Banach spaces defined on a neighborhood of $x \in F$. If $f \colon E \to \mathbb{R}$ is Lipschitz on a neighborhood of $y = h(x) \in E$, then
\[
    \partial_C (f \circ h)(x) \subset d_x h^* \left[ \partial_C f(y) \right].
\]
\end{proposition}

When $f = \min(f_1, f_2)$ for continuously differentiable functions $f_1, f_2$, the subdifferential at a crossing point $f_1(x) = f_2(x)$ is the convex hull of the individual gradients: $\partial_C f(x) = \mathrm{conv}\{ \nabla f_1(x), \nabla f_2(x) \}$. 

Finally, since our optimization is posed on a closed convex subset, we use Clarke's multiplier rule:
\begin{proposition}[Clarke's Multiplier Rule, {\cite[Theorem 12.4.1]{Sch}}] \label{prop:clarke_multiplier}
Let $E$ be a Banach space and $f \colon E \to \mathbb{R}$ be locally Lipschitz on a neighborhood of $x \in S$, where $S \subset E$ is a closed convex set. If $x$ is a local minimum of $f$ on $S$, then there exists $\xi \in \partial_C f(x)$ such that:
\[
    \langle \tilde{x} - x, \xi \rangle \ge 0 \quad \forall \tilde{x} \in S.
\]
\end{proposition}

\subsubsection*{Subdifferential of the Eigenvalue Operator}

Consider the Steklov eigenvalue $\sigma_k(\alpha, \mu)$ on $\Omega_L$ governed by a bulk metric density $\alpha \in L^q(\Omega_L)$ and a boundary measure $\mu \in \mathcal{M}_+(\partial \Omega_L)$. 

By standard perturbation theory for isolated eigenvalues of finite multiplicity (see \cite{FS12, Vin1, Vin2}), the map $(\alpha, \mu) \mapsto \sigma_k(\alpha, \mu)$ is locally Lipschitz. When $\sigma_k$ is simple, it is Fr\'echet differentiable, and for a normalized eigenfunction $u$, its derivative is given by:
\[
    \nabla_\alpha \sigma_k = |\nabla u|_{g_0}^2, \quad \nabla_\mu \sigma_k = -\sigma_k u^2 \big|_{\partial \Omega_L}.
\]
When $\sigma_k$ has multiplicity $m \ge 2$, the Clarke subdifferential is contained in the convex hull of the derivatives generated by the eigenspace $E_{\sigma_k}$:
\begin{equation}\label{eq:clarke_steklov}
    \partial_C \sigma_k(\alpha, \mu) \subset \mathrm{conv} \left\{ \left( |\nabla u|_{g_0}^2, -\sigma_k u^2 \big|_{\partial \Omega_L} \right) \;\middle|\; u \in E_{\sigma_k}, \int_{\partial \Omega_L} u^2 \, d\mu = 1 \right\}.
\end{equation}

In our $O(n)$-invariant setting on $\Omega_L = [-L, L] \times \Sbb^{n-1}$, the first non-zero eigenvalue is $\sigma_1 = \min\{\sigma^{(0)}_1, \sigma^{(1)}_1\}$, where $\sigma^{(0)}_1$ and $\sigma^{(1)}_1$ are each simple at the one-dimensional profile level. Consequently, if the two branches cross ($\sigma^{(0)}_1 = \sigma^{(1)}_1$), there exists $\nu \in [0, 1]$ such that any $\xi \in \partial_C \sigma_1$ is of the form
\[
    \xi = \nu \nabla \sigma^{(0)}_1 + (1-\nu) \nabla \sigma^{(1)}_1.
\]

\subsection{The Complete Metric Space of Admissible Parameters}
\label{sec:admissible_parameters} 

To apply Ekeland's Variational Principle, we introduce a family of complete metric spaces exhausting $\mathcal{M}_L$. Recall $q=\frac{n}{n-2}$. For each regularization parameter $\epsilon > 0$, we define:
\begin{equation}
    X_\epsilon := \left\{ (\alpha, \mu) \in L^q([-L, L]) \times \mathbb{R}_{\ge 0}^2 \;\middle|\; \alpha(t) \ge \epsilon \text{ a.e.} \right\}.
\end{equation}
Notice that $\mathcal{M}_L = \bigcup_{\epsilon > 0} \left\{ (\alpha, \mu) \in X_\epsilon \mid \mu_{-L} + \mu_L > 0 \right\}$.

\begin{proposition} \label{prop:completeness_X_epsilon}
For every $\epsilon > 0$, the space $X_\epsilon$ is a complete metric space when equipped with the distance:
\begin{equation}
    d((\alpha_1, \mu_1), (\alpha_2, \mu_2)) = \|\alpha_1 - \alpha_2\|_{L^q([-L, L])} + \sum_{\pm L} |\mu_{1, \pm L} - \mu_{2, \pm L}|.
\end{equation}
\end{proposition}

\begin{proof}
The ambient space $L^q([-L, L]) \times \mathbb{R}^2$ is a Banach space. Since $L^q$-convergence implies pointwise convergence almost everywhere along a subsequence, the condition $\alpha(t) \ge \epsilon$ a.e., together with $\mu_{\pm L} \ge 0$, defines a closed subset of $L^q([-L, L]) \times \mathbb{R}^2$. Hence $(X_\epsilon, d)$ is complete.
\end{proof}

\subsection{Application of Ekeland's Variational Principle}
\label{sec:ekeland_application}

Recall from Definition \ref{def:space_ML} and Theorem \ref{thm:uniform_bound} that
\[
    \Sigma_1(L) = \sup_{(\alpha, \mu) \in \mathcal{M}_L} \bar{\sigma}_1(\alpha, \mu) < \infty.
\]
By definition of $\Sigma_1(L)$ and the scale invariance of $\bar{\sigma}_1$ in $\alpha$ and in $\mu$, we may first choose a sequence of normalized admissible configurations $(\tilde{\alpha}_j, \tilde{\mu}_j) \in \mathcal{M}_L$ satisfying $\|\tilde{\alpha}_j\|_{L^q} = 1$ and $\tilde{\mu}_{j, -L} + \tilde{\mu}_{j, L} = 1$, together with a sequence of errors $\delta_j \downarrow 0$, such that
\[
    \bar{\sigma}_1(\tilde{\alpha}_j, \tilde{\mu}_j) > \Sigma_1(L) - \delta_j.
\]
Since $(\tilde{\alpha}_j, \tilde{\mu}_j) \in \mathcal{M}_L$, we have $\operatorname*{ess\,inf}_{[-L, L]} \tilde{\alpha}_j > 0$ for each $j$. We may therefore select a strictly decreasing sequence of obstacle parameters $0 < \epsilon_j < \operatorname*{ess\,inf}_{[-L, L]} \tilde{\alpha}_j$ such that $\epsilon_j \downarrow 0$. By construction, the normalized pair $(\tilde{\alpha}_j, \tilde{\mu}_j)$ belongs to $X_{\epsilon_j}$ for every $j$. To keep the notation light, we henceforth index this sequence directly by $\epsilon = \epsilon_j \downarrow 0$, writing $(\tilde{\alpha}_\epsilon, \tilde{\mu}_\epsilon) := (\tilde{\alpha}_j, \tilde{\mu}_j) \in X_\epsilon$ and $\delta_\epsilon := \delta_j \to 0$.

Define the functional $F \colon X_\epsilon \to \mathbb{R} \cup \{+\infty\}$ by $F(\alpha, \mu) := -\bar{\sigma}_1(\alpha, \mu)$ if $\mu_{-L} + \mu_L \ge 1/2$, and $F(\alpha, \mu) := +\infty$ otherwise. Then $F$ is bounded from below by $-\Sigma_1(L)$ and is lower semi-continuous on $(X_\epsilon, d)$ by continuity of $\sigma_1$ on $X_\epsilon \cap \{\mu_{-L} + \mu_L \ge 1/2\}$. Setting the displacement parameter
\[
    \gamma_\epsilon := \sqrt{\delta_\epsilon} \to 0,
\]
the inequality $\bar{\sigma}_1(\tilde{\alpha}_\epsilon, \tilde{\mu}_\epsilon) > \Sigma_1(L) - \delta_\epsilon$ implies $F(\tilde{\alpha}_\epsilon, \tilde{\mu}_\epsilon) < -\Sigma_1(L) + \delta_\epsilon \le \inf_{X_\epsilon} F + \delta_\epsilon$. Ekeland's Variational Principle applied to $F$ on $(X_\epsilon, d)$ with parameters $\delta_\epsilon$ and $\gamma_\epsilon$ (so that $\delta_\epsilon / \gamma_\epsilon = \gamma_\epsilon$) yields a pair $(\alpha_\epsilon, \mu_\epsilon) \in X_\epsilon$ satisfying:
\begin{enumerate}
    \item \textbf{Energy improvement:} 
    \[
        \bar{\sigma}_1(\alpha_\epsilon, \mu_\epsilon) \ge \bar{\sigma}_1(\tilde{\alpha}_\epsilon, \tilde{\mu}_\epsilon) > \Sigma_1(L) - \delta_\epsilon.
    \]

    \item \textbf{Proximity to the normalized sequence:} 
    \[
        \|\alpha_\epsilon - \tilde{\alpha}_\epsilon\|_{L^q} + \sum_{\pm L} |\mu_{\epsilon, \pm L} - \tilde{\mu}_{\epsilon, \pm L}| \le \gamma_\epsilon.
    \]

    \item \textbf{Strict minimality for the perturbed functional:} The pair $(\alpha_\epsilon, \mu_\epsilon)$ is the strict global minimizer on $X_\epsilon$ of
    \[
        F_\epsilon(\alpha, \mu) := F(\alpha, \mu) + \gamma_\epsilon \left( \|\alpha - \alpha_\epsilon\|_{L^q} + \sum_{\pm L} |\mu_{\pm L} - \mu_{\epsilon, \pm L}| \right).
    \]
\end{enumerate}
Since $(\alpha_\epsilon, \mu_\epsilon)$ is a strict global minimizer for $F_\epsilon$, zero must belong to its Clarke subdifferential plus the normal cone evaluated at this point. This critical algebraic condition allows us to extract the approximate $n$-harmonic map equation governing the conformal factor.

\subsubsection*{Uniform Boundary Non-Collapse}
\label{sec:boundary_noncollapse}

To ensure the strong limit of the sequence represents a valid Steklov problem on the cylinder, we must rigorously preclude boundary collapse ($\mu_{\epsilon, \pm L} \to 0$).

\begin{lemma}[Uniform Boundary Non-Collapse] \label{lem:boundary_noncollapse}
Let $L > 0$ be a finite conformal half-length such that $\Sigma_1(L) > \bar{\sigma}_1(\B^n)(1 - e^{-2nL})^{\frac{2}{n}}$. There exists a constant $C_L > 0$, independent of $\epsilon$, such that for $\epsilon$ sufficiently small, the boundary measures satisfy $\mu_{\epsilon, \pm L} \ge C_L$.

In particular, this condition holds for any $L$ in the stable regime $\Sigma_1(L) > \bar{\sigma}_1(\B^n)$.
\end{lemma}

\begin{proof}
By contradiction, suppose up to a subsequence that the left boundary measure collapses: $\mu_{\epsilon, -L} \to 0$. Since $\mu_{\epsilon, L} + \mu_{\epsilon, -L} = (\mu_{\epsilon, L}-\tilde\mu_{\epsilon, L}) + (\mu_{\epsilon, -L}-\tilde\mu_{\epsilon, -L}) + 1 \to 1$, the mass concentrates entirely on the right boundary: $\mu_{\epsilon, L} \to 1$. We conformally map the cylinder $\Omega_L = [-L, L] \times \Sbb^{n-1}$ to a Euclidean annulus $A_L \subset \R^n$ via the radial coordinate transformation $x(t, \xi) = e^{t-L}\xi$, where $\xi \in \Sbb^{n-1}$. Under this diffeomorphism, the right boundary $\Gamma_L$ ($t=L$) maps to the outer unit sphere $|x|=1$, and the left boundary $\Gamma_{-L}$ ($t=-L$) maps to the inner sphere $|x|=e^{-2L}$. The domain becomes the annulus $A_L = \{x \in \R^n \mid e^{-2L} \le |x| \le 1\}$. Let us explicitly compute the pushforward of the conformal density. The background cylindrical metric is $g_0 = dt^2 + g_{\Sbb^{n-1}}$. Since $|x| = e^{t-L}$, we have $d|x| = e^{t-L}dt = |x|dt$, which implies $dt^2 = |x|^{-2} d|x|^2$. Consequently, the Euclidean metric is given by $dx^2 = d|x|^2 + |x|^2 g_{\Sbb^{n-1}} = |x|^2 (dt^2 + g_{\Sbb^{n-1}})$. The invariant metric $g = \alpha_\epsilon(t)^{\frac{2}{n-2}} g_0$ on the cylinder can thus be expressed on the annulus as:
\[
    g = \alpha_\epsilon(t)^{\frac{2}{n-2}} |x|^{-2} dx^2 = \left(\alpha_\epsilon(t) |x|^{2-n}\right)^{\frac{2}{n-2}} dx^2.
\]

Therefore, the conformal density in Euclidean coordinates is precisely $\tilde{\alpha}_\epsilon(x) = \alpha_\epsilon(t) |x|^{2-n}$. 
The Euclidean volume element is $dx = |x|^{n-1} d|x| d\xi = |x|^n dt d\xi$. The Lebesgue $L^q$-norm transforms as:
\[
    \int_{A_L} \tilde{\alpha}_\epsilon^q \, dx = \int_{-L}^L \int_{\Sbb^{n-1}} \left(\alpha_\epsilon(t) |x|^{2-n}\right)^{\frac{n}{n-2}} |x|^n \, d\xi \, dt = \vol(\Sbb^{n-1}) \int_{-L}^L \alpha_\epsilon(t)^{\frac{n}{n-2}} \, dt,
\]
so that $\|\tilde{\alpha}_\epsilon\|_{L^q(A_L)} \to \vol(\Sbb^{n-1})^{\frac{n-2}{n}}$ since $\|\alpha_\epsilon\|_{L^q([-L,L])} \to 1$.

We bound the first Steklov eigenvalue $\sigma_\epsilon$ using the min-max principle. Consider the restricted Euclidean coordinate functions $u_i(x) = x_i$ for $i = 1, \dots, n$. To be admissible test functions, they must have zero mean on the boundary with respect to the sequence's boundary measure $d\mu_\epsilon = \mu_{\epsilon, L} d\xi|_{\Gamma_L} + \mu_{\epsilon, -L} d\xi|_{\Gamma_{-L}}$. Because $u_i = \xi_i$ on $\Gamma_L$ and $u_i = e^{-2L}\xi_i$ on $\Gamma_{-L}$, and the standard coordinate functions integrate to zero over the sphere ($\int_{\Sbb^{n-1}} \xi_i \, d\xi = 0$), we obtain $\int_{\partial \Omega_L} u_i \, d\mu_\epsilon = 0$. Thus, each $u_i$ is admissible. By the variational characterization of $\sigma_\epsilon$, for each $i$, we have $\sigma_\epsilon \int_{\partial \Omega_L} u_i^2 \, d\mu_\epsilon \le \int_{A_L} \tilde{\alpha}_\epsilon |\nabla u_i|^2 \, dx$. Summing these inequalities over all $i = 1, \dots, n$ yields:
\[
    \sigma_\epsilon \le \frac{\sum_{i=1}^n \int_{A_L} \tilde{\alpha}_\epsilon |\nabla u_i|^2 \, dx}{\sum_{i=1}^n \int_{\partial \Omega_L} u_i^2 \, d\mu_\epsilon}.
\]
For the numerator, since $\sum_{i=1}^n |\nabla x_i|^2 = n$, the total Dirichlet energy is precisely $n \int_{A_L} \tilde{\alpha}_\epsilon \, dx$. Applying H\"older's inequality with conjugate exponents $q = \frac{n}{n-2}$ and $p = \frac{n}{2}$, we obtain:
\[
    \sum_{i=1}^n \int_{A_L} \tilde{\alpha}_\epsilon |\nabla u_i|^2 \, dx \le n \|\tilde{\alpha}_\epsilon\|_{L^q(A_L)} \vol(A_L)^{\frac{2}{n}}.
\]
For the denominator, evaluating $\sum_{i=1}^n x_i^2 = |x|^2$ on the boundaries yields:
\[
    \sum_{i=1}^n \int_{\partial \Omega_L} u_i^2 \, d\mu_\epsilon = \vol(\Sbb^{n-1}) \left(\mu_{\epsilon, L}(1)^2 + \mu_{\epsilon, -L}(e^{-2L})^2\right).
\]
Because $\mu_{\epsilon, L} \to 1$ and $\mu_{\epsilon, -L} \to 0$, this denominator converges to $\vol(\Sbb^{n-1})$. Taking the limit $\epsilon \to 0$ of the ratio gives:
\[
    \limsup_{\epsilon \to 0} \sigma_\epsilon \le \frac{n \vol(\Sbb^{n-1})^{\frac{n-2}{n}} \vol(A_L)^{\frac{2}{n}}}{\vol(\Sbb^{n-1})} = n \vol(\Sbb^{n-1})^{-\frac{2}{n}} \vol(A_L)^{\frac{2}{n}}.
\]
Multiplying by $(\mu_{\epsilon, L} + \mu_{\epsilon, -L})\|\alpha_\epsilon\|_{L^q}^{-1}\vol(\Sbb^{n-1})^{\frac{2}{n}}$ to recover the scale-invariant normalized functional $\bar{\sigma}_1$ cancels the spherical volume terms:
\[
    \limsup_{\epsilon \to 0} \bar{\sigma}_1(\alpha_\epsilon, \mu_\epsilon) \le n \vol(A_L)^{\frac{2}{n}}.
\]

We now express the volume of the annulus in terms of the standard ball. Since $A_L = \B^n \setminus (e^{-2L}\B^n)$, its volume is $\vol(A_L) = \vol(\B^n)(1 - e^{-2nL})$. Substituting this relation into the upper bound gives:
\[
    \limsup_{\epsilon \to 0} \bar{\sigma}_1(\alpha_\epsilon, \mu_\epsilon) \le n \vol(\B^n)^{\frac{2}{n}} (1 - e^{-2nL})^{\frac{2}{n}}.
\]
As detailed in Proposition \ref{prop:ball_supremum}, the normalized Steklov eigenvalue of the standard Euclidean ball is $\bar{\sigma}_1(\B^n) = n \vol(\B^n)^{\frac{2}{n}}$. Therefore, the limit is bounded by:
\[
    \limsup_{\epsilon \to 0} \bar{\sigma}_1(\alpha_\epsilon, \mu_\epsilon) \le \bar{\sigma}_1(\B^n)(1 - e^{-2nL})^{\frac{2}{n}}.
\]

However, our Ekeland sequence satisfies $\lim_{\epsilon \to 0} \bar{\sigma}_1(\alpha_\epsilon, \mu_\epsilon) = \Sigma_1(L)$. Furthermore, by hypothesis, we restricted our lemma to lengths $L$ satisfying $\Sigma_1(L) > \bar{\sigma}_1(\B^n)(1 - e^{-2nL})^{\frac{2}{n}}$. This immediately produces the contradiction:
\[
    \Sigma_1(L) \le \bar{\sigma}_1(\B^n)(1 - e^{-2nL})^{\frac{2}{n}} < \Sigma_1(L).
\]

Therefore, the assumption must be false, and the left boundary measure cannot collapse. By a symmetric conformal mapping (using $x = e^{-t-L}\xi$) applied to the right boundary, an identical contradiction holds if $\mu_{\epsilon, L} \to 0$. Thus, there exists a uniform constant $C_L > 0$ such that $\mu_{\epsilon, \pm L} \ge C_L$ for all sufficiently small $\epsilon$, equivalently $\liminf_{\epsilon \to 0}\mu_{\epsilon, \pm L}>0$.
\end{proof}
From now on and for the remainder of this section, we assume that $\Sigma_1(L) > \bar{\sigma}_1(\B^n)$. 

\subsection{Clarke Subdifferential and Approximate Criticality}
\label{sec:approximate_criticality}

\begin{lemma}[Approximate Equivariant Eigenmap and $\epsilon$-Euler--Lagrange Conditions] \label{lem:approx_criticality_conditions}
Let $(\alpha_\epsilon, \mu_\epsilon) \in X_\epsilon$ be the almost-critical maximizing sequence obtained by Ekeland's Variational Principle in Section \ref{sec:ekeland_application}. Let $\sigma_\epsilon := \sigma_1(\alpha_\epsilon, \mu_\epsilon) \to \Sigma_1(L) \vol(\Sbb^{n-1})^{-2/n}$. There exists a $W^{1,n}(\Omega_L,\mathbb{R}^{m_\epsilon})$-equivariant map $\Phi_\epsilon \colon \Omega_L \to \R^{m_\epsilon} \subseteq \R^{n+1}$ (with target dimension $m_\epsilon \in \{1, n, n+1\}$) and an associated amplitude profile $U_\epsilon \colon [-L, L] \to \R^2$ satisfying $|\Phi_\epsilon(t, \xi)|^2 = |U_\epsilon(t)|^2$, such that:
\begin{itemize}
    \item[(i)] $\Phi_\epsilon$ is a weak $\sigma_\epsilon$-Steklov eigenmap for $(\alpha_\epsilon, \mu_\epsilon)$, i.e., for all $\Psi\in W^{1,n}(\Omega_L,\mathbb{R}^{m_\epsilon})$: 
    \[
        \int_{\Omega_L} \alpha_\epsilon \langle \nabla \Phi_\epsilon, \nabla \Psi \rangle \, dv_{g_0} = \sigma_\epsilon \sum_{\pm L} \mu_{\epsilon, \pm L} \int_{\Sbb^{n-1}} \langle \Phi_\epsilon(\pm L, \xi), \Psi(\pm L, \xi) \rangle \, d\xi;
    \]

    \item[(ii)] $\left\| C_\epsilon \alpha_\epsilon^{q-1} - e(U_\epsilon) \right\|_{L^{n/2}([-L, L])} \to 0$, where $C_\epsilon := \bar{\sigma}_1(\alpha_\epsilon, \mu_\epsilon) \vol(\Sbb^{n-1})^{-\frac{2}{n}} \|\alpha_\epsilon\|_{L^q}^{1-q}$ and the one-dimensional invariant energy density is $e(U_\epsilon)(t) := |d\Phi_\epsilon|_{g_0}^2$;

    \item[(iii)] The boundary terms satisfy the one-sided bound $|U_\epsilon(\pm L)|^2 \ge 1 + \tilde{\nu}_{\epsilon, \pm L}$ and the identity $\mu_{\epsilon, \pm L} |U_\epsilon(\pm L)|^2 = \mu_{\epsilon, \pm L} \left( 1 + \tilde{\nu}_{\epsilon, \pm L} \right)$, where $\tilde{\nu}_{\epsilon, \pm L} \to 0$.
\end{itemize}
\end{lemma}

\begin{proof} By the preceding boundary non-collapse
lemma, both boundary masses of the Ekeland pair are strictly positive and bounded away from
zero. The normalized spectral branches and the penalized functional are locally Lipschitz in a
neighborhood of this pair, so the constrained Clarke multiplier rule applies. The pair $(\alpha_\epsilon, \mu_\epsilon)$ minimizes the perturbed functional $F_\epsilon$ on the closed convex set $X_\epsilon$. Since the subset $\{\alpha \in L^q \mid \alpha \ge \epsilon \text{ a.e.}\}$ has empty interior in $L^q([-L, L])$, we restrict variations to $Y_\epsilon := (\alpha_\epsilon + L^\infty([-L, L])) \times \R^2$ endowed with the $L^\infty \times \R^2$ norm. Because $\alpha_\epsilon \ge \epsilon$, any perturbation $h \in L^\infty$ with $\|h\|_{L^\infty} < \epsilon/2$ satisfies $\alpha_\epsilon + h \ge \epsilon/2 > 0$. Moreover, Ekeland proximity ensures that the total boundary mass $M(\mu_\epsilon) := \mu_{\epsilon, -L} + \mu_{\epsilon, L} \to 1$, so for $\epsilon$ sufficiently small, $M(\mu_\epsilon) > 1/2$. 

On this $L^\infty \times \R^2$ neighborhood of $(\alpha_\epsilon, \mu_\epsilon)$, the functional $F_\epsilon$ is locally Lipschitz. Proposition \ref{prop:clarke_multiplier} then yields a subgradient $\xi \in \partial_C F_\epsilon(\alpha_\epsilon, \mu_\epsilon)$ such that $\langle (\tilde{\alpha}, \tilde{\mu}) - (\alpha_\epsilon, \mu_\epsilon), \xi \rangle \ge 0$ for all $(\tilde{\alpha}, \tilde{\mu}) \in Y_\epsilon \cap X_\epsilon$. By the sum rule \cite[Section 2.3]{Vin2}, we have $\xi = -G + \gamma_\epsilon(\eta_\epsilon, \nu_\epsilon)$, where $G \in \partial_C \bar{\sigma}_1(\alpha_\epsilon, \mu_\epsilon)$, and $\eta_\epsilon \in L^{n/2}([-L, L])$, $\nu_\epsilon = (\nu_{\epsilon, -L}, \nu_{\epsilon, L}) \in \R^2$ (given by Proposition \ref{prop:subdiff_norm}) satisfy $\|\eta_\epsilon\|_{L^{n/2}} \le 1$ and $|\nu_{\epsilon, \pm L}| \le 1$. To express $G$ explicitly, we write $\bar{\sigma}_1$ via the normalization map $h(\alpha, \mu) = \left( \frac{\alpha}{\|\alpha\|_{L^q}}, \frac{\mu}{M(\mu)} \right) =: (\hat{\alpha}, \hat{\mu})$, so that
\[
    \bar{\sigma}_1(\alpha, \mu) = \vol(\Sbb^{n-1})^{\frac{2}{n}} \sigma_1(\hat{\alpha}, \hat{\mu}).
\]
By the chain rule (Proposition \ref{prop:subdiff_chain}), any $G \in \partial_C \bar{\sigma}_1(\alpha, \mu)$ belongs to the pullback set $(d_{(\alpha, \mu)} h)^* \left[ \vol(\Sbb^{n-1})^{\frac{2}{n}} \partial_C \sigma_1(\hat{\alpha}, \hat{\mu}) \right]$. For $M(\mu)>0$, the Fr\'echet differentials of $h$ in directions $v \in L^\infty([-L,L])$ and $w \in \mathbb{R}^2$ are:
\[
    d_\alpha \hat{\alpha}[v] = \frac{v}{\|\alpha\|_{L^q}} - \frac{\alpha}{\|\alpha\|_{L^q}^{q+1}} \int_{-L}^L \alpha^{q-1} v \, dt \quad \text{and} \quad d_\mu \hat{\mu}[w] = \frac{w}{M(\mu)} - \frac{\mu}{M(\mu)^2} \sum_{c \in \{\pm L\}} w_c.
\]

Let $\zeta \in \partial_C \sigma_1(\hat{\alpha}, \hat{\mu})$ be represented by $(e(U), -\sigma_1(\hat{\alpha}, \hat{\mu}) |U(\pm L)|^2)$ as in \eqref{eq:clarke_steklov}, where $U$ is the amplitude profile normalized by $\sum_{\pm L} \hat{\mu}_{\pm L} |U(\pm L)|^2 = 1$. Using the identities $\int_{-L}^L e(U) \alpha \, dt = \|\alpha\|_{L^q} \sigma_1(\hat{\alpha}, \hat{\mu})$ and $\sum_{\pm L} \mu_{\pm L} |U(\pm L)|^2 = M(\mu)$, we compute $G = (d_{(\alpha, \mu)} h)^* [ \vol(\Sbb^{n-1})^{\frac{2}{n}} \zeta ] = (G_\alpha, G_\mu)$:
\[
    G = \left( \frac{\vol(\Sbb^{n-1})^{\frac{2}{n}}}{\|\alpha\|_{L^q}} e(U) - \frac{\bar{\sigma}_1(\alpha, \mu)}{\|\alpha\|_{L^q}^q} \alpha^{q-1}, \quad \frac{\bar{\sigma}_1(\alpha, \mu)}{M(\mu)} \left( 1 - |U(\pm L)|^2 \right) \right).
\]

At $(\alpha_\epsilon, \mu_\epsilon)$, the eigenspace $E_1(\alpha_\epsilon, \mu_\epsilon)$ decomposes into longitudinal modes ($k=0$) and transverse modes ($k=1$). By Theorem \ref{thm:spectrum_computation} and Corollary \ref{cor:multiplicity}, the transverse eigenspace is spanned by $\psi_{1,\epsilon}(t)\xi_i$ for $i=1,\dots,n$. Assuming the 1D profiles $\psi_{0,\epsilon}$ and $\psi_{1,\epsilon}$ are normalized such that $\sum_{\pm L} \mu_{\epsilon, \pm L} \psi_{k,\epsilon}(\pm L)^2 = 1$ (for $k=0,1$), a convex combination of the longitudinal and transverse branches with weights $c_0, c^\perp \ge 0$ satisfying $c_0 + c^\perp = 1$ yields the one-dimensional profile
\[
    U_\epsilon(t) := \sqrt{M(\mu_\epsilon)} \left( \sqrt{c_0}\psi_{0,\epsilon}(t), \sqrt{c^\perp}\psi_{1,\epsilon}(t) \right) =: \left( U_\epsilon^0(t), U_\epsilon^\perp(t) \right).
\]
We associate to $U_\epsilon$ the equivariant map $\Phi_\epsilon(t, \xi) := \left( U_\epsilon^0(t), U_\epsilon^\perp(t)\xi_1, \dots, U_\epsilon^\perp(t)\xi_n \right)$, whose non-zero components are weak Steklov eigenfunctions for $\sigma_\epsilon$, establishing (i). On the cylinder metric $g_0 = dt^2 + g_{\Sbb^{n-1}}$, the energy density $|d\Phi_\epsilon|_{g_0}^2$ is given by:
\[
    e(U_\epsilon)(t) = |\partial_t \Phi_\epsilon|^2 + |\nabla_{\Sbb^{n-1}} \Phi_\epsilon|^2 = |\dot{U}_\epsilon(t)|^2 + (n-1)|U_\epsilon^\perp(t)|^2.
\]

Substituting $G = (G_\alpha, G_\mu)$ at $U_\epsilon$ into the inequality $\langle -G + \gamma_\epsilon(\eta_\epsilon, \nu_\epsilon), (v, \delta\mu) \rangle \ge 0$ for admissible directions $v = \tilde{\alpha} - \alpha_\epsilon \in L^\infty$ and $\delta\mu_{\pm L} = \tilde{\mu}_{\pm L} - \mu_{\epsilon, \pm L}$ gives:
\[
\begin{aligned}
    \int_{-L}^L \left[ -\frac{\vol(\Sbb^{n-1})^{\frac{2}{n}}}{\|\alpha_\epsilon\|_{L^q}} e(U_\epsilon) + \frac{\bar{\sigma}_1(\alpha_\epsilon, \mu_\epsilon)}{\|\alpha_\epsilon\|_{L^q}^q} \alpha_\epsilon^{q-1} + \gamma_\epsilon \eta_\epsilon \right] v \, dt & \\
    + \sum_{\pm L} \left[ \frac{\bar{\sigma}_1(\alpha_\epsilon, \mu_\epsilon)}{M(\mu_\epsilon)} \left( |U_\epsilon(\pm L)|^2 - 1 \right) + \gamma_\epsilon \nu_{\epsilon, \pm L} \right] \delta\mu_{\pm L} &\ge 0.
\end{aligned}
\]
Multiplying by $\frac{\|\alpha_\epsilon\|_{L^q}}{\vol(\Sbb^{n-1})^{2/n}}$ and setting
\[
    C_\epsilon := \bar{\sigma}_1(\alpha_\epsilon, \mu_\epsilon) \vol(\Sbb^{n-1})^{-\frac{2}{n}} \|\alpha_\epsilon\|_{L^q}^{1-q}, \qquad B_\epsilon := \frac{\|\alpha_\epsilon\|_{L^q} \bar{\sigma}_1(\alpha_\epsilon, \mu_\epsilon)}{\vol(\Sbb^{n-1})^{2/n} M(\mu_\epsilon)} > 0,
\]
we obtain the coupled variational inequality:
\begin{equation}
\label{eq:full_variational_inequality}
    \int_{-L}^L \left( C_\epsilon \alpha_\epsilon^{q-1} - e(U_\epsilon) + \tilde{\eta}_\epsilon \right) v \, dt + B_\epsilon \sum_{\pm L} \left( |U_\epsilon(\pm L)|^2 - 1 - \tilde{\nu}_{\epsilon, \pm L} \right) \delta\mu_{\pm L} \ge 0,
\end{equation}
where $\tilde{\eta}_\epsilon := \frac{\|\alpha_\epsilon\|_{L^q} \gamma_\epsilon}{\vol(\Sbb^{n-1})^{2/n}} \eta_\epsilon \in L^{n/2}$ and $\tilde{\nu}_{\epsilon, \pm L} := - \frac{M(\mu_\epsilon) \gamma_\epsilon}{\bar{\sigma}_1(\alpha_\epsilon, \mu_\epsilon)} \nu_{\epsilon, \pm L} \in \mathbb{R}$ both converge to $0$ as $\epsilon \to 0$.

First, setting $v = 0$ and choosing arbitrary $\delta\mu_{\pm L} \ge 0$ in \eqref{eq:full_variational_inequality} yields
\begin{equation} \label{eq:onesided_bound}
    |U_\epsilon(\pm L)|^2 \ge 1 + \tilde{\nu}_{\epsilon, \pm L}.
\end{equation}
Next, testing radial scalings $\tilde{\mu} = (1 \pm t)\mu_\epsilon$ for $t \in (0, 1)$, i.e., $\delta\mu_{\pm L} = \pm t \mu_{\epsilon, \pm L}$ with $v = 0$, gives
\[
    \sum_{\pm L} \mu_{\epsilon, \pm L} \left( |U_\epsilon(\pm L)|^2 - 1 - \tilde{\nu}_{\epsilon, \pm L} \right) = 0.
\]
Since $\mu_{\epsilon, \pm L} \ge 0$ and each parenthesis is non-negative by \eqref{eq:onesided_bound}, each term vanishes separately, proving (iii):
\begin{equation} \label{eq:euler_identity}
    \mu_{\epsilon, \pm L} |U_\epsilon(\pm L)|^2 = \mu_{\epsilon, \pm L} \left( 1 + \tilde{\nu}_{\epsilon, \pm L} \right).
\end{equation}

Finally, we fix $\delta\mu = 0$ and vary the bulk density. Choosing any $v \in L^\infty_+([-L,L])$ gives $\tilde{\alpha} = \alpha_\epsilon + v \ge \epsilon$, which forces
\[
    C_\epsilon \alpha_\epsilon^{q-1} - e(U_\epsilon) + \tilde{\eta}_\epsilon \ge 0 \quad \text{a.e. on } [-L, L].
\]
For the reverse inequality, testing $v = \max\{\epsilon - \alpha_\epsilon, -K\}$ for $K > 0$ (so that $\tilde{\alpha} = \max\{\epsilon, \alpha_\epsilon - K\} \ge \epsilon$) shows that $C_\epsilon \alpha_\epsilon^{q-1} - e(U_\epsilon) + \tilde{\eta}_\epsilon \le 0$ a.e. on $\{\alpha_\epsilon > \epsilon\}$. On the contact set $\{\alpha_\epsilon = \epsilon\}$, since $e(U_\epsilon) \ge 0$, we have $C_\epsilon \alpha_\epsilon^{q-1} - e(U_\epsilon) \le C_\epsilon \epsilon^{q-1}$. Combining both regions yields pointwise a.e.:
\[
    |C_\epsilon \alpha_\epsilon^{q-1} - e(U_\epsilon)| \le |\tilde{\eta}_\epsilon(t)| + C_\epsilon \epsilon^{q-1}.
\]
Taking the $L^{n/2}([-L, L])$ norm gives
\[
    \left\| C_\epsilon \alpha_\epsilon^{q-1} - e(U_\epsilon) \right\|_{L^{n/2}} \le \|\tilde{\eta}_\epsilon\|_{L^{n/2}} + C_\epsilon \epsilon^{q-1} (2L)^{2/n} \longrightarrow 0 \quad \text{as } \epsilon \to 0,
\]
which proves (ii).
\end{proof}

\subsection{A Priori Uniform Estimates}
\label{sec:a_priori_estimates}

By restricting to a $1$-dimensional profile, the conformal metric density $\alpha$ couples linearly to the invariant energy density, reducing the $n$-dimensional Dirichlet energy to the weighted integral $\int_{-L}^L \alpha(t) e(U)(t) \, dt$. Since $n \ge 3$, Morrey's inequality yields a continuous embedding of the 1D Sobolev space $W^{1,n}([-L, L])$ into the H\"older space $C^{0, 1-\frac{1}{n}}([-L, L])$, and hence a compact embedding into $C^0([-L, L])$ by the Arzel\`a--Ascoli theorem. This subcritical one-dimensional embedding is a key analytical advantage of the equivariant framework, as it precludes the bubbling phenomena inherent to unconstrained $n$-harmonic map limits.

\begin{corollary}[Exact Amplitude Bounds] \label{cor:exact_amplitude_bounds}
The boundary amplitude of the $\epsilon$-eigenmap satisfies $|U_\epsilon(\pm L)|^2 = 1 + \tilde{\nu}_{\epsilon, \pm L} \to 1$. Furthermore, we have the global bound $\|U_\epsilon\|_{L^\infty([-L, L])} \le C_{2,\epsilon} \to 1$.
\end{corollary}

\begin{proof}
Because Lemma \ref{lem:boundary_noncollapse} guarantees that the boundary measures are strictly positive ($\mu_{\epsilon, \pm L} > 0$) for small enough $\epsilon$, we may divide the boundary identities \eqref{eq:euler_identity} established in Lemma \ref{lem:approx_criticality_conditions} by $\mu_{\epsilon, \pm L}$ to obtain $|U_\epsilon(\pm L)|^2 = 1 + \tilde{\nu}_{\epsilon, \pm L}$. Since the Ekeland errors $\tilde{\nu}_{\epsilon, \pm L} \to 0$, the boundary amplitude converges to $1$. By Lemma \ref{lem:approx_criticality_conditions}(i), $\Phi_\epsilon$ is a weak Steklov eigenmap satisfying 
\[
    -\operatorname{div}(\alpha_\epsilon \nabla \Phi_\epsilon) = \sigma_\epsilon \mu_\epsilon \Phi_\epsilon.
\]
Because the boundary measure $\mu_\epsilon$ is supported on $\partial \Omega_L$, the components of $\Phi_\epsilon$ are $\alpha_\epsilon$-harmonic in the interior bulk $(-L, L)$, satisfying $\operatorname{div}(\alpha_\epsilon \nabla \Phi_\epsilon) = 0$. Define $W_\epsilon(t) := |\Phi_\epsilon(t, \xi)|^2 = |U_\epsilon(t)|^2$. Applying the chain rule for the inner product and multiplying by $\alpha_\epsilon$, the divergence evaluates to:
\[
    \operatorname{div}\left( \alpha_\epsilon \nabla W_\epsilon \right) = 2 \langle \operatorname{div}(\alpha_\epsilon \nabla \Phi_\epsilon), \Phi_\epsilon \rangle + 2 \alpha_\epsilon |\nabla \Phi_\epsilon|^2.
\]
The first term vanishes identically by interior harmonicity. Since $|\nabla \Phi_\epsilon|^2 = e(U_\epsilon)(t)$, and the gradient of $W_\epsilon(t)$ is purely longitudinal, this reduces to the 1-dimensional subharmonic identity:
\[
    \frac{d}{dt} \left( \alpha_\epsilon(t) \dot{W}_\epsilon(t) \right) = 2 \alpha_\epsilon(t) e(U_\epsilon)(t) \ge 0.
\] 

This inequality implies that the 1-dimensional flux $J_\epsilon(t) := \alpha_\epsilon(t) \dot{W}_\epsilon(t)$ is monotonically non-decreasing on $[-L, L]$. If $W_\epsilon(t)$ were to achieve a strict interior maximum at some $t_0 \in (-L, L)$, its derivative $\dot{W}_\epsilon(t)$ would locally transition from positive to negative. Because $\alpha_\epsilon(t) \ge \epsilon > 0$, the flux $J_\epsilon(t)$ would similarly transition from positive to negative, violating its monotonicity. Therefore, $W_\epsilon(t)$ cannot possess a strict interior maximum and must attain its global supremum on the boundary $\partial \Omega_L$. This guarantees the uniform global bound:
\[
    \|U_\epsilon\|_{L^\infty}^2 = \max_{\pm L} W_\epsilon(\pm L) = 1 + \max_{\pm L} \tilde{\nu}_{\epsilon, \pm L} =: C_{2,\epsilon}^2 \to 1.
\]
\end{proof}

\subsubsection*{Uniform $W^{1,n}$ Bounds on the Amplitude Profile}

To maintain the validity of the $\epsilon$-Euler--Lagrange equations, we do not rescale the conformal factor $\alpha_\epsilon$; instead, we exploit the Ekeland distance estimate to bound the sequence.

\begin{lemma} \label{lem:W1n_bound}
The 1D amplitude profiles $U_\epsilon$ are uniformly bounded in $W^{1,n}([-L, L], \R^2)$, and the bulk conformal densities $\alpha_\epsilon$ are uniformly bounded in $L^q([-L, L])$. Consequently, $U_\epsilon$ is uniformly equicontinuous on $[-L, L]$.
\end{lemma}

\begin{proof}
Ekeland's Variational Principle guarantees the proximity bound:
\[
    \|\alpha_\epsilon - \tilde{\alpha}_\epsilon\|_{L^q} + \sum_{\pm L} |\mu_{\epsilon, \pm L} - \tilde{\mu}_{\epsilon, \pm L}| \le \gamma_\epsilon.
\]
Since $\gamma_\epsilon \to 0$, the triangle inequality enforces $\lim_{\epsilon \to 0} \|\alpha_\epsilon\|_{L^q} = 1$ and $\lim_{\epsilon \to 0} \sum_{\pm L} \mu_{\epsilon, \pm L} = 1$. Thus, $\alpha_\epsilon$ is uniformly bounded in $L^q([-L, L])$. From Lemma \ref{lem:approx_criticality_conditions}(ii), we have the $L^{n/2}$ convergence:
\[
    \left\| C_\epsilon \alpha_\epsilon^{q-1} - e(U_\epsilon) \right\|_{L^{n/2}} \to 0,
\]
where $C_\epsilon \to \Sigma_1(L) \vol(\Sbb^{n-1})^{-\frac{2}{n}}$. Since $q = \frac{n}{n-2}$ implies $\frac{n}{2}(q-1) = q$, raising $\alpha_\epsilon^{q-1}$ to the power $n/2$ integrates to $\|\alpha_\epsilon\|_{L^q}^q$. This uniformly bounds the invariant energy density:
\[
    \limsup_{\epsilon \to 0} \|e(U_\epsilon)\|_{L^{n/2}} \le \Sigma_1(L) \vol(\Sbb^{n-1})^{-\frac{2}{n}}.
\]
Because the energy density satisfies $e(U_\epsilon)(t) \ge |\dot{U}_\epsilon(t)|^2$, the longitudinal derivatives $\dot{U}_\epsilon$ are uniformly bounded in $L^n([-L, L], \R^2)$. Combined with the uniform bound $\|U_\epsilon\|_{L^\infty} \le C_{2,\epsilon} \to 1$ from Corollary \ref{cor:exact_amplitude_bounds}, the full sequence $U_\epsilon$ is uniformly bounded in $W^{1,n}([-L, L], \R^2)$. By Morrey's inequality in dimension one, the continuous embedding $W^{1,n}([-L,L]) \hookrightarrow C^{0, 1-\frac{1}{n}}([-L,L])$ bounds the $C^{0, 1-\frac{1}{n}}$ norm of $U_\epsilon$ uniformly, which guarantees uniform equicontinuity.
\end{proof}

\subsection{Strong Convergence and Exact Criticality}
\label{sec:strong_convergence}

With the uniform analytic bounds established on the 1-dimensional domain, we can now pass to the limit as $\epsilon \to 0$. Because the target geometry is constrained to a 1D profile, the analysis is subcritical, preventing the formation of defect measures and bubbling typical of unconstrained conformal maximization. 

\begin{lemma}[Strong Compactness] \label{lem:strong_compactness}
There exists a subsequence such that the bulk density $\alpha_\epsilon \to \alpha_0$ strongly in $L^q([-L, L])$, the boundary weights $\mu_{\epsilon, \pm L} \to \mu_{0, \pm L} > 0$, and the amplitude profile $U_\epsilon \to U_0$ strongly in $W^{1,n}([-L, L], \R^2)$. 

Furthermore, the limit $(\alpha_0, \mu_0)$ is a configuration with $\|\alpha_0\|_{L^q} = 1$ and $\mu_{0, L} + \mu_{0, -L} = 1$, maximizing the functional: $\bar{\sigma}_1(\alpha_0, \mu_0) = \Sigma_1(L)$.
\end{lemma}

\begin{proof}
Because $\alpha_\epsilon$ and $U_\epsilon$ are uniformly bounded in the reflexive spaces $L^q([-L, L])$ and $W^{1,n}([-L, L], \R^2)$ respectively, Kakutani's theorem guarantees weak limits $\alpha_0 \in L^q$ and $U_0 \in W^{1,n}$. For $n \ge 3$, the Arzel\`a--Ascoli theorem together with the uniform equicontinuity of Lemma \ref{lem:W1n_bound} ensures that $U_\epsilon \to U_0$ uniformly in $C^0([-L, L], \R^2)$. This uniform convergence preserves the boundary evaluations: $|U_0(\pm L)|^2 = \lim_{\epsilon \to 0} |U_\epsilon(\pm L)|^2 = 1$. The boundary weights yield limits $\mu_{0, \pm L} \ge C_L > 0$ satisfying the mass conservation $\mu_{0, L} + \mu_{0, -L} = 1$.

To upgrade the weak convergence in the bulk to strong convergence, we decouple the density from the energy. Let $C_0 := \Sigma_1(L) \vol(\Sbb^{n-1})^{-2/n}$. Testing the weak Steklov equation of Lemma \ref{lem:approx_criticality_conditions}(i) against $\Phi_\epsilon$ itself and dividing by $\vol(\Sbb^{n-1})$, we obtain:
\begin{equation}\label{eq:alpha_eps_U_eps}
    \int_{-L}^L \alpha_\epsilon e(U_\epsilon) \, dt = \sigma_\epsilon \sum_{\pm L} \mu_{\epsilon, \pm L} |U_\epsilon(\pm L)|^2.
\end{equation}
Because $\sigma_\epsilon \to C_0$, $\sum_{\pm L} \mu_{\epsilon, \pm L} \to 1$, and $|U_\epsilon(\pm L)|^2 \to 1$, the right-hand side converges to $C_0$. Thus, $\lim_{\epsilon \to 0} \int_{-L}^L \alpha_\epsilon e(U_\epsilon) \, dt = C_0$.

Furthermore, because $\Phi_0 \in W^{1,n}$ is a valid test function, testing the weak equation against the fixed limit $\Phi_0$ yields:
\begin{equation}\label{eq:alpha_eps_Phi_0}
    \int_{-L}^L \alpha_\epsilon \langle \nabla \Phi_\epsilon, \nabla \Phi_0 \rangle \, dt = \sigma_\epsilon \sum_{\pm L} \mu_{\epsilon, \pm L} \langle U_\epsilon(\pm L), U_0(\pm L) \rangle.
\end{equation}
Because $U_\epsilon \to U_0$ uniformly on $[-L, L]$, we have $U_\epsilon(\pm L) \to U_0(\pm L)$. Coupled with $|U_0(\pm L)|^2 = 1$, the boundary pairing converges to $1$. Thus,
\[
    \lim_{\epsilon \to 0} \int_{-L}^L \alpha_\epsilon \langle \nabla \Phi_\epsilon, \nabla \Phi_0 \rangle \, dt = C_0.
\]

We now use the non-negative expansion $|\nabla \Phi_\epsilon - \nabla \Phi_0|^2 \ge 0$, which rearranges algebraically to $e(U_\epsilon) \ge 2\langle \nabla \Phi_\epsilon, \nabla \Phi_0 \rangle - e(U_0)$. Multiplying by the non-negative density $\alpha_\epsilon$ and integrating over $[-L, L]$ yields:
\[
    \int_{-L}^L \alpha_\epsilon e(U_\epsilon) \, dt \ge 2 \int_{-L}^L \alpha_\epsilon \langle \nabla \Phi_\epsilon, \nabla \Phi_0 \rangle \, dt - \int_{-L}^L \alpha_\epsilon e(U_0) \, dt.
\]
By Equations \eqref{eq:alpha_eps_U_eps} and \eqref{eq:alpha_eps_Phi_0}, taking the limit as $\epsilon \to 0$, the first two integrals tend to $C_0$ and $2C_0$ respectively. For the third integral, because $\alpha_\epsilon \rightharpoonup \alpha_0$ weakly in $L^q$ and the fixed function $e(U_0) \in L^{n/2}$ lies in its dual space, the integral converges to $\int_{-L}^L \alpha_0 e(U_0)\,dt$. This enforces the lower bound:
\begin{equation} \label{eq:energy_lower_bound}
    C_0 \ge 2C_0 - \int_{-L}^L \alpha_0 e(U_0) \, dt \implies \int_{-L}^L \alpha_0 e(U_0) \, dt \ge C_0.
\end{equation}

To bound this from above, we use Lemma \ref{lem:approx_criticality_conditions}(ii): $\left\| C_\epsilon \alpha_\epsilon^{q-1} - e(U_\epsilon) \right\|_{L^{n/2}} \to 0$. Because $\|\alpha_\epsilon\|_{L^q} \to 1$ and $(q-1)\frac{n}{2} = q$, the $L^{n/2}$-norm of the first term evaluates to $C_\epsilon \|\alpha_\epsilon\|_{L^q}^{q-1} \to C_0$. The triangle inequality forces $\lim_{\epsilon \to 0} \|e(U_\epsilon)\|_{L^{n/2}} = C_0$. By the weak lower semicontinuity (WLSC) of the $L^n$-norm on $\dot{U}_\epsilon \rightharpoonup \dot{U}_0$ and the uniform convergence $U_\epsilon^\perp \to U_0^\perp$, one has $\|e(U_0)\|_{L^{n/2}} \le \liminf_{\epsilon \to 0} \|e(U_\epsilon)\|_{L^{n/2}} = C_0$. 

We now apply H\"older's inequality to the limiting configuration $(\alpha_0, U_0)$ with conjugate exponents $q = \frac{n}{n-2}$ and $q' = \frac{n}{2}$, followed by the WLSC bounds of the individual norms:
\begin{align} \label{eq:holder_squeeze}
    C_0 &\le \int_{-L}^L \alpha_0 e(U_0) \, dt \quad \text{(by inequality \eqref{eq:energy_lower_bound})} \nonumber\\ 
    &\le \|\alpha_0\|_{L^q} \|e(U_0)\|_{L^{n/2}} \nonumber \\
    &\le \left( \liminf_{\epsilon \to 0} \|\alpha_\epsilon\|_{L^q} \right) \left( \liminf_{\epsilon \to 0} \|e(U_\epsilon)\|_{L^{n/2}} \right) \nonumber \\
    &= 1 \cdot C_0 = C_0. \nonumber
\end{align}

Because the extreme ends of the inequality chain are identically $C_0$, every intermediate inequality must collapse to equality. This forces two conclusions. First, $\|\alpha_0\|_{L^q} = 1$; no mass is lost in the weak limit. Because $q \in (1, \infty)$, the Lebesgue space $L^q$ is uniformly convex. By the Radon--Riesz property, the combination of weak convergence $\alpha_\epsilon \rightharpoonup \alpha_0$ and norm preservation $\|\alpha_\epsilon\|_{L^q} \to \|\alpha_0\|_{L^q}$ forces strong convergence: $\alpha_\epsilon \to \alpha_0$ strongly in $L^q([-L, L])$. 

Second, the collapse to equality forces $\|e(U_0)\|_{L^{n/2}} = \lim_{\epsilon \to 0} \|e(U_\epsilon)\|_{L^{n/2}} = C_0$, and hence $\|\nabla \Phi_0\|_{L^n} = \lim_{\epsilon \to 0} \|\nabla \Phi_\epsilon\|_{L^n}$. Since $n \ge 3$, the Lebesgue space $L^n$ is uniformly convex. Again by the Radon--Riesz property, weak convergence $\nabla \Phi_\epsilon \rightharpoonup \nabla \Phi_0$ in $L^n$ together with norm preservation forces $\nabla \Phi_\epsilon \to \nabla \Phi_0$ strongly in $L^n$. Therefore, the sequence of profiles $U_\epsilon$ converges strongly to $U_0$ in $W^{1,n}([-L, L], \R^2)$, completing the proof.
\end{proof}

\begin{theorem}[Exact Criticality and Smoothness] \label{thm:exact_criticality}
The strong limits $\alpha_0$ and $U_0$ satisfy the 1-dimensional $n$-harmonic map equations. Moreover, the energy density $e(U_0)$ is strictly positive on $[-L, L]$, the conformal density $\alpha_0$ is $C^\infty$ smooth and strictly positive on $[-L, L]$, and $\Phi_0$ defines a smooth $O(n)$-equivariant free boundary $n$-harmonic map into $\B^{m_0} \subseteq \B^{n+1}$. Consequently, setting $f_0 := \alpha_0^{\frac{1}{n-2}} \in C^\infty_{>0}([-L, L])$ and $\rho_{0, \pm L} := \mu_{0, \pm L} f_0(\pm L)^{1-n} > 0$, the smooth $O(n)$-invariant metric $g_L = f_0(t)^2(dt^2 + g_{\Sbb^{n-1}})$ and weight $\rho_0$ achieve $\Sigma_1(L)$.
\end{theorem}

\begin{proof}
Recall that due to Lemma \ref{lem:approx_criticality_conditions}(ii) and Corollary \ref{cor:exact_amplitude_bounds}, one has:
\[\left\Vert{} \bar{\sigma}_1(\alpha_\epsilon, \mu_\epsilon) \vol(\Sbb^{n-1})^{-\frac{2}{n}} \Vert{}\alpha_\epsilon\Vert{}_{L^q}^{1-q} \alpha_\epsilon^{q-1} - e(U_\epsilon) \right\Vert{}_{L^{n/2}} \to 0 \quad \text{and} \quad \vert{}U_\epsilon(\pm L)\vert{}^2 \to 1.\]
Hence, it follows from Lemma \ref{lem:strong_compactness} that almost everywhere on $[-L, L]$:
\[\Sigma_1(L)\vol(\Sbb^{n-1})^{-\frac{2}{n}} \alpha_0(t)^{q-1} = e(U_0)(t) \quad \text{and} \quad \vert{}U_0(\pm L)\vert{}^2 = 1.\]
Isolating the conformal factor with $q = \frac{n}{n-2}$ gives $\alpha_0 = \Sigma_1(L)^{\frac{2-n}{2}} e(U_0)^{\frac{n-2}{2}}\vol(\Sbb^{n-1})^{\frac{n-2}{n}}$. Now, we pass to the limit in the weak Steklov formulation. For any test function $\Psi \in W^{1,n}(\Omega_L, \R^{m_0})$, we have:
\[\int_{\Omega_L} \alpha_\epsilon \langle \nabla \Phi_\epsilon, \nabla \Psi \rangle \, dv_{g_0} = \sigma_\epsilon \sum_{\pm L} \mu_{\epsilon, \pm L} \int_{\Sbb^{n-1}} \langle \Phi_\epsilon(\pm L, \xi), \Psi(\pm L, \xi) \rangle \, d\xi.\]
Because $\alpha_\epsilon \to \alpha_0$ strongly in $L^q(\Omega_L)$ and $\nabla \Phi_\epsilon \to \nabla \Phi_0$ strongly in $L^n(\Omega_L)$, by H\"older's inequality with $p=\frac{n-1}{n-2}$ and $p'=n-1$, the product $\alpha_\epsilon \nabla \Phi_\epsilon$ converges strongly in the dual space $L^{\frac{n}{n-1}}(\Omega_L)$. Since $\nabla \Psi \in L^n(\Omega_L)$, the bulk integral passes to the limit. On the boundary, the uniform convergence of $\Phi_\epsilon$ ensures the trace evaluations converge. Thus, the limit satisfies:
\[\int_{\Omega_L} \alpha_0 \langle \nabla \Phi_0, \nabla \Psi \rangle \, dv_{g_0} = \Sigma_1(L)\vol(\Sbb^{n-1})^{-\frac{2}{n}} \sum_{\pm L} \mu_{0, \pm L} \int_{\Sbb^{n-1}} \langle \Phi_0(\pm L, \xi), \Psi(\pm L, \xi) \rangle \, d\xi.\]
Substituting the expression for $\alpha_0$ leads to the weak $n$-harmonic map equation:
\begin{equation} \label{eq:nharmonic_eq}
    \int_{\Omega_L} e(U_0)^{\frac{n-2}{2}} \langle \nabla \Phi_0, \nabla \Psi \rangle \, dv_{g_0}= \frac{\Sigma_1(L)^{\frac{n}{2}}}{\vol(\Sbb^{n-1})}  \sum_{\pm L} \mu_{0, \pm L} \int_{\Sbb^{n-1}} \langle \Phi_0(\pm L, \xi), \Psi(\pm L, \xi) \rangle \, d\xi.
\end{equation}
Interior regularity theory for $p$-harmonic maps implies that $\Phi_0$ is of class $C^{1,\alpha}$ \cite{Lin} (in fact, in our one-dimensional setting, the weak ODE \eqref{eq:nharmonic_eq} directly implies that the flux $e(U_0)^{\frac{n-2}{2}}\dot{U}_0$ is absolutely continuous on $[-L, L]$, so $U_0 \in C^1([-L, L], \R^2)$ up to the boundary). Standard elliptic bootstrapping requires the operator to be non-degenerate ($e(U_0) > 0$). Let $\mathcal{O} = \{ t \in [-L, L] \mid e(U_0)(t) > 0 \}$ be the open set where the operator is strictly elliptic. On $\mathcal{O}$, regularity theory guarantees that $U_0$ is $C^\infty$, allowing the bulk equation \eqref{eq:nharmonic_eq} to restrict to the 1-dimensional coupled ODE system:
\begin{align}
    \frac{d}{dt} \left( e(U_0)^{\frac{n-2}{2}} \dot{U}_0^0 \right) &= 0, \label{eq:longitudinal_ode} \\
    \frac{d}{dt} \left( e(U_0)^{\frac{n-2}{2}} \dot{U}_0^\perp \right) - (n-1)e(U_0)^{\frac{n-2}{2}} U_0^\perp &= 0. \label{eq:transverse_ode}
\end{align}

On $\mathcal{O}$, multiplying \eqref{eq:longitudinal_ode} by $\dot{U}_0^0$ and \eqref{eq:transverse_ode} by $\dot{U}_0^\perp$, summing the results, and applying the identity $e(U_0) = |\dot{U}_0^0|^2 + |\dot{U}_0^\perp|^2 + (n-1)|U_0^\perp|^2$ yields a total derivative. Integrating this defines the Hamiltonian of the system, which is constant on each connected component of $\mathcal{O}$:
\begin{equation} \label{eq:hamiltonian}
    H(t) := (n-1) e(U_0)^{\frac{n-2}{2}} \left( |\dot{U}_0^0|^2 + |\dot{U}_0^\perp|^2 - |U_0^\perp|^2 \right) \equiv \text{constant}.
\end{equation}

Because $U_0 \in C^1([-L, L])$, the algebraic expression for $H(t)$ defines a globally continuous function on the entire interval $[-L, L]$. Therefore, $H$ is constant everywhere. Suppose, for contradiction, that the system degenerates at some point $t_0 \in [-L, L]$ such that $e(U_0)(t_0) = 0$. By definition, this forces $\dot{U}_0^0(t_0) = 0$, $U_0^\perp(t_0) = 0$, and $\dot{U}_0^\perp(t_0) = 0$, meaning $H(t_0) = 0$.

Evaluating the longitudinal flux from \eqref{eq:longitudinal_ode} in a neighborhood of $t_0$ forces the constant flux on $\mathcal{O}$ to vanish, meaning $\dot{U}_0^0 \equiv 0$ on $\mathcal{O}$. Dividing out the positive leading term $e(U_0)^{\frac{n-2}{2}}$ from the Hamiltonian then reduces the system to the transverse component: $|\dot{U}_0^\perp|^2 = |U_0^\perp|^2$. This equality extends to $[-L,L]$. By elementary ODE theory, the only solution to this differential equation that patches continuously to the degenerate point $U_0^\perp(t_0) = \dot{U}_0^\perp(t_0) = 0$ is the trivial solution $U_0^\perp(t) \equiv 0$. Because both $\dot{U}_0^0 \equiv 0$ and $U_0^\perp \equiv 0$, the energy density vanishes globally: $e(U_0)(t) \equiv 0$. However, $\Sigma_1(L)\vol(\Sbb^{n-1})^{-\frac{2}{n}} \alpha_0^{q-1} = e(U_0)$ and $\Sigma_1(L) > 0$, which would imply that $\alpha_0\equiv 0$, contradicting Lemma \ref{lem:strong_compactness}.

Thus, $e(U_0)(t) > 0$ everywhere on $[-L, L]$. The operator is uniformly non-degenerate elliptic on $[-L, L]$, and standard 1D bootstrapping implies that $U_0 \in C^\infty([-L, L], \R^2)$ and that $\alpha_0$ is strictly positive and $C^\infty$ smooth on $[-L, L]$. Finally, the boundary identity $|U_0(\pm L)|^2 = 1$ ensures that $\Phi_0(\partial \Omega_L) \subset \partial \B^{m_0}$, while integrating \eqref{eq:nharmonic_eq} by parts yields the classical Steklov boundary condition
\[
    e(U_0)(\pm L)^{\frac{n-2}{2}} \partial_\nu \Phi_0 = \frac{\Sigma_1(L)^{\frac{n}{2}}}{\vol(\Sbb^{n-1})} \mu_{0, \pm L} \Phi_0 \quad \text{on } \Gamma_{\pm L}.
\]
Since $e(U_0)(\pm L) > 0$ and $\mu_{0, \pm L} > 0$, the outward conormal derivative $\partial_\nu \Phi_0$ is positively proportional to the position vector $\Phi_0$ on $\partial \Omega_L$, confirming that $\Phi_0$ meets $\partial \B^{m_0}$ orthogonally.
\end{proof}

\begin{corollary}[Injectivity of the Equivariant Eigenmap] \label{cor:embedding_property}
Whenever the maximizing map $\Phi_0(t, \xi) = (U_0^0(t), U_0^\perp(t)\xi) \colon \Omega_L \to \B^{n+1}$ of Theorem \ref{thm:exact_criticality} contains both longitudinal and transverse components ($U_0^0 \not\equiv 0$ and $U_0^\perp \not\equiv 0$), $\Phi_0$ is a smooth embedding.
\end{corollary}

\begin{proof}
Since $U_0^\perp$ is a non-trivial first transverse eigenprofile (minimizing the 1D Rayleigh quotient), we may assume $U_0^\perp(t) \ge 0$. If $U_0^\perp$ were to vanish at an interior point or at the boundary, the local minimum condition or the Steklov boundary conditions would force $\dot{U}_0^\perp = 0$ at that point. By Cauchy-Lipschitz uniqueness for the regular linear ODE $(\alpha_0 \dot{U}_0^\perp)' = (n-1)\alpha_0 U_0^\perp$, this would imply $U_0^\perp \equiv 0$, a contradiction. Thus $U_0^\perp(t) > 0$ for all $t \in [-L, L]$. Furthermore, by \eqref{eq:longitudinal_ode}, the quantity $e(U_0)(t)^{\frac{n-2}{2}} \dot{U}_0^0(t) \equiv c$ is constant on $[-L, L]$. Since $U_0^0$ is a non-constant longitudinal eigenprofile (orthogonal to constants on $\partial \Omega_L$), we must have $c \neq 0$. Because $e(U_0)(t) > 0$ everywhere on $[-L, L]$ by Theorem \ref{thm:exact_criticality}, it follows that $\dot{U}_0^0(t) \neq 0$ for all $t \in [-L, L]$, so $U_0^0$ is strictly monotone on $[-L, L]$. 

Now if $\Phi_0(t_1, \xi_1) = \Phi_0(t_2, \xi_2)$, the strict monotonicity of $U_0^0$ forces $t_1 = t_2$, and then $U_0^\perp(t_1) > 0$ forces $\xi_1 = \xi_2$. Thus $\Phi_0$ is injective. Similarly, its differential $d\Phi_0$ has full rank $n$ everywhere because $\dot{U}_0^0(t) \neq 0$ and $U_0^\perp(t) > 0$. Therefore, $\Phi_0$ is a smooth embedding.
\end{proof}

\section{Continuity and Asymptotics over the Moduli Space}
\label{sec:moduli_space}

In Section \ref{sec:existence}, we established the existence of a smooth maximizing metric within the $O(n)$-invariant conformal class for a standard cylinder of a \emph{fixed} conformal half-length $L$ satisfying $\Sigma_1(L) > \bar{\sigma}_1(\B^n)$. To complete the proof of the Main Theorem, we must now elevate this fixed-length result to a global optimization over the entire moduli space of finite conformal cylinders, parameterized continuously by $L \in (0, \infty)$. Let us recall that the fixed-length supremum map is given by:
\begin{equation}
    \Sigma_1(L) := \sup \left\{ \bar{\sigma}_1(\alpha, \mu) \;\middle|\; (\alpha, \mu) \in \mathcal{M}_L \right\},
\end{equation}
where $\mathcal{M}_L$ denotes the space of admissible $O(n)$-invariant bulk densities and boundary measures on the 1-dimensional domain $[-L, L]$. 

To demonstrate that the global Steklov supremum is achieved at a finite, non-zero conformal modulus $L^* \in (0, \infty)$, we must analyze the continuity and asymptotic behavior of the Steklov functional at the extreme boundaries of the moduli space. 

Specifically, we will show that the map $\Sigma_1(L)$ is continuous on $(0, \infty)$, vanishes as the cylinder collapses ($L \to 0$), and is bounded above by the normalized Steklov eigenvalue of the standard $n$-ball as the cylinder stretches to infinity ($L \to \infty$). Combined with the explicit geometry of the critical $n$-catenoid in Section \ref{sec:spectral_domination}, which strictly dominates the standard ball, these asymptotic bounds force the global maximum to be achieved in the interior of the moduli space.

\subsection{Continuity of the Supremum Map}

\begin{lemma} \label{lem:continuity_lambda}
The map $L \mapsto \Sigma_1(L)$ is continuous on $(0, \infty)$.
\end{lemma}

\begin{proof}
Let $L > 0$ and consider a perturbed length $L' = L + \delta$ for small $\delta$. There is a canonical diffeomorphism between the 1-dimensional domains $[-L, L]$ and $[-L', L']$ given by the uniform rescaling $F(t) = \rho t$, where $\rho = \frac{L'}{L} = 1 + \frac{\delta}{L}$. 

For any invariant metric density $\alpha \in L^q([-L', L'])$ and boundary measure $\mu$, pulling back via $F$ yields a configuration $(\tilde{\alpha}, \tilde{\mu})$ on $[-L, L]$ defined by $\tilde{\alpha}(t) = \alpha(\rho t)$ and $\tilde{\mu}_{\pm L} = \mu_{\pm L'}$. 

We evaluate the Steklov eigenvalue using the Rayleigh quotient. For any test profile $\psi$ on $[-L', L']$, define the pulled-back profile $\tilde{\psi}(t) = \psi(\rho t)$ on $[-L, L]$. Applying the change of variables $s = \rho t$, the energy transforms as:
\begin{align*}
    \int_{-L'}^{L'} \left( |\dot{\psi}(s)|^2 + \lambda_k |\psi(s)|^2 \right) \alpha(s) \, ds &= \int_{-L}^L \left( \rho^{-2} |\dot{\tilde{\psi}}(t)|^2 + \lambda_k |\tilde{\psi}(t)|^2 \right) \tilde{\alpha}(t) \, \rho \, dt \\
    &\le \max(\rho^{-1}, \rho) \int_{-L}^L \left( |\dot{\tilde{\psi}}(t)|^2 + \lambda_k |\tilde{\psi}(t)|^2 \right) \tilde{\alpha}(t) \, dt.
\end{align*}
Since the boundary masses are preserved ($\mu_{\pm L'} = \tilde{\mu}_{\pm L}$), this energy bound transfers directly to the infimum of the Rayleigh quotient, yielding the eigenvalue bound $\sigma_1(\alpha, \mu) \le \max(\rho^{-1}, \rho) \, \sigma_1(\tilde{\alpha}, \tilde{\mu})$.

Simultaneously, the $L^q$ norm of the metric density scales as:
\[
    \|\alpha\|_{L^q([-L', L'])} = \left( \int_{-L}^L \tilde{\alpha}(t)^q \, \rho \, dt \right)^{1/q} = \rho^{1/q} \|\tilde{\alpha}\|_{L^q([-L, L])}.
\]

Dividing the eigenvalue by the $L^q$ norm, the normalized functional satisfies:
\[
    \bar{\sigma}_1(\alpha, \mu) \le \frac{\max(\rho^{-1}, \rho)}{\rho^{1/q}} \bar{\sigma}_1(\tilde{\alpha}, \tilde{\mu}).
\]
Because $\rho \to 1$ as $\delta \to 0$, the multiplicative prefactor is bounded by $1 + C|\delta|$ for some constant $C$ depending only on $L$ and $q$. Taking the supremum over all configurations on $L'$ yields $\Sigma_1(L') \le (1 + C|\delta|) \Sigma_1(L)$. 

By symmetric reasoning pulling back from $L$ to $L'$, the reverse inequality $\Sigma_1(L) \le (1 + C'|\delta|) \Sigma_1(L')$ holds. These bilateral bounds force $\Sigma_1(L') \to \Sigma_1(L)$ as $\delta \to 0$, establishing continuity.
\end{proof}

\subsection{Degeneration of Short Cylinders ($L \to 0$)}

As the cylinder length collapses to zero, the geometric band flattens, and the two boundary components collide.

\begin{lemma} \label{lem:limit_L_zero}
The supremum map vanishes as the cylinder collapses: $\lim_{L \to 0} \Sigma_1(L) = 0$.
\end{lemma}

\begin{proof}
This follows directly from the uniform upper bounds established in Section \ref{sec:bounds}. For cylinders of length $L < 1$, evaluating the Rayleigh quotient of the transverse spectral branch ($k=1$) with a constant profile function yielded the upper bound:
\[
    \bar{\sigma}_1(\alpha, \mu) \le \lambda_1(\Sbb^{n-1}) \vol(\Sbb^{n-1})^{\frac{2}{n}} 2^{\frac{2}{n}} L^{\frac{2}{n}}.
\]
The exponent $2/n$ is strictly positive. Taking the supremum over the admissible class gives $\Sigma_1(L) \le C \cdot L^{2/n}$. Letting $L \to 0$, we immediately obtain $\lim_{L \to 0} \Sigma_1(L) = 0$. Geometrically, the vanishing of the longitudinal distance forces the transverse eigenvalues to tend to zero.
\end{proof}

\subsection{Spectral Bound for Infinite Cylinders ($L \to \infty$)}

The asymptotic case $L \to \infty$ corresponds to an infinitely stretching cylinder. Because the functional optimizes over all admissible boundary measures, we distinguish between regimes where the mass remains balanced and regimes where one boundary loses its mass. When a boundary loses its mass, the metric could concentrate near the vanishing boundary. To bound the spectrum without assuming uniform geometric control, we exploit the conformal equivalence of the cylinder to a punctured Euclidean ball. Because the $H^1$-capacity of a point vanishes in dimensions $n \ge 3$, we construct cutoff functions that isolate the energy from any potential metric degeneration.

\begin{lemma} \label{lem:limit_L_infty}
As the cylinder stretches infinitely, the supremum map satisfies:
\[
    \limsup_{L \to \infty} \Sigma_1(L) \le \bar{\sigma}_1(\B^n),
\]
where $\bar{\sigma}_1(\B^n)$ is the normalized Steklov eigenvalue of the standard $n$-dimensional unit ball.
\end{lemma}

\begin{proof}
Let $L_m \to \infty$ be a sequence of half-lengths, and $(\alpha_m, \mu_m) \in \mathcal{M}_{L_m}$ a normalized sequence of configurations approaching $\Sigma_1(L_m)$ (with $\|\alpha_m\|_{L^q([-L_m, L_m])} = 1$ and $\mu_{m, L_m} + \mu_{m, -L_m} = 1$). Passing to a subsequence, assume the boundary measures converge: $\mu_{m, \pm L_m} \to \mu_{\pm}$, with $\mu_+ + \mu_{-} = 1$. We analyze $\sigma_1 = \min(\sigma^{(0)}_1, \sigma^{(1)}_1)$ in two asymptotic cases.

\vspace{0.3cm}

\textbf{Case 1 (Balanced Mass):} Assume $\mu_+ > 0$ and $\mu_{-} > 0$. 
We bound the eigenvalue using the longitudinal branch with the test function $\psi(t) = \frac{t}{L_m} - (\mu_{m, L_m} - \mu_{m, -L_m})$. The shift ensures $\mu_{m, L_m}\psi(L_m) + \mu_{m, -L_m}\psi(-L_m) = 0$. The Steklov denominator evaluates to:
\[
    D_m = \mu_{m, L_m} \psi(L_m)^2 + \mu_{m, -L_m} \psi(-L_m)^2 = 4 \mu_{m, L_m} \mu_{m, -L_m}.
\]
Because the mass is balanced, $\lim_{m \to \infty} D_m = 4\mu_+ \mu_{-} > 0$. With $\dot{\psi}(t) = \frac{1}{L_m}$, applying H\"older's inequality (exponents $q = \frac{n}{n-2}$ and $p = \frac{n}{2}$) bounds the numerator:
\[
    \int_{-L_m}^{L_m} \dot{\psi}^2 \alpha_m \, dt \le \frac{1}{L_m^2} \|\alpha_m\|_{L^q} (2L_m)^{\frac{2}{n}}.
\]
Substituting this into the normalized functional cancels the density norm:
\[
    \bar{\sigma}_1(\alpha_m, \mu_m) \le \vol(\Sbb^{n-1})^{\frac{2}{n}}\frac{2^{\frac{2}{n}}}{4 \mu_{m, L_m} \mu_{m, -L_m}} L_m^{\frac{2-2n}{n}}.
\]
Since $n \ge 3$, the exponent $\frac{2-2n}{n} < 0$. With the denominator bounded away from zero, the upper bound vanishes. Thus, $\lim_{m \to \infty} \bar{\sigma}_1(\alpha_m, \mu_m) = 0$.

\vspace{0.3cm}

\textbf{Case 2 (Vanishing Boundary Mass):} Assume all mass concentrates on one boundary, e.g., $\mu_{-} = 0$ and $\mu_+ = 1$. 
To shield the test functions, we map the cylinder conformally to a Euclidean annulus. The map $x(t, \xi) = e^{t-L_m} \xi$ diffeomorphically sends $[-L_m, L_m] \times \Sbb^{n-1}$ onto $A_m = \B^n \setminus \overline{\B_{r_m}}$, where $r_m = e^{-2L_m} \to 0$. The boundary $\Gamma_{L_m}$ maps to the outer sphere $\partial \B^n$.

Let $\tilde{\alpha}_m(x) = \alpha_m(t)|x|^{2-n}$ be the conformal pushforward to $A_m$, which (as computed in the proof of Lemma \ref{lem:boundary_noncollapse}) satisfies
\[
    \|\tilde{\alpha}_m\|_{L^q(A_m)} = \vol(\Sbb^{n-1})^{\frac{n-2}{n}} \|\alpha_m\|_{L^q([-L_m, L_m])} = \vol(\Sbb^{n-1})^{\frac{n-2}{n}}.
\]
We isolate the origin using a radial logarithmic cutoff $\eta_m(|x|)$:
\[
    \eta_m(r) = 
    \begin{cases} 
      1 & \text{for } r \in [\sqrt{r_m}, 1], \\
      2 \frac{\log r - \log r_m}{-\log r_m} & \text{for } r \in [r_m, \sqrt{r_m}].
    \end{cases}
\]
We test the Rayleigh quotient against the shielded transverse functions $\tilde{u}_i(x) = \eta_m(|x|) x_i$ ($i = 1, \dots, n$). Because $\tilde{u}_i = 0$ on $|x| = r_m$ and $\tilde{u}_i = \xi_i$ on $|x| = 1$, orthogonality to constants on the boundary is satisfied. The boundary denominator evaluates to:
\[
    \sum_{i=1}^n \int_{\partial \B^n} \tilde{u}_i^2 \mu_{m, L_m} \, d\mathcal{H}_{g_0} = \mu_{m, L_m} \vol(\Sbb^{n-1}).
\]
The summed Dirichlet energy density of the shielded functions is:
\[
    E_m(x) := \sum_{i=1}^n |\nabla (\eta_m x_i)|^2 \le n + \frac{4}{(\log r_m)^2} + \frac{4}{-\log r_m}.
\]
Because the logarithmic terms vanish as $r_m \to 0$, $E_m(x)$ is uniformly bounded and converges pointwise on $\B^n \setminus \{0\}$ to $n$. By H\"older's inequality and the Dominated Convergence Theorem:
\[
    \limsup_{m \to \infty} \int_{A_m} \tilde{\alpha}_m E_m \, dx \le \|\tilde{\alpha}_m\|_{L^q(A_m)} \lim_{m \to \infty} \|E_m\|_{L^{n/2}(A_m)} = \vol(\Sbb^{n-1})^{\frac{n-2}{n}} \left(n \vol(\B^n)^{2/n}\right).
\]
Bounding the first eigenvalue by the quotient of these sums yields:
\[
    \limsup_{m \to \infty} \sigma_1(\alpha_m, \mu_m) \le \frac{\vol(\Sbb^{n-1})^{\frac{n-2}{n}} n \vol(\B^n)^{2/n}}{\mu_{m, L_m} \vol(\Sbb^{n-1})} = \frac{n \vol(\B^n)^{2/n}}{\mu_{m, L_m} \vol(\Sbb^{n-1})^{\frac{2}{n}}}.
\]
Substituting this bound into the normalized functional $\bar{\sigma}_1$, the parameters $\mu_{m, L_m} \to 1$, $\|\alpha_m\|_{L^q} = 1$, and $\vol(\Sbb^{n-1})^{\frac{2}{n}}$ cancel out, leaving:
\begin{align*}
    \limsup_{m \to \infty} \bar{\sigma}_1(\alpha_m, \mu_m) &\le \frac{n \vol(\B^n)^{2/n}}{\vol(\Sbb^{n-1})^{\frac{2}{n}}} \vol(\Sbb^{n-1})^{\frac{2}{n}} \\
    &= n \left( \frac{1}{n}\vol(\Sbb^{n-1}) \right)^{2/n} = n^{\frac{n-2}{n}} \vol(\Sbb^{n-1})^{\frac{2}{n}} \equiv \bar{\sigma}_1(\B^n).
\end{align*}
This completes the proof.
\end{proof}

\section{Geometry and Spectral Domination of the Critical $n$-Catenoid}
\label{sec:spectral_domination}

In this section, we explicitly construct the critical $n$-catenoid $\mathcal{C}_n$, prove that its first non-zero Steklov eigenvalue is $\sigma_1(\mathcal{C}_n) = 1$, and show that its normalized Steklov eigenvalue strictly dominates that of the standard $n$-ball $\B^n$ in all dimensions $n \ge 3$. Combining this strict inequality with the asymptotic bounds of Section \ref{sec:moduli_space}, we conclude that the global supremum of $\bar{\sigma}_1$ over the moduli space is achieved at a finite conformal modulus $L^* \in (0, \infty)$.

\subsection{The $n$-Dimensional Catenoid in $\R^{n+1}$}

Let $M(n) \subset \R^{n+1}$ be the $O(n)$-invariant hypersurface of revolution parameterized by $\Phi(t, \xi) = (t, L(t)\xi)$ for $t \in I$ and $\xi \in \Sbb^{n-1}$. The minimal hypersurface equation (vanishing mean curvature) is \cite{SSTZ}
\[
    L L'' = (n-1)\left(1 + (L')^2\right).
\]
After normalizing the neck radius via an ambient homothety, we prescribe $L(0) = 1$ and $L'(0) = 0$, extending evenly by $L(-t) = L(t)$. Differentiating $\left(1 + (L')^2\right) L^{-2(n-1)}$ shows that it is constant, and the initial conditions give the first integral
\[
    1 + (L')^2 = L^{2(n-1)}.
\]
Thus $L''(t) = (n-1)L(t)^{2n-3} > 0$, establishing strict convexity. In particular, $L'(t) > 0$ for $t > 0$, yielding $L'(t) = \sqrt{L(t)^{2(n-1)} - 1}$ for $t \ge 0$. For dimensions $n \ge 3$, the profile blows up at a finite maximal half-height:
\[
    t_{\max} = \int_1^\infty \frac{dx}{\sqrt{x^{2(n-1)} - 1}} < \infty.
\]
As $t \to \pm t_{\max}$, the radius $L(t) \to \infty$ and the tangent becomes vertical.

\subsection{The Free Boundary Condition and First Steklov Eigenvalue}

To realize a piece of $M(n)$ as a free boundary minimal hypersurface in a Euclidean ball, it must intersect the bounding sphere orthogonally at some truncation heights $t = \pm t_*$. At $t = t_*$, the outward normal to the sphere is proportional to the position vector $\vec{X} = (t_*, L(t_*)\xi)$, while the normal to the catenoid is proportional to $(-L'(t_*), \xi)$. Orthogonality ($\vec{X} \cdot \vec{\nu}_{\text{surf}} = 0$) enforces the transversality condition:
\begin{equation} \label{eq:fbms_condition}
    t_* L'(t_*) = L(t_*).
\end{equation}

\begin{lemma}
For $n \ge 3$, there exists a unique $t_* \in (0, t_{\max})$ satisfying \eqref{eq:fbms_condition}.
\end{lemma}

\begin{proof}
Define $h(t) = t L'(t) - L(t)$ on $[0, t_{\max})$. We have $h(0) = -L(0) = -1 < 0$. By strict convexity of $L$, $h'(t) = t L''(t) > 0$ for $t > 0$, so $h$ is strictly increasing on $[0, t_{\max})$. Using $L'(t) = \sqrt{L(t)^{2n-2}-1}$, we have $h(t) = L(t)\left(t\sqrt{L(t)^{2n-4} - L(t)^{-2}} - 1\right) \to +\infty$ as $t \to t_{\max}$. By the Intermediate Value Theorem, there exists a unique root $t_* \in (0, t_{\max})$.
\end{proof}

\begin{remark}
While the two-dimensional catenoid extends to infinite axial height, the profile of the $n$-catenoid ($n \ge 3$) is confined to the finite interval $(-t_{\max}, t_{\max})$.
\end{remark}

Let $t_* \in (0, t_{\max})$ be the unique root of \eqref{eq:fbms_condition}. The truncated catenoid over $[-t_*, t_*]$ meets the sphere of radius $R_n := \sqrt{t_*^2 + L(t_*)^2}$ orthogonally. Rescaling by $1/R_n$ embeds it into the unit ball $\B^{n+1}$, defining the \emph{critical $n$-catenoid}:
\[
    \mathcal{C}_n := \left\{ \frac{1}{R_n} (t, L(t)\xi) \;\middle|\; t \in [-t_*, t_*], \; \xi \in \Sbb^{n-1} \right\}.
\]
Let $g_{\mathrm{cat}}$ be the induced Euclidean metric on $\mathcal{C}_n$. Because $\mathcal{C}_n$ is a free boundary minimal hypersurface in $\B^{n+1}$, its ambient coordinate functions $x_1, \dots, x_{n+1}$ are harmonic ($\Delta_{g_{\mathrm{cat}}} x_i = 0$) and satisfy $\partial_\nu x_i = x_i$ on $\partial \mathcal{C}_n$, so $\sigma = 1$ is a Steklov eigenvalue of $(\mathcal{C}_n, g_{\mathrm{cat}}, 1)$ of multiplicity at least $n+1$.

Let us verify that $\sigma = 1$ is indeed the \emph{first} non-zero Steklov eigenvalue $\sigma_1(\mathcal{C}_n, g_{\mathrm{cat}}, 1) = 1$. By Theorem \ref{thm:spectrum_computation}, the first non-zero Steklov eigenvalue of an $O(n)$-invariant cylinder is $\sigma_1 = \min\{\sigma_1^{(0)}, \sigma_1^{(1)}\}$, where $\sigma_1^{(0)}$ is the unique non-zero longitudinal eigenvalue ($k=0$) and $\sigma_1^{(1)}$ is the lowest transverse eigenvalue ($k=1$, characterized by a radial eigenprofile without interior zeros). On $\mathcal{C}_n$, the longitudinal coordinate function $x_1 = t/R_n$ is independent of $\xi \in \Sbb^{n-1}$ and non-constant, so it spans the one-dimensional non-zero longitudinal eigenspace, yielding $\sigma_1^{(0)} = 1$. Simultaneously, the transverse coordinate functions $x_{i+1} = \frac{L(t)}{R_n}\xi_i$ ($i = 1, \dots, n$) belong to the $k=1$ branch, and their radial profile $Y_1(t) = \frac{L(t)}{R_n} > 0$ never vanishes on $[-t_*, t_*]$. By Sturm--Liouville oscillation theory, the unique eigenprofile of \eqref{eq:profile_ode}--\eqref{eq:profile_bc} with no zeros on $[-t_*, t_*]$ corresponds to the lowest eigenvalue of the $k=1$ branch, which proves $\sigma_1^{(1)} = 1$. Therefore,
\[
    \sigma_1(\mathcal{C}_n, g_{\mathrm{cat}}, 1) = \min\left\{\sigma_1^{(0)}, \sigma_1^{(1)}\right\} = 1.
\]

\subsection{Spectral Domination over the Standard Ball}

\begin{proposition} \label{prop:catenoid_beats_ball}
For all dimensions $n \ge 3$, the normalized first Steklov eigenvalue of the critical $n$-catenoid strictly exceeds that of the standard $n$-dimensional unit ball:
\[
    \bar{\sigma}_1(\mathcal{C}_n) > \bar{\sigma}_1(\B^n).
\]
\end{proposition}

\begin{proof}
For the unit ball $\B^n$, the first Steklov eigenvalue is $\sigma_1(\B^n) = 1$, the boundary area is $\vol(\Sbb^{n-1})$, and the volume is $\frac{1}{n} \vol(\Sbb^{n-1})$, so
\begin{equation} \label{eq:value_ball}
    \bar{\sigma}_1(\B^n) = 1 \cdot \vol(\Sbb^{n-1}) \cdot \left( \frac{\vol(\Sbb^{n-1})}{n} \right)^{\frac{2-n}{n}} = n^{\frac{n-2}{n}} \vol(\Sbb^{n-1})^{\frac{2}{n}}.
\end{equation}

For the critical $n$-catenoid $\mathcal{C}_n$, we have shown above that $\sigma_1(\mathcal{C}_n) = 1$. Let $R_* := L(t_*)/R_n$ denote the radius of each of its two boundary spheres. The boundary area is $\mathcal{H}_{g_{\mathrm{cat}}}(\partial \mathcal{C}_n) = 2 \vol(\Sbb^{n-1}) R_*^{n-1}$. Since $\mathcal{C}_n \subset \B^{n+1}$ is minimal and meets $\partial \B^{n+1}$ orthogonally, applying the divergence theorem to the tangential projection of the position vector $\overrightarrow{X}$ (which satisfies $\operatorname{div}_{\mathcal{C}_n}(\overrightarrow{X}^\top) = n$ on $\mathcal{C}_n$ and $\langle \overrightarrow{X}, \nu \rangle = 1$ on $\partial \mathcal{C}_n$) yields $\vol(\mathcal{C}_n) = \frac{1}{n}\mathcal{H}_{g_{\mathrm{cat}}}(\partial \mathcal{C}_n) = \frac{2}{n} \vol(\Sbb^{n-1}) R_*^{n-1}$. Substituting these into \eqref{eq:KM_normalization} gives:
\begin{equation} \label{eq:value_catenoid}
    \bar{\sigma}_1(\mathcal{C}_n) = n^{\frac{n-2}{n}} \left( 2 \vol(\Sbb^{n-1}) R_*^{n-1} \right)^{\frac{2}{n}}.
\end{equation}
Dividing \eqref{eq:value_catenoid} by \eqref{eq:value_ball} gives:
\[
    \frac{\bar{\sigma}_1(\mathcal{C}_n)}{\bar{\sigma}_1(\B^n)} = \left( 2 R_*^{n-1} \right)^{\frac{2}{n}}.
\]

Set $L_* := L(t_*)$. By \eqref{eq:fbms_condition}, we have $t_* = L_* / L'(t_*) = L_* / \sqrt{L_*^{2n-2}-1}$, and hence
\[
    R_* = \frac{L_*}{\sqrt{t_*^2 + L_*^2}} = \frac{\sqrt{L_*^{2n-2}-1}}{L_*^{n-1}}.
\]
Setting $X := L_*^{2n-2} > 1$, we have $R_* = \left( \frac{X-1}{X} \right)^{\frac{1}{2}}$, and the ratio simplifies to:
\begin{equation} \label{eq:ratio_X}
    \frac{\bar{\sigma}_1(\mathcal{C}_n)}{\bar{\sigma}_1(\B^n)} = \left( 2 \left( \frac{X-1}{X} \right)^{\frac{n-1}{2}} \right)^{\frac{2}{n}}.
\end{equation}
Thus $\bar{\sigma}_1(\mathcal{C}_n) > \bar{\sigma}_1(\B^n)$ if and only if $2 \left( \frac{X-1}{X} \right)^{\frac{n-1}{2}} > 1$, which is equivalent to:
\begin{equation} \label{eq:domination_threshold}
    X > T_n := \frac{2^{\frac{2}{n-1}}}{2^{\frac{2}{n-1}} - 1}.
\end{equation}
By Lemma \ref{lem:analytic_root_bound} below, $X > T_n$ holds for all $n \ge 3$.
\end{proof}

\begin{lemma}[Analytic Bound on the Free Boundary Root] \label{lem:analytic_root_bound}
For all $n \ge 3$, the unique root $L_* = L(t_*)$ of the $n$-catenoid free boundary equation satisfies $X = L_*^{2n-2} > T_n$.
\end{lemma}

\begin{proof}
Since $t_* = \int_1^{L_*} \frac{ds}{\sqrt{s^{2n-2}-1}}$ and $t_* = \frac{L_*}{\sqrt{L_*^{2n-2}-1}}$, the root $L_* > 1$ is uniquely determined by:
\[
    \frac{L_*}{\sqrt{L_*^{2n-2}-1}} = \int_1^{L_*} \frac{1}{\sqrt{s^{2n-2}-1}} \, ds.
\]
Substituting $x = s^{n-1}$ ($ds = \frac{1}{n-1}x^{-\frac{n-2}{n-1}}dx$) and using $x^{\frac{n-2}{n-1}} > 1$ for $x > 1$ yields the strict upper bound:
\[
    \int_1^{L_*} \frac{1}{\sqrt{s^{2n-2}-1}} \, ds = \frac{1}{n-1}\int_1^{L_*^{n-1}} \frac{dx}{x^{\frac{n-2}{n-1}}\sqrt{x^2-1}} < \frac{1}{n-1} \int_1^{L_*^{n-1}} \frac{dx}{\sqrt{x^2-1}} = \frac{\operatorname{arccosh}(L_*^{n-1})}{n-1}.
\]
Setting $Y := L_*^{n-1} = \sqrt{X} > 1$ and using $L_* = Y^{\frac{1}{n-1}} > 1$, we obtain:
\[
    \frac{1}{\sqrt{Y^2-1}} < \frac{Y^{\frac{1}{n-1}}}{\sqrt{Y^2-1}} < \frac{\operatorname{arccosh}(Y)}{n-1},
\]
and therefore $Y$ satisfies the strict inequality:
\begin{equation} \label{eq:arccosh_bound}
    \operatorname{arccosh}(Y) > \frac{n-1}{\sqrt{Y^2-1}}.
\end{equation}

Let us evaluate both sides of \eqref{eq:arccosh_bound} at the threshold value $Y_{\text{thresh}} := \sqrt{T_n}$. Set $c := 2^{\frac{2}{n-1}} \in (1, 2]$ (since $n \ge 3$) and $u := \sqrt{c-1} \in (0, 1]$. Then $T_n = \frac{c}{c-1} = \frac{1+u^2}{u^2}$, so $Y_{\text{thresh}} = \frac{\sqrt{1+u^2}}{u}$ and $\sqrt{Y_{\text{thresh}}^2 - 1} = \frac{1}{u}$. Moreover, $c = 1 + u^2 = 2^{\frac{2}{n-1}}$ gives $n - 1 = \frac{2\ln 2}{\ln(1+u^2)}$. Using the elementary inequality $\ln(1+u^2) < u^2$ for $u > 0$, the right-hand side of \eqref{eq:arccosh_bound} at $Y_{\text{thresh}}$ satisfies:
\[
    \frac{n-1}{\sqrt{Y_{\text{thresh}}^2 - 1}} = \frac{(2\ln 2) u}{\ln(1+u^2)} > \frac{2\ln 2}{u} > \frac{1.38}{u}.
\]
On the other hand, since $u \in (0, 1]$, $\ln(1/u) \le \frac{1}{u} - 1$ and $\operatorname{arccosh}(Y_{\text{thresh}}) = \ln\left(\frac{1 + \sqrt{1+u^2}}{u}\right)$, the left-hand side of \eqref{eq:arccosh_bound} at $Y_{\text{thresh}}$ satisfies:
\[
    \operatorname{arccosh}(Y_{\text{thresh}}) = \ln\left(1 + \sqrt{1+u^2}\right) + \ln\left(\frac{1}{u}\right)
    \le \ln\left(1+\sqrt{2}\right) - 1 + \frac{1}{u} < \frac{1}{u},
\]
because $\ln(1+\sqrt{2}) \approx 0.8814 < 1$. Combining both estimates gives
\[
    \operatorname{arccosh}(Y_{\text{thresh}}) < \frac{1}{u} < \frac{1.38}{u} < \frac{n-1}{\sqrt{Y_{\text{thresh}}^2 - 1}},
\]
so \eqref{eq:arccosh_bound} fails at $Y_{\text{thresh}}$. Since $y \mapsto \operatorname{arccosh}(y)$ is strictly increasing and $y \mapsto \frac{n-1}{\sqrt{y^2-1}}$ is strictly decreasing on $(1, \infty)$, the root $Y$ satisfying \eqref{eq:arccosh_bound} must satisfy $Y > Y_{\text{thresh}}$, i.e., $X = Y^2 > T_n$.
\end{proof}

\subsection{Existence of a Global Maximizer at a Finite Conformal Modulus}

\begin{proposition} \label{prop:existence_L_star}
Set $\Sigma_1^{\mathrm{global}} := \sup_{L \in (0, \infty)} \Sigma_1(L)$. Then $\Sigma_1^{\mathrm{global}} > \bar{\sigma}_1(\B^n)$, and there exists a finite conformal modulus $L^* \in (0, \infty)$ such that $\Sigma_1(L^*) = \Sigma_1^{\mathrm{global}}$. Consequently, by Theorem \ref{thm:exact_criticality}, there exists a smooth $O(n)$-invariant metric $g_{\max}$ and a smooth $O(n)$-invariant boundary weight $\rho_{\max}$ on $\Omega_{L^*} = [-L^*, L^*] \times \Sbb^{n-1}$ achieving the global supremum $\Sigma_1^{\mathrm{global}}$ over all $O(n)$-invariant metrics and weights on the cylinder.
\end{proposition}

\begin{proof}
Since the critical $n$-catenoid $\mathcal{C}_n$ is an $O(n)$-invariant cylinder of some finite conformal modulus $L_{\mathrm{cat}} \in (0, \infty)$, Proposition \ref{prop:catenoid_beats_ball} implies:
\[
    \Sigma_1^{\mathrm{global}} \ge \Sigma_1(L_{\mathrm{cat}}) \ge \bar{\sigma}_1(\mathcal{C}_n) > \bar{\sigma}_1(\B^n) > 0.
\]
By Lemma \ref{lem:limit_L_zero} and Lemma \ref{lem:limit_L_infty}, we have $\lim_{L \to 0} \Sigma_1(L) = 0$ and $\limsup_{L \to \infty} \Sigma_1(L) \le \bar{\sigma}_1(\B^n) < \Sigma_1^{\mathrm{global}}$. Since $L \mapsto \Sigma_1(L)$ is continuous on $(0, \infty)$ by Lemma \ref{lem:continuity_lambda}, it attains its global maximum at some interior point $L^* \in (0, \infty)$. Since $\Sigma_1(L^*) = \Sigma_1^{\mathrm{global}} > \bar{\sigma}_1(\B^n)$, Theorem \ref{thm:exact_criticality} yields a smooth $O(n)$-invariant maximizer $(g_{\max}, \rho_{\max})$ on $\Omega_{L^*}$ achieving $\Sigma_1(L^*) = \Sigma_1^{\mathrm{global}}$.
\end{proof}

\section{Criticality and Characterization of the Maximizer}
\label{sec:criticality}

To complete the proof of the Main Theorem (Theorem \ref{thm:main_theorem}), it remains to show that the smooth global maximizer $(g_{\max}, \rho_{\max})$ obtained in Proposition \ref{prop:existence_L_star} is isometric (up to scaling) to the critical $n$-catenoid $\mathcal{C}_n$.

\subsection{Criticality Conditions in the Space of $O(n)$-Invariant Metrics}

Let $\Omega = [-1, 1] \times \Sbb^{n-1}$ be equipped with the $O(n)$ action of Definition \ref{def:On_invariant}.

\begin{theorem}\label{thm:Crit_On}
Let $(g_0,\rho_0)$ be a smooth pair of $O(n)$-invariant metric and boundary weight on $\Omega$ that are critical for $\bar{\sigma}_1$ in the space of all smooth $O(n)$-invariant metrics and weights. Then there exist $m>0$, an orthogonal representation $A\colon O(n)\to O(m+1)$, and a collection of $\sigma_1(g_0, \rho_0)$-eigenfunctions $u_1, \dots, u_{m+1}$ such that $\Phi = (u_1, \dots, u_{m+1})$ defines an $O(n)$-equivariant free boundary minimal immersion 
\[
    \Phi\colon (\Omega, \partial \Omega) \to \B^{m+1},\quad \Phi(\kappa\cdot x) = A_\kappa(\Phi(x)),
\]
such that $g_0 = f_0 \Phi^* g_{\mathrm{Eucl}}$ for a positive constant $f_0 > 0$ on $\Omega$, $\rho_0$ is constant on $\partial \Omega$, and $\partial_\nu\Phi=\sigma_1(g_0,\rho_0)\rho_0\Phi$ on $\partial\Omega$.
\end{theorem}

\begin{proof}
By scaling $g_0$ and $\rho_0$, we may assume without loss of generality that $\vol(\Omega, g_0) = 1$, $\int_{\partial \Omega} \rho_0 \, d\mathcal{H}_{g_0} = 1$, and consequently $\sigma_1(g_0, \rho_0) = \bar{\sigma}_1(\Omega, g_0, \rho_0) =: \sigma_1$. Let $E_1(g_0, \rho_0)$ be the $\sigma_1$-eigenspace.

We consider the class of smooth $O(n)$-invariant variations $(g_t, \rho_t)$ with $(g_0, \rho_0)$ at $t=0$ generated by a \emph{pointwise trace-free} $O(n)$-invariant symmetric $(0,2)$-tensor $h \in C^\infty_{O(n)}(S^2(\Omega))$ and an $O(n)$-invariant boundary function $\eta \in C^\infty_{O(n)}(\partial \Omega)$ satisfying
\begin{equation} \label{eq:zero_mean_constraints}
    \operatorname{tr}_{g_0} h = \langle g_0, h \rangle \equiv 0 \quad \text{on } \Omega \qquad \text{and} \qquad \int_{\partial \Omega} \eta \rho_0 \, d\mathcal{H}_{g_0} = 0.
\end{equation}
Specifically, setting $A := g_0^{-1}h$, we may take $g_t(v, w) := g_0(e^{tA}v, w)$ (which preserves the volume form $dv_{g_t} = dv_{g_0}$ pointwise because $\operatorname{tr} A = 0$, and satisfies $\dot{g}_0 = h$) and define $\rho_t$ by $\rho_t \, d\mathcal{H}_{g_t} = (1 + t\eta)\rho_0 \, d\mathcal{H}_{g_0}$. For sufficiently small $t$, both $g_t$ and $\rho_t$ remain smooth, positive, and $O(n)$-invariant, with $\vol(\Omega, g_t) = 1$ and $\int_{\partial \Omega} \rho_t \, d\mathcal{H}_{g_t} = 1$.

By standard perturbation theory (analogous to Lemma 4.1 in \cite{KN}; see also \cite{FS12, KM}), the directional derivative of $\sigma_1(g_t, \rho_t)$ at $t=0$ is governed by the quadratic form
\[
    Q_{h, \eta}(u) := -\int_\Omega \langle T(u), h \rangle \, dv_{g_0} - \sigma_1 \int_{\partial \Omega} u^2 \eta \rho_0 \, d\mathcal{H}_{g_0}
\]
for $u \in E_1(g_0, \rho_0)$ with $\int_{\partial \Omega} u^2 \rho_0 \, d\mathcal{H}_{g_0} = 1$, where $T(u) := du \otimes du - \frac{1}{2}|\nabla u|_{g_0}^2 g_0$ is the usual stress-energy tensor. Let
\[
    S(u) := du \otimes du - \frac{1}{n}|\nabla u|_{g_0}^2 g_0
\]
denote the pointwise trace-free stress-energy tensor. Because $\operatorname{tr}_{g_0} h \equiv 0$ pointwise on $\Omega$, we have $\langle T(u), h \rangle = \langle S(u), h \rangle$ pointwise on $\Omega$, and therefore
\[
    Q_{h, \eta}(u) = -\int_\Omega \langle S(u), h \rangle \, dv_{g_0} - \sigma_1 \int_{\partial \Omega} u^2 \eta \rho_0 \, d\mathcal{H}_{g_0}.
\]

Let $\mathcal{H} := L^2(S^2(\Omega), dv_{g_0}) \times L^2(\partial \Omega, \rho_0 d\mathcal{H}_{g_0})$ and define the convex hull
\[
    \mathcal{K} := \operatorname{conv} \left\{ \left(S(u), \, \sigma_1 u^2\big|_{\partial \Omega} \right) \;\middle|\; u \in E_1(g_0, \rho_0), \; \int_{\partial \Omega} u^2 \rho_0 \, d\mathcal{H}_{g_0} = 1 \right\} \subset \mathcal{H}.
\]
Since $E_1(g_0, \rho_0)$ is finite-dimensional, its unit sphere in $L^2(\partial \Omega, \rho_0 d\mathcal{H}_{g_0})$ is compact, and its image under the continuous quadratic map $u \mapsto (S(u), \sigma_1 u^2|_{\partial \Omega})$ is a compact subset of a finite-dimensional subspace of $\mathcal{H}$. By Carathéodory's theorem, the convex hull $\mathcal{K}$ is therefore a compact convex subset of $\mathcal{H}$.

\textbf{Claim:} $\left(0, \sigma_1\right) \in \mathcal{K}$.

Suppose for contradiction that $(0, \sigma_1) \notin \mathcal{K}$. By the Hahn--Banach separation theorem applied to the compact convex set $\mathcal{K}$ and the point $(0, \sigma_1)$, there exists $(\omega, \xi) \in \mathcal{H}$ and a constant $c \in \R$ such that
\[
    \langle (\omega, \xi), (S(u), \sigma_1 u^2) \rangle_{\mathcal{H}} \le c < \langle (\omega, \xi), (0, \sigma_1) \rangle_{\mathcal{H}} = \sigma_1 \int_{\partial \Omega} \xi \rho_0 \, d\mathcal{H}_{g_0}
\]
for all $u \in E_1(g_0, \rho_0)$ with $\int_{\partial \Omega} u^2 \rho_0 \, d\mathcal{H}_{g_0} = 1$. Replacing $\xi$ by $\xi - c/\sigma_1$ (which subtracts $c$ from both sides since $\int_{\partial \Omega} u^2 \rho_0 \, d\mathcal{H}_{g_0} = \int_{\partial \Omega} \rho_0 \, d\mathcal{H}_{g_0} = 1$), we may take $c = 0$, so that $\sigma_1 \int_{\partial \Omega} \xi \rho_0 \, d\mathcal{H}_{g_0} > 0$ and
\[
    \int_\Omega \langle S(u), \omega \rangle \, dv_{g_0} + \sigma_1 \int_{\partial \Omega} u^2 \xi \rho_0 \, d\mathcal{H}_{g_0} \le 0
\]
for all $u \in E_1(g_0, \rho_0)$ with $\int_{\partial \Omega} u^2 \rho_0 \, d\mathcal{H}_{g_0} = 1$. Since $g_0$ and $\rho_0$ are $O(n)$-invariant, for every $\kappa \in O(n)$ and normalized $u \in E_1(g_0, \rho_0)$, the pullback $u \circ R_{\kappa^{-1}}$ is also a normalized eigenfunction in $E_1(g_0, \rho_0)$ and satisfies $S(u \circ R_{\kappa^{-1}}) = R_{\kappa^{-1}}^* S(u)$. Replacing $u$ by $u \circ R_{\kappa^{-1}}$, changing variables $y = R_{\kappa^{-1}}(x)$, and averaging over $\kappa \in O(n)$ with respect to the normalized Haar measure $d\mu(\kappa)$, we obtain the $O(n)$-invariant fields
\[
    W := \int_{O(n)} R_\kappa^* \omega \, d\mu(\kappa) \in L^2_{O(n)}(S^2(\Omega)), \qquad B := \int_{O(n)} (\xi \circ R_\kappa) \, d\mu(\kappa) \in L^2_{O(n)}(\partial \Omega),
\]
which satisfy $\sigma_1 \int_{\partial \Omega} B \rho_0 \, d\mathcal{H}_{g_0} > 0$ and $\int_\Omega \langle S(u), W \rangle \, dv_{g_0} + \sigma_1 \int_{\partial \Omega} u^2 B \rho_0 \, d\mathcal{H}_{g_0} \le 0$ for all normalized $u \in E_1(g_0, \rho_0)$. Approximating $W$ and $B$ in $L^2$ by smooth $O(n)$-invariant fields (which preserves the strict inequality $\int_\Omega \langle S(u), W \rangle \, dv_{g_0} + \sigma_1 \int_{\partial \Omega} u^2 B \rho_0 \, d\mathcal{H}_{g_0} - \sigma_1 \int_{\partial \Omega} B \rho_0 \, d\mathcal{H}_{g_0} < 0$ uniformly on the compact unit sphere of $E_1(g_0, \rho_0)$), we define the smooth $O(n)$-invariant variation $(h, \eta)$ by the \emph{pointwise} trace-free projection of $W$ and the zero-mean projection of $B$:
\[
    h := W - \frac{1}{n}\bigl(\operatorname{tr}_{g_0} W\bigr)g_0, \qquad \eta := B - \int_{\partial \Omega} B \rho_0 \, d\mathcal{H}_{g_0}.
\]
Then $(h, \eta)$ satisfies \eqref{eq:zero_mean_constraints}, and since $S(u)$ is pointwise trace-free ($\langle S(u), g_0 \rangle \equiv 0$), we have $\langle S(u), h \rangle = \langle S(u), W \rangle$ pointwise on $\Omega$. Consequently, for every $u \in E_1(g_0, \rho_0)$ with $\int_{\partial \Omega} u^2 \rho_0 \, d\mathcal{H}_{g_0} = 1$,
\begin{align*}
    -Q_{h, \eta}(u) &= \int_\Omega \langle S(u), h \rangle \, dv_{g_0} + \sigma_1 \int_{\partial \Omega} u^2 \eta \rho_0 \, d\mathcal{H}_{g_0} \\
    &= \int_\Omega \langle S(u), W \rangle \, dv_{g_0} + \sigma_1 \int_{\partial \Omega} u^2 B \rho_0 \, d\mathcal{H}_{g_0} - \sigma_1 \int_{\partial \Omega} B \rho_0 \, d\mathcal{H}_{g_0} < 0.
\end{align*}
Thus $Q_{h, \eta}(u) > 0$ for all normalized $u \in E_1(g_0, \rho_0)$ (and $Q_{-h, -\eta}(u) < 0$), which contradicts the criticality of $(g_0, \rho_0)$ for $\bar{\sigma}_1$ (cf. Lemma 4.2 in \cite{KN}). This proves the claim $(0, \sigma_1) \in \mathcal{K}$.

Consequently, there exist convex weights $\alpha_1, \dots, \alpha_{m+1} > 0$ with $\sum_{i=1}^{m+1} \alpha_i = 1$ and normalized eigenfunctions $\tilde{u}_1, \dots, \tilde{u}_{m+1} \in E_1(g_0, \rho_0)$ such that $\sum_{i=1}^{m+1} \alpha_i S(\tilde{u}_i) = 0$ in $\Omega$ and $\sum_{i=1}^{m+1} \alpha_i \sigma_1 \tilde{u}_i^2 = \sigma_1$ on $\partial \Omega$. Setting $u_i := \sqrt{\alpha_i}\tilde{u}_i \in E_1(g_0, \rho_0)$, the map $\Phi = (u_1, \dots, u_{m+1}) \colon \Omega \to \R^{m+1}$ satisfies
\begin{equation} \label{eq:sum_Su_zero}
    \sum_{i=1}^{m+1} S(u_i) = \sum_{i=1}^{m+1} du_i \otimes du_i - \frac{1}{n}|d\Phi|_{g_0}^2 g_0 = 0 \quad \text{in } \Omega, \qquad \sum_{i=1}^{m+1} u_i^2 = 1 \quad \text{on } \partial \Omega.
\end{equation}
By averaging $\Phi$ over the action of $O(n)$ on $L^2(O(n), \R^{m+1})$ exactly as in the proof of Theorem 4.3 in \cite{KN}, we may assume that $\Phi$ is $O(n)$-equivariant with respect to an orthogonal representation $A \colon O(n) \to O(m+1)$.

Since each $u_i$ is $g_0$-harmonic ($\Delta_{g_0} u_i = 0$), we have $-\Delta_{g_0} |\Phi|^2 = -2|d\Phi|_{g_0}^2 \le 0$, so $|\Phi|^2$ is subharmonic on $\Omega$ and equals $1$ on $\partial \Omega$. By the maximum principle, $|\Phi|^2 \le 1$ on $\Omega$, so $\Phi(\Omega) \subset \B^{m+1}$, and by the Hopf lemma, $\partial_\nu |\Phi|^2 = 2\sigma_1 \rho_0 |\Phi|^2 = 2\sigma_1 \rho_0 > 0$ on $\partial \Omega$.

Finally, let us show that $|d\Phi|_{g_0}^2$ is constant on $\Omega$ and that $\Phi$ is a free boundary minimal immersion. By \eqref{eq:sum_Su_zero}, the pullback of the Euclidean metric is
\[
    \Phi^* g_{\mathrm{Eucl}} = \sum_{i=1}^{m+1} du_i \otimes du_i = \frac{1}{n}|d\Phi|_{g_0}^2 g_0.
\]
Recall that on any Riemannian manifold $(\Omega, g_0)$, the divergence of the stress-energy tensor of a harmonic function $u_i$ satisfies $\operatorname{div}_{g_0}\left(du_i \otimes du_i - \frac{1}{2}|\nabla u_i|_{g_0}^2 g_0\right) = (\Delta_{g_0} u_i) du_i = 0$. Summing over $i = 1, \dots, m+1$ and substituting $\sum du_i \otimes du_i = \frac{1}{n}|d\Phi|_{g_0}^2 g_0$, we obtain:
\[
    0 = \operatorname{div}_{g_0} \left( \left(\frac{1}{n} - \frac{1}{2}\right)|d\Phi|_{g_0}^2 g_0 \right) = \frac{2-n}{2n} \, d\left(|d\Phi|_{g_0}^2\right).
\]
Since $n \ge 3$, we have $\frac{2-n}{2n} \neq 0$, which forces $d\left(|d\Phi|_{g_0}^2\right) \equiv 0$ on $\Omega$. Therefore, $|d\Phi|_{g_0}^2 \equiv n c_0$ is a positive constant on $\Omega$, and $\Phi^* g_{\mathrm{Eucl}} = c_0 g_0$ (so $g_0 = f_0 \Phi^* g_{\mathrm{Eucl}}$ with the global constant $f_0 = c_0^{-1} > 0$). Because $\Phi^* g_{\mathrm{Eucl}}$ is homothetic to $g_0$, the $g_0$-harmonicity $\Delta_{g_0} \Phi = 0$ immediately implies $\Delta_{\Phi^* g_{\mathrm{Eucl}}} \Phi = 0$, which means that $\Phi \colon \Omega \to \B^{m+1}$ is a minimal immersion. On $\partial \Omega$, the condition $\partial_\nu \Phi = \sigma_1 \rho_0 \Phi$ shows that $\Phi(\Omega)$ meets $\partial \B^{m+1}$ orthogonally, and taking lengths with respect to $\Phi^* g_{\mathrm{Eucl}} = c_0 g_0$ yields $\sqrt{c_0} = |\partial_\nu \Phi|_{\mathrm{Eucl}} = \sigma_1 \rho_0 |\Phi| = \sigma_1 \rho_0$, so $\rho_0$ is constant on $\partial \Omega$.
\end{proof}

\subsection{Dimension of the Target Ball and Identification of the Global Maximizer}

\begin{theorem}\label{thm:Crit_Catenoid} 
Let $\Omega = [-1, 1] \times \Sbb^{n-1}$ with $n \ge 3$. Suppose that $(g_0, \rho_0)$ is a smooth $O(n)$-invariant configuration that is critical for $\bar{\sigma}_1$ in the space of all $O(n)$-invariant metrics and weights. Then, the corresponding equivariant free boundary minimal immersion 
\[
    \Phi \colon (\Omega, \partial \Omega) \to \B^{m+1}
\]
given by Theorem \ref{thm:Crit_On}, is valid for $m=n$, with $O(n)$ acting on $\R^{n+1} = \R \times \R^n$ by the trivial action on the first coordinate and the standard orthogonal action on $\R^n$.
\end{theorem}

\begin{proof}
Identifying $(\Omega, g_0)$ with its conformal cylinder model $(\Omega_L, g_L)$ via Proposition \ref{prop:moduli_equivalence}, since $\Phi$ is $O(n)$-equivariant and its coordinate functions lie in the $\sigma_1(g_0, \rho_0)$-eigenspace, the span of its coordinate functions decomposes into irreducible representations of $O(n)$ appearing in the $k=0$ (longitudinal, dimension $1$) and $k=1$ (transverse, dimension $n$) Steklov branches (by Theorem \ref{thm:spectrum_computation} and Corollary \ref{cor:multiplicity}). 

If the $n$-dimensional standard representation ($k=1$) were absent, $\Phi$ would be constant along each cross-section $\{t\} \times \Sbb^{n-1}$, contradicting that $\Phi$ is an immersion ($\Phi^* g_{\mathrm{Eucl}} = c_0 g_0$ with $c_0 > 0$). Thus the $k=1$ representation must appear with multiplicity $1$, contributing $n$ coordinates $a(t)\xi_1, \dots, a(t)\xi_n$ with $a(t) > 0$. Similarly, if the trivial representation ($k=0$) were absent, then $\Phi(t, \xi) = (a(t)\xi_1, \dots, a(t)\xi_n) \in \B^n$ would satisfy $a(\pm L) = 1$, so by Rolle's theorem there would exist $t_0 \in (-L, L)$ with $\dot{a}(t_0) = 0$, implying $d\Phi(\partial_t) = 0$ at $t_0$ and contradicting that $\Phi$ is an immersion. Therefore, both the $k=0$ and $k=1$ branches must appear in $\Phi$. Hence, $m = n$ and up to an orthogonal transformation of $\R^{n+1}$,
\[
    \Phi(t, \xi) = \left( b(t), a(t)\xi_1, \dots, a(t)\xi_n \right), \quad (t, \xi) \in [-L, L] \times \Sbb^{n-1}.
\]
\end{proof}

We now combine Proposition \ref{prop:existence_L_star}, Theorem \ref{thm:Crit_On}, and Theorem \ref{thm:Crit_Catenoid} to complete the proof of our Main Theorem.

\begin{theorem}[Completion of the Proof of Theorem \ref{thm:main_theorem}] \label{thm:global_maximizer_extraction}
The smooth global maximizer $(g_{\max}, \rho_{\max})$ achieved at the optimal conformal modulus $L^* \in (0, \infty)$ is uniquely characterized, up to scaling and diffeomorphism, as the metric induced by the critical $n$-catenoid $\mathcal{C}_n \subset \B^{n+1}$ with constant boundary weight.
\end{theorem}

\begin{proof}
By Proposition \ref{prop:moduli_equivalence} and Proposition \ref{prop:existence_L_star}, pulling back $(g_{\max}, \rho_{\max})$ from $\Omega_{L^*}$ to $\Omega = [-1, 1] \times \Sbb^{n-1}$ yields a smooth $O(n)$-invariant metric and weight that globally maximize $\bar{\sigma}_1$ over the space of all smooth $O(n)$-invariant metrics and weights on $\Omega$. In particular, $(g_{\max}, \rho_{\max})$ is critical for $\bar{\sigma}_1$ under all smooth $O(n)$-invariant variations of metric and weight.

Applying Theorem \ref{thm:Crit_On} and Theorem \ref{thm:Crit_Catenoid} to $(g_{\max}, \rho_{\max}) = (f(t)^2(dt^2 + g_{\Sbb^{n-1}}), \rho_{\max})$ on $\Omega_{L^*}$, we obtain an $O(n)$-equivariant free boundary minimal immersion $\Phi \colon (\Omega_{L^*}, \partial \Omega_{L^*}) \to \B^{n+1}$ of the form $\Phi(t, \xi) = (b(t), a(t)\xi)$ such that $g_{\max} = f_0 \Phi^* g_{\mathrm{Eucl}}$ for a constant $f_0 > 0$ and $\rho_{\max}$ is constant on $\partial \Omega_{L^*}$. The longitudinal equation $(f(t)^{n-2} \dot{b}(t))' = 0$ combined with the non-constancy of $b$ implies that $b$ is strictly monotone on $[-L^*, L^*]$, while the first transverse eigenprofile satisfies $a(t) > 0$ everywhere by the one-dimensional ground-state property (as in Corollary \ref{cor:embedding_property}). Therefore, the immersion $\Phi$ is a smooth embedding.

Thus, the image $\Phi(\Omega_{L^*}) \subset \B^{n+1}$ is a connected, $O(n)$-invariant, embedded free boundary minimal hypersurface of revolution in $\B^{n+1}$. By the classification of rotational free boundary minimal hypersurfaces in the unit ball (see, e.g., \cite[Theorem 2]{FGM17}), any such hypersurface with the topology of a cylinder is a piece of an $n$-catenoid. By the uniqueness of the free boundary truncation height $t_* \in (0, t_{\max})$ established in Section \ref{sec:spectral_domination}, $\Phi(\Omega_{L^*})$ coincides with the critical $n$-catenoid $\mathcal{C}_n$ up to an ambient orthogonal transformation of $\R^{n+1}$. This completes the proof of Theorem \ref{thm:main_theorem}.
\end{proof}

\appendix

\section{Numerical Tables and High-Dimensional Asymptotics}
\label{sec:numerics_3D}

\begin{figure}[htpb]
    \centering
    \begin{tikzpicture}
    \begin{axis}[
        xlabel={Dimension $n$},
        ylabel={Ratio $\bar{\sigma}_1(\mathcal{C}_n)/\bar{\sigma}_1(\mathbb{B}^n)$},
        grid=major,
        ymin=1, ymax=1.5,
        xmin=2.5, xmax=13,
        width=0.7\textwidth, height=0.5\textwidth
    ]
    \addplot[mark=*, blue, thick] coordinates {
        (3, 1.4228) (4, 1.3159) (5, 1.2538) (6, 1.2128) (7, 1.1834) 
        (8, 1.1613) (9, 1.1441) (10, 1.1302) (11, 1.1188) (12, 1.1093)
    };
    \draw[dashed, red, thick] (axis cs:2.5,1) -- (axis cs:13,1);
    \node[anchor=west] at (axis cs:9, 1.1) {Limit $n \to \infty = 1$};
    \end{axis}
    \end{tikzpicture}
    \caption{Asymptotic convergence of the Steklov ratio.}
    \label{fig:asymptotic_ratio}
\end{figure}
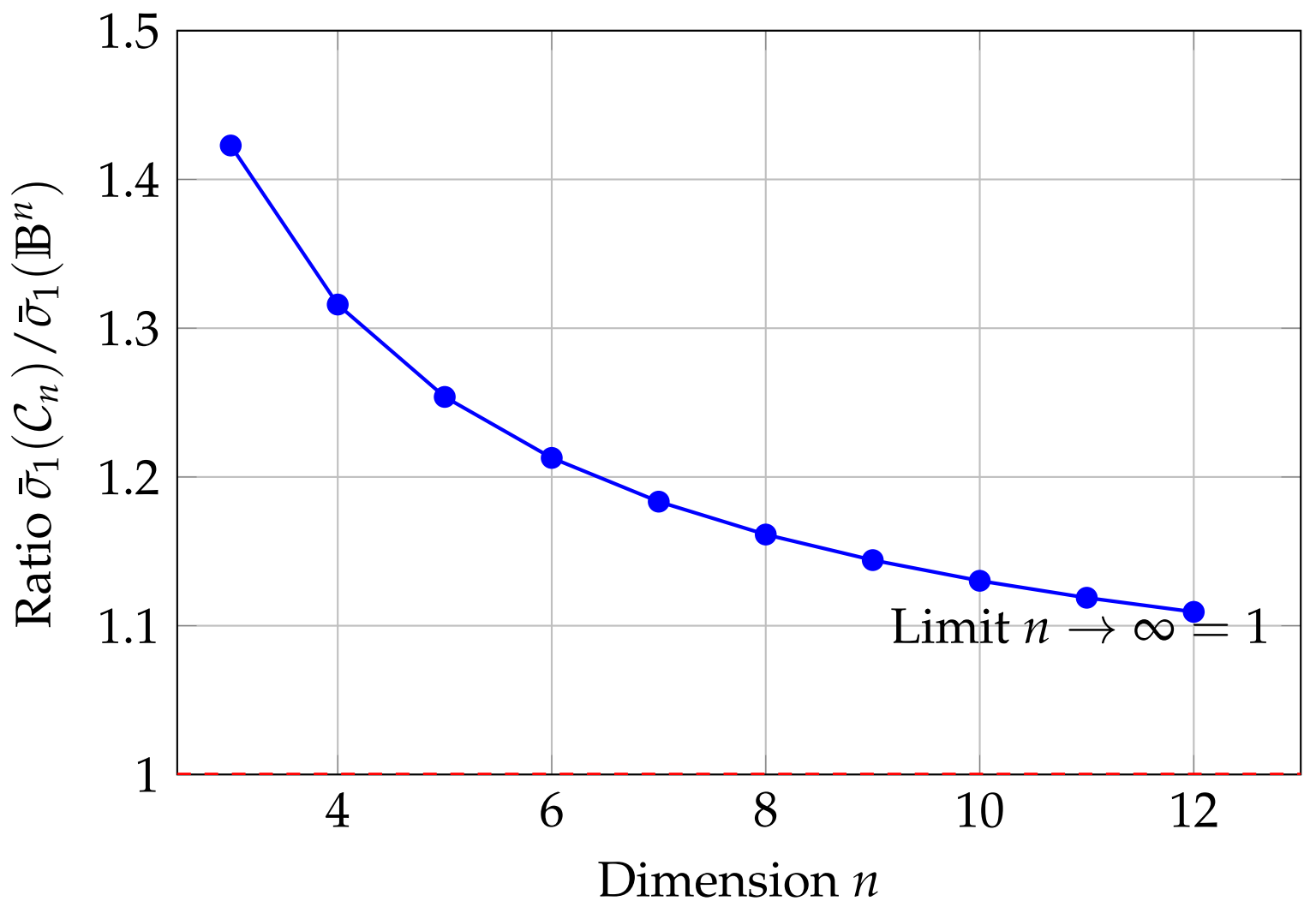

To quantify the spectral domination margin, we numerically evaluate the free boundary root $L_*$ across dimensions $n \in \{3, \dots, 12\}$. As shown in Table \ref{tab:catenoid_numerics} and Figure \ref{fig:asymptotic_ratio}, the geometric parameter $X = L_*^{2n-2}$ strictly exceeds the critical threshold $T_n$ for all evaluated dimensions. For $n=3$, the critical $3$-catenoid exceeds the normalized Steklov eigenvalue of the standard ball by over $42\%$.

\begin{table}[htpb]
\centering
\begin{tabular}{ccccc}
\hline
$n$ & $L_*$ & $X$ & Threshold $T_n$ & Ratio $\bar{\sigma}_1(\mathcal{C}_n)/\bar{\sigma}_1(\B^n)$ \\
\hline
3 & 1.603120 & 6.604876 & 2.000000 & 1.422834 \\
4 & 1.489373 & 10.914909 & 2.702414 & 1.315881 \\
5 & 1.416086 & 16.170270 & 3.414214 & 1.253813 \\
6 & 1.364344 & 22.348015 & 4.129813 & 1.212762 \\
7 & 1.325568 & 29.432223 & 4.847322 & 1.183426 \\
8 & 1.295261 & 37.411134 & 5.565925 & 1.161346 \\
9 & 1.270822 & 46.275741 & 6.285214 & 1.144095 \\
10 & 1.250632 & 56.018815 & 7.004959 & 1.130227 \\
11 & 1.233630 & 66.634457 & 7.725024 & 1.118827 \\
12 & 1.219085 & 78.117577 & 8.445322 & 1.109284 \\
\hline
\end{tabular}
\vspace{0.2cm}
\caption{Numerical evaluations of the critical $n$-catenoid geometry and the Steklov ratio.}
\label{tab:catenoid_numerics}
\end{table}

\begin{proposition}
 \label{cor:infinite_limit}
As the dimension $n \to \infty$, we have $X = L_*^{2n-2} \to +\infty$, and the normalized Steklov ratio $\bar{\sigma}_1(\mathcal{C}_n)/\bar{\sigma}_1(\B^n)$ converges to $1$.
\end{proposition}

\begin{proof}
Set $Y_n := L_*^{n-1} = \sqrt{X} > 1$. By inequality \eqref{eq:arccosh_bound} in the proof of Lemma \ref{lem:analytic_root_bound}, $Y_n$ satisfies
\[
    \sqrt{Y_n^2 - 1} \, \operatorname{arccosh}(Y_n) > n - 1 \longrightarrow +\infty \quad \text{as } n \to \infty.
\]
Since the function $y \mapsto \sqrt{y^2 - 1}\,\operatorname{arccosh}(y)$ is continuous and strictly increasing on $[1, +\infty)$, this forces $Y_n \to +\infty$ and hence $X = Y_n^2 \to +\infty$ as $n \to \infty$.

Consequently, $1 - \frac{1}{X} \to 1$. Since $\frac{n-1}{n} \in (0, 1)$, we have $1 - \frac{1}{X} \le \left(1 - \frac{1}{X}\right)^{\frac{n-1}{n}} \le 1$, so $\lim_{n \to \infty} \left(1 - \frac{1}{X}\right)^{\frac{n-1}{n}} = 1$. Since $\lim_{n \to \infty} 2^{\frac{2}{n}} = 1$, passing to the limit in \eqref{eq:ratio_X} gives:
\[
    \lim_{n \to \infty} \frac{\bar{\sigma}_1(\mathcal{C}_n)}{\bar{\sigma}_1(\B^n)} = \lim_{n \to \infty} \left( 2^{\frac{2}{n}} \left( 1 - \frac{1}{X} \right)^{\frac{n-1}{n}} \right) = 1.
\]
\end{proof}

\noindent\textbf{Geometric Context:} This analytic convergence reflects measure concentration in high dimensions. As $n \to \infty$, the volume of an $n$-sphere concentrates around its equator. Correspondingly, the critical $n$-catenoid flattens toward a double equatorial disk, and its normalized Steklov functional approaches that of the Euclidean unit ball $\B^n$.

\section{Conformal Maximization on the \texorpdfstring{$n$}{n}-Ball}\label{sec:maximization_of_Bn}

While the general bounds of Proposition \ref{prop:steklov_isop_bounds} ensure that $\bar{\sigma}_k$ is finite within any given conformal class, they do not guarantee the existence of a smooth maximizer in general. Even for Laplace eigenvalues, relatively few explicit examples of conformally maximal metrics are known in dimensions $n \ge 3$; notable results include \cite[Corollary 7]{LY1} and \cite{Berger73, ESvol} for $(\mathbb{S}^n, g_{\text{std}})$, \cite{Kim} for the second eigenvalue on $\mathbb{S}^n$, \cite{EG_new} for the second eigenvalue on $(\mathbb{RP}^n, g_{\text{std}})$, and \cite{ESTore} for flat tori. 

Brock \cite{Brock} proved sharp bounds for Euclidean domains, and Fraser and Schoen \cite[Theorem 6.2]{FS11} established conformal bounds via the reduced conformal volume. For completeness, we record here a short direct proof characterizing the unique conformal maximizer of the scale-invariant functional \eqref{eq:KM_normalization} on the topological $n$-ball, modeled by the upper hemisphere $\Omega = \mathbb{S}^n_+ := \mathbb{S}^n \cap \{x_{n+1} > 0\}$. Let $\pi_S \colon \mathbb{S}^n_+ \to \mathbb{B}^n$ be a conformal diffeomorphism onto the Euclidean unit ball, and define the pullback metric $g_0 := \pi_S^*(g^{\mathbb{B}^n}_{\mathrm{std}}) \in [g^{\mathbb{S}^n_+}_{\mathrm{std}}]$.

\begin{theorem}\label{thm:ball_supremum}
Let $n \geq 3$. For all $g \in [g_0]$ and $\rho \in C^\infty_{>0}(\partial \mathbb{S}^n_+)$:
\begin{equation} \label{eq:St_W}
\bar{\sigma}_1(\mathbb{S}^n_+, g, \rho) \leq \bar{\sigma}_1(\mathbb{S}^n_+, g_0, 1),
\end{equation}
with equality if and only if $\rho$ is constant and $(\mathbb{S}^n_+, g)$ is isometric to $(\mathbb{S}^n_+, g_0)$ up to scaling.
\end{theorem}

\begin{proof}
Write $M := \mathbb{S}^n_+$. Composing $\pi_S$ with a suitable conformal automorphism of $\mathbb{B}^n$, we choose a balanced conformal diffeomorphism $F = (F_1, \dots, F_n) \colon (M, g) \to (\mathbb{B}^n, g_{\mathrm{std}})$ such that
\[
    \int_{\partial M} F_i \rho \, d\mathcal{H}_g = 0 \qquad (i = 1, \dots, n).
\]
Each coordinate function $F_i$ is thus an admissible test function for $\sigma_1(M, g, \rho)$. Since $\sum_{i=1}^n F_i^2 = 1$ on $\partial M$, summing the Rayleigh inequalities over $i = 1, \dots, n$ and applying H\"older's inequality yields
\[
    \sigma_1(M, g, \rho) \int_{\partial M} \rho \, d\mathcal{H}_g \le \int_M |dF|_g^2 \, dv_g \le \vol(M, g)^{\frac{n-2}{n}} \left( \int_M |dF|_g^n \, dv_g \right)^{\frac{2}{n}}.
\]
Conformality and surjectivity of $F \colon (M, g) \to (\mathbb{B}^n, g_{\mathrm{std}})$ imply
\[
    \int_M |dF|_g^n \, dv_g = n^{\frac{n}{2}} \vol(\mathbb{B}^n, g_{\mathrm{std}}),
\]
and therefore $\bar{\sigma}_1(M, g, \rho) \le n \vol(\mathbb{B}^n, g_{\mathrm{std}})^{\frac{2}{n}}$. Since $\sigma_1(\mathbb{B}^n, g_{\mathrm{std}}, 1) = 1$ \cite{GP} and $\mathcal{H}_{g_{\mathrm{std}}}(\partial \mathbb{B}^n) = n \vol(\mathbb{B}^n, g_{\mathrm{std}})$, the right-hand side equals $\bar{\sigma}_1(\mathbb{B}^n, g_{\mathrm{std}}, 1) = \bar{\sigma}_1(\mathbb{S}^n_+, g_0, 1)$, proving \eqref{eq:St_W}.

If equality holds in \eqref{eq:St_W}, then equality in H\"older's inequality forces $|dF|_g$ to be constant on $M$. By conformality, $F^* g_{\mathrm{std}} = c \, g$ for some constant $c > 0$, so $(M, g)$ is homothetic to $(\mathbb{B}^n, g_{\mathrm{std}})$ and hence isometric to $(\mathbb{S}^n_+, g_0)$ up to scaling. Moreover, equality in the Rayleigh inequalities implies that $\partial_{\nu_g} F = \sigma_1(M, g, \rho) \rho F$ on $\partial M$. Since $F$ is a homothety onto the unit ball $\mathbb{B}^n$, we also have $\partial_{\nu_g} F = \sqrt{c} \, F$ on $\partial M$, so $\rho$ is constant. The converse follows immediately from the scale invariance of $\bar{\sigma}_1$ in $g$ and $\rho$.
\end{proof}

\begin{remark}
The round hemisphere $(\mathbb{S}^n_+, g_{\mathrm{std}}, 1)$ is not a maximizer for $n \ge 3$. Its first nonzero Steklov eigenvalue is $\sigma_1 = \frac{2}{\sqrt{\pi}} \frac{\Gamma((n + 1)/2)}{\Gamma(n/2)}$, with eigenfunctions expressed in terms of associated Legendre functions $P^\mu_\nu$ (see \cite{Bri3, Nist}).
\end{remark}

\begin{proposition}\label{prop:ball_supremum}
The conformal supremum of the normalized first Steklov eigenvalue on $\mathbb{B}^n$ is uniquely achieved (up to conformal automorphisms of $\mathbb{B}^n$ and scaling) by the standard flat metric $g_{\mathrm{std}}$ with constant boundary weight:
\[
\sup_{\substack{g \in [g_0] \\ \rho > 0}}\bar{\sigma}_1(\mathbb{S}^n_+, g, \rho) = \bar{\sigma}_1(\mathbb{B}^n, g_{\mathrm{std}}, 1) = n^{\frac{n-2}{n}} \vol(\mathbb{S}^{n-1})^{\frac{2}{n}}.
\]
\end{proposition}

\begin{proof}
This follows directly from Theorem \ref{thm:ball_supremum} and the evaluation of \eqref{eq:KM_normalization} on $(\mathbb{B}^n, g_{\mathrm{std}}, 1)$.
\end{proof}

\begin{remark}
Since any spherical cap on $\mathbb{S}^n$ is conformally diffeomorphic to $\mathbb{S}^n_+$, the bound in Proposition \ref{prop:ball_supremum} holds for any spherical cap.
\end{remark}



\begin{thebibliography}{KKMS24}

\bibitem[AM12]{AM} W.~Arendt and R.~Mazzeo. {\em Friedlander's eigenvalue inequalities and the Dirichlet-to-Neumann semigroup}. Commun. Pure Appl. Anal., {\bf 11}:6 (2012), 2275--2285.


\bibitem[Ber73]{Berger73} M.~Berger. {\em Sur les premières valeurs propres des variétés riemanniennes}. Compositio Mathematica, {\bf 26}:2 (1973), 129--149.

\bibitem[Ble83]{Ble} D.~D.~Bleecker. {\em The spectrum of a Riemannian manifold with a unit Killing vector field}. Transactions of the American Mathematical Society, {\bf 275} (1983), 409--416.

\bibitem[Bri22]{Bri1} J.~Brisson. {\em Tubular excision and Steklov eigenvalues}. The Journal of Geometric Analysis, {\bf 32}:166 (2022).

\bibitem[Bri23]{Bri3} J.~Brisson. {\em Excisions tubulaires et valeurs propres de Steklov de boules géodésiques}. Thèse de doctorat, Université Laval, 2023.

\bibitem[Bri24]{Bri2} J.~Brisson. {\em Multiple tubular excisions and large Steklov eigenvalues}. Annals of Global Analysis and Geometry, {\bf 65}:2 (2024), 18.

\bibitem[Bro01]{Brock} F.~Brock. {\em An isoperimetric inequality for eigenvalues of the Stekloff problem}. ZAMM - Journal of Applied Mathematics and Mechanics, {\bf 81}:1 (2001), 69--71.

\bibitem[Cla83]{Clarke} F.~H.~Clarke. {\em Optimization and Nonsmooth Analysis}. Wiley-Interscience Series in Discrete Mathematics and Optimization, John Wiley \& Sons, New York, 1983.

\bibitem[Cod12]{Cod} E.~A.~Coddington. {\em An Introduction to Ordinary Differential Equations}. Dover Publications, 2012.

\bibitem[CD94]{CD} B.~Colbois and J.~Dodziuk. {\em Riemannian metrics with large $\lambda_1$}. Proceedings of the American Mathematical Society, {\bf 122}:3 (1994), 905--906.

\bibitem[CEG19]{CEG} B.~Colbois, A.~El Soufi, and A.~Girouard. {\em Compact manifolds with fixed boundary and large Steklov eigenvalues}. Proceedings of the American Mathematical Society, {\bf 147}:9 (2019), 3813--3827.

\bibitem[CG23]{CG} B.~Colbois and A.~Girouard. {\em Metric upper bounds for Steklov and Laplace eigenvalues}. Annales de l'Institut Fourier, {\bf 73}:2 (2023), 811--839.

\bibitem[CM08]{ColMae} B.~Colbois and D.~Maerten. {\em Eigenvalue estimates for the Neumann problem of a bounded domain}. The Journal of Geometric Analysis, {\bf 18}:4 (2008), 1022--1032.

\bibitem[EG25]{EG_new} M.~Eddaoudi and A.~Girouard. {\em Upper bounds for the second nonzero eigenvalue of the Laplacian via folding and conformal volume}. Preprint arXiv:2501.08761, 2025.

\bibitem[Eke74]{Ekeland} I.~Ekeland. {\em On the variational principle}. Journal of Mathematical Analysis and Applications, {\bf 47}:2 (1974), 324--353.


\bibitem[ESI84]{ESvol} A.~El Soufi and S.~Ilias. {\em Le volume conforme et ses applications d'après Li et Yau}. Séminaire de théorie spectrale et géométrie, Université de Grenoble I (1984), 1--15.

\bibitem[ESI03]{ESTore} A.~El Soufi and S.~Ilias. {\em Extremal metrics for the first eigenvalue of the Laplacian in a conformal class}. Proceedings of the American Mathematical Society, {\bf 131}:5 (2003), 1611--1618.

\bibitem[FS11]{FS11} A.~Fraser and R.~Schoen. {\em The first Steklov eigenvalue, conformal geometry, and minimal surfaces}. Advances in Mathematics, {\bf 226}:5 (2011), 4011--4030.

\bibitem[FS12]{FS12} A.~Fraser and R.~Schoen. {\em Minimal surfaces and eigenvalue problems}. Contemporary Mathematics, {\bf 599} (2012), 105--121.

\bibitem[FS16]{FS16} A.~Fraser and R.~Schoen. {\em Sharp eigenvalue bounds and minimal surfaces in the ball}. Inventiones Mathematicae, {\bf 203}:3 (2016), 823--890.

\bibitem[FS19]{FS19} A.~Fraser and R.~Schoen. {\em Shape optimization for the Steklov problem in higher dimensions}. Advances in Mathematics, {\bf 348} (2019), 146--162.

\bibitem[FGM17]{FGM17} M. Freidin, M. Gulian, and P. McGrath. {\em Free boundary minimal surfaces in the unit ball with low cohomogeneity}. Proceedings of the American Mathematical Society, {\bf 145}:4 (2017), 1671--1683.



\bibitem[GKL21]{GKL} A. Girouard, M. Karpukhin, and J. Lagac\'e. {\em Continuity of eigenvalues and shape optimization for Laplace and Steklov problems}. Geometric and Functional Analysis, {\bf 31}:3 (2021), 513--561.

\bibitem[GL21]{GL} A.~Girouard and J. Lagac\'e. {\em Large Steklov eigenvalues via homogenisation on manifolds}. Inventiones Mathematicae, {\bf 226}:3 (2021), 1011--1056.

\bibitem[GP10]{GP} A.~Girouard and I.~Polterovich. {\em On the Hersch-Payne-Schiffer inequalities for Steklov eigenvalues}. Functional Analysis and its Applications, {\bf 44}:2 (2010), 106--117. 

\bibitem[GP17]{GP17} A.~Girouard and I.~Polterovich. {\em Spectral geometry of the Steklov problem (Survey article)}. Journal of Spectral Theory, {\bf 7}:2 (2017), 321--359.

\bibitem[GP24]{GP2} A.~Girouard and P.~Polymerakis. {\em Large Steklov eigenvalues under volume constraints}. The Journal of Geometric Analysis, {\bf 34}:10 (2024), 316.

\bibitem[GNY04]{GNY} A. Grigor’yan, Y. Netrusov, and S.-T. Yau. {\em Eigenvalues of elliptic operators and geometric applications}. Surveys in Differential Geometry, {\bf 9}:1 (2004), 147--217.

\bibitem[Has11]{Hass} A.~Hassannezhad. {\em Conformal upper bounds for the eigenvalues of the Laplacian and Steklov problem}. Journal of Functional Analysis, {\bf 261}:12 (2011), 3419--3436.

\bibitem[KKMS24]{KKMS} M. Karpukhin, R. Kusner, P. McGrath, and D. Stern. {\em Embedded minimal surfaces in $\mathbb{S}^3$ and $\mathbb{B}^3$ via equivariant eigenvalue optimization}. Preprint arXiv:2402.13121, 2024. 

\bibitem[KM22]{KM} M. Karpukhin and A. Métras. {\em Laplace and Steklov extremal metrics via $n$-harmonic maps}. The Journal of Geometric Analysis, {\bf 32}:5 (2022), 154.

\bibitem[KNPP20]{KNPP} M. Karpukhin, N. Nadirashvili, A. V. Penskoi, and I. Polterovich. {\em Conformally maximal metrics for Laplace eigenvalues on surfaces}. Surveys in Differential Geometry, {\bf 24}:1 (2020), 205--256.

\bibitem[KN25]{KN} M.~Karpukhin and A.~D.~Noutchegueme. {\em Optimisation of Steklov transmission eigenvalues and minimal surfaces}. Preprint arXiv:2509.24055, 2025.

\bibitem[Kim22]{Kim} H.~Kim. {\em Maximization of the second Laplacian eigenvalue on the sphere}. Proceedings of the American Mathematical Society, {\bf 150}:8 (2022), 3449--3460.

\bibitem[Kok14]{Kokarev} G.~Kokarev. {\em Variational aspects of Laplace eigenvalues on Riemannian surfaces}. Advances in Mathematics, {\bf 258} (2014), 191--239.

\bibitem[Kor93]{Kor} N.~Korevaar. {\em Upper bounds for eigenvalues of conformal metrics}. Journal of Differential Geometry, {\bf 37} (1993), 73--93.

\bibitem[LMP23]{LMP} M.~Levitin, D.~Mangoubi, and I.~Polterovich. {\em Topics in Spectral Geometry}. Graduate Studies in Mathematics, Vol. {\bf 237}, American Mathematical Society, 2023.

\bibitem[LY80]{LY1} P.~Li and S.-T.~Yau. {\em Estimates of eigenvalues of a compact Riemannian manifold}. Proceedings of Symposia in Pure Mathematics, {\bf 36} (1980), 205--239.

\bibitem[LY82]{LY2} P.~Li and S.-T.~Yau. {\em A new conformal invariant and its applications to the Willmore conjecture and the first eigenvalue of compact surfaces}. Inventiones Mathematicae, {\bf 69}:2 (1982), 269--291.

\bibitem[Lin06]{Lin} P.~Lindqvist. {\em Notes on the $p$-Laplace equation}. Report Univ. Jyväskylä Dep. Math. Stat., {\bf 102}, 2006.

\bibitem[Mat05]{Mat2} A.-M.~Matei. {\em Boundedness of the first eigenvalue of the $p$-Laplacian}. Proceedings of the American Mathematical Society, {\bf 133}:10 (2005), 3013--3019.

\bibitem[Mat13]{Mat1} A.-M.~Matei. {\em Conformal bounds for the first eigenvalue of the $p$-Laplacian}. Nonlinear Analysis, {\bf 80} (2013), 11--21.

\bibitem[Mor26]{Mor} E.~Morozov. {\em Rotationally symmetric critical metrics for Laplace eigenvalues on tori in a conformal class}. Preprint arXiv:2601.19196, 2026.

\bibitem[Nis10]{Nist} F.~W.~J.~Olver et al. {\em NIST Handbook of Mathematical Functions}. Cambridge University Press, 2010.

\bibitem[Pet14]{Pet14} R.~Petrides. {\em Existence and regularity of maximal metrics for the first Laplace eigenvalue on surfaces}. Geometric and Functional Analysis, {\bf 24}:4 (2014), 1336--1376.

\bibitem[Pet19]{Pet19} R.~Petrides. {\em Maximizing Steklov eigenvalues on surfaces}. Journal of Differential Geometry, {\bf 113}:1 (2019), 95--188.

\bibitem[Pet24]{Pet24} R.~Petrides. {\em Geometric spectral optimization on surfaces}. Preprint arXiv:2410.13347, 2024.

\bibitem[Sch07]{Sch} W.~Schirotzek. {\em Nonsmooth Analysis}. Universitext, Springer, 2007.

\bibitem[SSTZ21]{SSTZ} G.~Smith, A.~Stern, H.~Tran, and D.~Zhou. {\em On the Morse index of higher-dimensional free boundary minimal catenoids}. Calculus of Variations and Partial Differential Equations, {\bf 60}:5 (2021), 177.

\bibitem[Vin25]{Vin1} D.~Vinokurov. {\em Conformal optimization of eigenvalues on surfaces with symmetries}. Preprint arXiv:2502.03756, 2025.

\bibitem[Vin26]{Vin2} D.~Vinokurov. {\em Eigenvalue optimization in higher dimensions and $p$-harmonic maps}. Preprint arXiv:2601.17896, 2026.

\bibitem[Zet05]{Zet} A.~Zettl. {\em Sturm-Liouville Theory}. Mathematical Surveys and Monographs, Vol. {\bf 121}, American Mathematical Society, 2005.

\end{thebibliography}
\end{document}